\documentclass[leqno,11pt,oneside]{amsart}

\usepackage{bm,cancel}
\usepackage{amssymb,amsfonts}
\usepackage{amsmath,latexsym,amsthm}
\usepackage{pdfsync} % this screws up the marginpar command
\usepackage{xcolor,graphicx}
\usepackage{mathrsfs}
\usepackage{ulem}
\usepackage[unicode,breaklinks=true,colorlinks=true]{hyperref}

\usepackage[top=1in, bottom=1in, left=1.25in, right=1.25in, marginparwidth=1.1in, marginparsep=0.05in]{geometry}

\vspace{10pt}

\newtheorem{thm}{Theorem}[section]

\newtheorem{lem}[thm]{Lemma}
\newtheorem{prop}[thm]{Proposition}
\theoremstyle{definition}
\newtheorem{defn}[thm]{Definition}

\newtheorem{rem}[thm]{Remark}
\numberwithin{equation}{section}
\newcommand{\eqdefa}{\,\overset{\hbox{\scalebox{0.6}{def}}}{\scalebox{1.3}[1]{=}}\,}

\newcommand{\R}{\mathbb R}
\newcommand\C{\mathcal{C}}

\newcommand{\CF}{\mathbb{F}}
\newcommand{\FF}{\mathbb{F}}
\newcommand{\GG}{\mathbb{G}}
\newcommand{\cH}{\mathcal{H}}
\newcommand{\cW}{\mathcal{W}}

\newcommand{\cA}{\mathcal{A}}

\newcommand{\dy}{\,{\rm d}y}

\newcommand{\dx}{\,{\rm d}x}
\newcommand{\dSx}{\,{\rm d}S_x}

\newcommand{\dt}{\,{\rm d}t}
\newcommand{\ds}{\,{\rm d}s}

\newcommand{\andf}{\quad\hbox{and}\quad}

\renewcommand{\div}{\mathop{\rm div}\nolimits}
\newcommand{\curl}{\mathop{\rm curl}\nolimits}
\renewcommand{\leq}{\leqslant}
\renewcommand{\geq}{\geqslant}

\def\XXint#1#2#3{{\setbox0=\hbox{$#1{#2#3}{\int}$}\vcenter{\hbox{$#2#3$}}\kern-.5\wd0}}

\newcommand{\lec}{\lesssim}
\newcommand{\norm}[1]{\left\|#1\right\|}
\newcommand{\tsum}{{\textstyle \sum}}
\newcommand{\bke}[1]{\left ( #1 \right )}
\newcommand{\bkt}[1]{\left [ #1 \right ]}
\newcommand{\bket}[1]{\left \{ #1 \right \}}

\newcommand{\pd}{\partial}
\newcommand{\DD}{\mathbb{D}}
\newcommand{\ep}{\varepsilon}

\newcommand{\supp}{\operatorname{supp}}

\begin{document}
\title{The local regularity theory for the Stokes and Navier--Stokes equations near the curved boundary}%
\author[H. Chen]{Hui Chen}%
\address[H. Chen]
{School of Science, Zhejiang University of Science and Technology, Hangzhou, 310023, China }
\email{chenhui@zust.edu.cn}
\author[S. Liang]{Su Liang}
\address[S. Liang]
{Department of Mathematics, University of British Columbia, Vancouver, BC V6T1Z2, Canada }
\email{liangsu96@gmail.com}
\author[T.-P. Tsai]{Tai-Peng Tsai}
\address[T.-P. Tsai]
{Department of Mathematics, University of British Columbia, Vancouver, BC V6T1Z2, Canada }
\email{ttsai@math.ubc.ca}
\date{}
%
%%--------------------------------------------------------------------------------------
\begin{abstract}
In this paper, we study local regularity of the solutions to the Stokes equations near a curved boundary under no-slip or Navier boundary conditions. We extend previous  boundary estimates near a flat boundary to that near a curved boundary, under very low starting regularity assumptions. Compared with the flat case, the proof for the curved case is more complicated and we adapt new techniques such as the  ``normal form" after the mollification, recovering vertical derivative estimates from horizontal derivative estimates, and transferring temporal derivatives to spatial derivatives, to deal with the higher order perturbation terms generated by boundary straightening. 
As an application, we propose a new definition of boundary regular points for the incompressible Navier--Stokes equations  that guarantees higher spatial regularity.
\end{abstract}
\maketitle


\vskip 0.2cm


\tableofcontents

\section{Introduction}

The main goal of this paper is the improvement of regularity of solutions to the Stokes equations near a curved boundary under no-slip or Navier boundary conditions. Based on these results, we also propose an alternative definition of boundary regular points.

Let $\Omega \subset\R^d$, $d\geq 2$, be a domain with $\pd \Omega\in C^{1,1}$. 
For simplicity we assume the origin is located on the boundary $\pd\Omega$.  Denote  $B_R=B_R(0)$, $B^+_R=B_R\cap \Omega$ (could be disconnected) and $Q^+_R=B^+_R\times(-R^2,0)$. Consider the non-stationary Stokes equations in $Q^+_1$
\begin{equation}\tag{S}\label{Stokes-Eqn}
\left\{\begin{aligned}
\,\partial_{t}\bm{u}-\Delta\bm{u}+\nabla p&=\bm{f}+\div \CF,\\[6pt]
\div \bm{u}&=0,
\end{aligned}\right.
\end{equation}
where $\bm{u}=\left(u_{1},\cdots,u_{d}\right)$ is the velocity field,  $p$ is the pressure, $\bm{f}=\left(f_{1},\cdots,f_{d}\right)$ is the body force and the matrix $\CF=\left(F_{ij}\right)_{1\leq i,j\leq d}$ is the surface force.  Replacing $\CF$ with $-\bm{u}\otimes \bm{u}$ we recover the Navier--Stokes equations (and assume $\bm{f}=0$)
\begin{equation}\tag{NS}\label{NS}
\left\{
\begin{aligned}
\,\partial_{t}\bm{u}+\bm{u}\cdot\nabla\bm{u}-\Delta\bm{u}+\nabla p&=0,\\[6pt]
 \div \bm{u}&=0.
\end{aligned}
\right.
\end{equation}
These equations describe the flow of incompressible viscous fluids. We refer to the books \cite{Lem16,Tsa18} for more background of the Navier--Stokes equations \eqref{NS}.

In contrast to equations of parabolic type, weak solutions of the Stokes equations \eqref{Stokes-Eqn} may not be smooth near boundary even for zero external forces. Its exact behavior varies according to different boundary conditions on the  boundary portion $\Sigma=\Gamma\times(-1,0)$ with  $\Gamma=B_1\cap\pd\Omega$. In this paper, we are concerned with two types of boundary conditions (BC): no-slip boundary condition
\begin{align}\label{0-BC0}
\bm{u}  = 0 \quad \hbox{on }\Sigma,
\end{align}
and Navier boundary condition
\begin{equation}\label{NavierBC0}
 [\left(2\DD\bm{u}+\CF\right)\bm{n}]_{tan}+\alpha\bm{u}_{tan}=0\andf\bm{u}\cdot\bm{n}=0 \quad \hbox{on }\Sigma,
\end{equation}
where $\bm{n}$ is the unit outward normal vector to the boundary $\partial \Omega$, $\DD\bm{u}=\frac{1}{2}\left(\nabla\bm{u}+(\nabla\bm{u})^{T}\right)$ is the deformation tensor,  $\bm{v}_{tan}=\bm{v}-\left(\bm{v}\cdot\bm{n}\right)\bm{n}$ is the tangential component of a vector on the boundary, and  
$\alpha \geq 0$ is the \emph{friction coefficient}.  For $\alpha=2\kappa_\tau$ and $\kappa_\tau$ the principal curvature tensor of $\pd \Omega$, the Navier boundary condition \eqref{NavierBC0} is reduced to the following boundary condition (see e.g. \cite[Proposition A.1]{DK26-2})
\begin{equation*}
\bm{\omega}\bm{n}=0\andf \bm{u}\cdot\bm{n}=0\quad \hbox{on }\Sigma,
\end{equation*}
where $\bm{\omega}$ is  the skew-symmetric part of the gradient $\nabla \bm{u}$, i.e., the vorticity matrix. For dimension $d=3$, the condition $\bm{\omega}\bm{n}=0$ is equivalent to $(\curl \bm{u}) \times \bm{n}=0$.

We recall recent developments of local regularity theory for Stokes equations \eqref{Stokes-Eqn} near a flat boundary, i.e., when $\Gamma$ is a plane. We denote $f\in L^{q,r}(\Omega\times I)$ if $f\in L^{r}\left(I;L^{q}(\Omega)\right)$ for $I\subset \R$.  

Under no-slip boundary condition \eqref{0-BC0}, Seregin-Shilkin \cite[Theorem 2.5]{SS14} showed spatial smoothing ($\nabla^2 \bm{u},\nabla p\in L^{q,r}(Q_{1/2}^+)$, $1<q,r<\infty$)  if one assumes $\bm{u}, \nabla \bm{u}, p \in L^{q,r}(Q_1^+)$. (Also see the {\it a priori} bound of Seregin \cite[Proposition 1]{Ser00}.)
However, without any assumption on the pressure, the smoothing in spatial variables fails due to non-local effect of the pressure. The first counter example was constructed by K. Kang \cite{Kan04a} using a boundary flux supported away from $\Sigma$, and was refined in \cite{KLLT22}.
Seregin-\v Sver\'ak \cite{SS10} found another counter example in the form of shear flows. Recently, Chang-Kang \cite[Theorem 1.1]{CK23} proved that for any $1<q<\infty$, a bounded very weak solution (depending on $q$) can be constructed whose derivatives are unbounded in $L^{q}$,
\begin{equation}\label{CK}
\|\nabla \bm{u}\|_{L^q (Q_{1/2}^{+})}=\infty, \quad\|\bm{u}\|_{L^{\infty}\left(Q_1^{+}\right)}<\infty.
\end{equation}
Kang-Min \cite{KangMin} constructed examples of bounded $\bm{u}$ and $\nabla \bm{u}$ with $\nabla^2 \bm{u}$ unbounded in $L^q(Q_{1}^{+})$. Also see Chang-Kang \cite{CK26} for construction of unbounded velocity, although the velocity and its gradient are locally square integrable.

Under Navier boundary condition \eqref{NavierBC0}, we  proved in \cite[Theorem 1.1]{CLT24} that $\nabla \bm{u}\in L^{q,r}(Q^{+}_{1/2})$ assuming only $\bm{u}\in L^{q,r}(Q^{+}_1)$. It is a very important distinction in comparison to \eqref{CK} for no-slip boundary condition. 
Moreover, we also proved in \cite[Theorem 1.2]{CLT24} that $\nabla^3\bm{u},\nabla^2p\in L^{q,r}(Q^+_{1/2})$ if  $\bm{u},p\in L^{q,r}(Q^{+}_1)$. We refer to \cite{CK26,CLT24,CLT25,DKP22,DK26, DK26-2} and the references therein for more literature under various boundary conditions.

In this paper, we are interested in the local regularity of the solutions near a \emph{curved boundary}. Throughout this paper, let
\begin{equation}\label{q_*}
1<q,r<\infty,\quad 
q_{*}=\left\{\begin{array}{lll}
\max\{1,\frac{dq}{q+d}\},&\text{if }\, q\neq \frac{d}{d-1}\\
1+,&\text{if }\,q=\frac{d}{d-1}
\end{array}\right..
\end{equation}

Our first set of results, Theorems \ref{thm_0BC}--\ref{thm2_navierBC}, is for the Stokes equations \eqref{Stokes-Eqn}. They are optimal in the sense that we do not expect higher derivative estimates.

\begin{thm}[Derivative estimates for no-slip BC]\label{thm_0BC}
\sl{
Assume \eqref{q_*}. Suppose that the  boundary portion $\Gamma \in C^{1,1}$, $\bm{u}, p, \CF\in L^{q,r}(Q^+_1)$, $\bm{f}\in L^{q_*,r}(Q^+_1)$, and $(\bm{u},p)$ is a very weak solution pair (see Definition \ref{def:vwsp-0BC}) to the Stokes equations \eqref{Stokes-Eqn} with no-slip boundary condition \eqref{0-BC0}, then $(\bm{u},p)$ is a weak solution pair (see Definition \ref{def:wsp-0BC}) satisfying $\nabla\bm{u}\in L^{q,r}(Q^+_{1/2})$, and 
\begin{equation}\label{gra1}
\| \nabla \bm{u}\|_{L^{q,r}\left(Q^+_{1/2}\right)}\lec \|\bm{u}\|_{L^{q,r}(Q^+_1)}+\|p\|_{L^{q,r}(Q^+_1)}+\|\bm{f}\|_{L^{q_*,r}(Q^+_1)}+\|\CF\|_{L^{q,r}(Q^+_1)}.
\end{equation}
If in addition $ \bm{f}+\div \CF \in L^{q,r}(Q^+_1)$, we have $\pd_{t}\bm{u},\nabla^2\bm{u},\nabla p\in L^{q,r}(Q^+_{1/2})$, and
\begin{equation}\label{gra0}
 \| \pd_{t} \bm{u},\nabla^2\bm{u},\nabla p\|_{L^{q,r}\left(Q^+_{1/2}\right)}\lec \|\bm{u}\|_{L^{q,r}(Q^+_1)}+\|p\|_{L^{q,r}(Q^+_1)}+\|\bm{f}+\div \CF\|_{L^{q,r}(Q^+_1)}.
\end{equation}
}
\end{thm}

\emph{Comments on Theorem \ref{thm_0BC}}:
\begin{enumerate}
\renewcommand{\theenumi}{\roman{enumi}}
\item
That any strong solution pair (with $\nabla^2\bm{u}, \pd_t \bm{u}, p, \nabla p\in L^{q,r}(Q^+_1)$) near a curved boundary satisfies estimate \eqref{gra0} was shown by Solonnikov \cite[Theorem 1.2]{Sol04}. Later, Vialov-Shilkin \cite[Theorem 2.3]{VS13} proved that any weak solution pair near a $C^{2,1}$ curved boundary 
is a strong solution.
The contribution of Theorem \ref{thm_0BC} is to show that any \emph{very weak solution pair} is a strong solution near a $C^{1,1}$ curved boundary. 

\item Assuming that the higher derivatives already exist in $L^{q,r}$,
the spatial exponents on the LHS of \eqref{gra0} can be improved assuming better regularity of $\bm{f}$ and $\CF$, see Seregin \cite[Proposition 2]{Ser00}, \cite{Ser09} and Chang-Kang \cite[Theorem 1.5]{CK23} for the case of weak solution pair near a flat boundary, and \cite[Theorem 2.4]{VS13} for that near a  curved boundary.

\item The presence of $\norm{p}_{L^{q,r}}$ on the RHS of \eqref{gra1} and \eqref{gra0} is important. It excludes counter-examples of \cite{CK23,CLT25,Kan04a,KLLT22,KangMin,SS10}.

\item Estimate \eqref{gra1} can indeed be derived from the corresponding local estimate for the heat equation, applied componentwise. The reason for including a proof here is that our argument works directly at the level of very weak solutions, with rather low initial regularity. Although boundary Caccioppoli-type estimates for the heat equation are classical in spirit, we were not able to locate a reference treating this estimate in the framework of very weak solutions with the low regularity considered here. As a result, we do not need  $\norm{\nabla \bm{u}}_{L^{q,r}}$ on the RHS of \eqref{gra0}, unlike \cite{CK23,Ser00,SS14,VS13}.

\item As far as we know, estimate \eqref{gra0} does not seem to follow directly from the standard local estimates for the heat equation and the Laplace equation alone; see Remark \ref{rem:seregin-prop1}.

\item Unlike \eqref{second-est} in Theorem \ref{thm2_navierBC}, estimate  \eqref{gra0} does not need $\Gamma \in C^{2,1}$, and does not need boundary condition for $\FF$.

\item 
We do not expect third derivative estimates  in Theorem  \ref{thm_0BC}, by similar reasons of \cite[Remark 6.3]{CLT25} for Navier BC case.

\end{enumerate}

\begin{thm}[Boundary gradient estimate for Navier BC]\label{thm_navierBC}
\sl{
Assume \eqref{q_*}. Suppose that the  boundary portion $\Gamma \in C^{1,1}$, $\bm{u}, \CF\in L^{q,r}(Q^+_1)$, $\bm{f}\in L^{q_*,r}(Q^+_1)$,   $\nabla\bm{u}\in L^{1}(Q^+_1)$, and $\bm{u}$ is a weak solution (see Definition \ref{def:ws}) to the Stokes equations \eqref{Stokes-Eqn} with Navier boundary condition \eqref{NavierBC0}, then we have $\nabla \bm{u}\in L^{q,r}(Q^+_{1/2})$ and
\begin{equation}\label{main1}
\| \nabla\bm{u}\|_{L^{q,r}\left(Q^+_{1/2}\right)}\lec \|\bm{u}\|_{L^{q,r}\left(Q^+_1\right)}+\|\bm{f}\|_{L^{q_{*},r}\left(Q^+_1\right)}+\|\CF\|_{L^{q,r}\left(Q^+_1\right)}.
\end{equation}
}
\end{thm}
\emph{Comments on Theorem \ref{thm_navierBC}}:
\begin{enumerate}
\item[(i)] Theorem \ref{thm_navierBC} extends \cite[Theorem 1.1]{CLT24} from a flat boundary to a curved boundary, and from 3D to any dimension $d\geq 2$.

\item[(ii)] However, \cite[Theorem 1.1]{CLT24} shows the gradient estimate \eqref{main1} for very weak solutions near a flat boundary. For technical reasons, we cannot prove the corresponding result near a curved boundary. Here we assume $\nabla \bm{u}\in L^{1}(Q^+_1)$ as our starting regularity, so Theorem \ref{thm_navierBC} is only for weak solutions. On the other hand, the norm $\| \nabla\bm{u}\|_{L^{1}\left(Q^+_{1}\right)}$ is not needed on the right side of \eqref{main1}.

\item[(iii)]

Unlike \eqref{gra1} for no-slip BC case, we do not need 
$\norm{p}_{L^{q,r}}$ on the RHS of   \eqref{main1}. This is already known in \cite{CLT24}.

\item[(iv)]
The theorem remains true if we assume $\alpha$ is a nonnegative function on $\pd \Omega$ with sufficient smoothness.

\end{enumerate}
\begin{thm}[Higher-order derivative estimates for Navier BC]\label{thm2_navierBC}
\sl{
Assume \eqref{q_*}. Suppose that the  boundary portion $\Gamma \in C^{2,1}$, $\bm{u}, p, \bm{f}\in L^{q,r}(Q^+_1)$, $\CF=0$, $\nabla\bm{u}\in L^{1}(Q^+_1)$, and $(\bm{u},p)$ is a weak solution pair (see Definition \ref{def:wsp}) of the Stokes equations \eqref{Stokes-Eqn} with Navier boundary condition \eqref{NavierBC0}.
Then we have $\pd_t\bm{u}, \nabla^2 \bm{u}, \nabla p\in L^{q,r}(Q^+_{1/2})$, and
\begin{equation}\label{second-est}
\|\pd_t\bm{u}, \nabla^2\bm{u},\nabla p\|_{L^{q,r}\left(Q^+_{1/2}\right)}\lec \|\bm{u}\|_{L^{q,r}(Q^+_1)}+\|p\|_{L^{q,r}(Q^+_1)}+\|\bm{f}\|_{L^{q,r}(Q^+_1)}.
\end{equation}
If in addition $\Gamma \in C^{3,1}$ and  $\nabla \bm{f}\in L^{q,r}(Q^+_1)$, we have 
\begin{align}\label{third-est}
\| \nabla\pd_t\bm{u}, \nabla^3\bm{u},  \nabla^2 p\|_{L^{q,r}\left(Q^+_{1/2}\right)}\lec\|\bm{u}\|_{L^{q,r}(Q^+_1)}+\|p\|_{L^{q,r}(Q^+_1)}+\|\bm{f}\|_{L^r\left(-1,0;W^{1,q}(B^+_1)\right)}.
\end{align}
}
\end{thm}

\emph{Comments on Theorem \ref{thm2_navierBC}}
\begin{enumerate}
\renewcommand{\theenumi}{\roman{enumi}}

\item Theorem \ref{thm2_navierBC} extends \cite[Theorem 1.2]{CLT24} from a flat boundary to a curved boundary, and from 3D to any dimension $d\geq 2$.
We need $\norm{p}_{L^{q,r}}$ on the RHS of \eqref{second-est} and \eqref{third-est} to avoid counter-examples in  \cite{CLT25}.

\item Dong-Kwon \cite[Theorem 2.5]{DK26-2}  study the local estimates near curved boundary for the Stokes equations  under Navier boundary condition with variable viscosity coefficients, and obtain the  second derivative estimate \eqref{second-est} assuming that $\pd_{t}\bm{u},\nabla^2\bm{u},\nabla p$ are bounded in $L^{q,r}(Q^+_{1})$. In contrast, we prove \eqref{second-est} without this assumption. 
\item The reason we require $\CF=0$ is that the higher order estimates might fail if $[\CF \bm{n}]_{tan}|_{\Sigma}\neq 0$.
 See our previous work \cite[Theorem 1.2]{CLT24} for the flat boundary case. If we assume $[\CF \bm{n}]_{tan}|_{\Sigma}= 0$, then $\div \CF$ has the same position with $\bm{f}$, and we can just replace $\bm{f}$ in 
\eqref{second-est} and \eqref{third-est} with $\bm{f}+\div \CF$.

\item
We do not have fourth derivative estimates  in Theorem  \ref{thm2_navierBC} under Navier BC as shown in \cite[Remark 6.3]{CLT25}.

\end{enumerate}

One of the contributions and motivations of this paper is to relax the initial regularity of the solution to the very weak solution class. Unlike in the strong solution setting, the finiteness of certain perturbation terms cannot be assumed and absorbed into the left-hand side, and we have to introduce several techniques, such as normal forms, to deal with this.

In the proof of the previous three theorems, we construct some approximations of the velocity fields by carefully choosing the test functions which satisfy the definitions of (very) weak solutions with no-slip or Navier BC. These approximations satisfy Stokes equations in the whole space and the tricky part is to deal with the 
higher order perturbation terms on the RHS of their equations, which are generated by boundary straightening. 
We introduce several tools that play an important role in dealing with the higher order perturbation terms. Firstly we define the normal form (such as \eqref{normal-form1} and \eqref{heatE-hatw}) of velocity and vorticity, which is more regular than the directly straightened velocity and vorticity, and it can serve as the intermediate vector field for our regularity estimate. 
Secondly, we often prove the horizontal derivative estimates first, and then recover vertical derivative estimates by their bindings with the equations (for example, the estimate of $\pd_d^2 u$ will be the same as $\nabla'\nabla u$ if we have $\Delta u=0$), such as \eqref{hori-1}, \eqref{galdi-trick2}, \eqref{galdi-trick4}. Thirdly, to deal with the lower order terms with temporal derivative, we transfer them to terms with only spatial derivatives, such as \eqref{treat-boundary-t}, \eqref{pd_t-to-pdx}, and then use tricks like \eqref{bad1} to make it small on the RHS, so that it can be absorbed to the LHS.

\medskip

As an application of our previous results, we propose a new definition of boundary regular points for the incompressible Navier--Stokes equations \eqref{NS} that guarantees higher spatial regularity using Theorems \ref{thm_0BC}--\ref{thm2_navierBC}.
Recall that an interior point $(x_0,t_0)\in \Omega \times (0,T]$ is called a  ``regular point'' of a solution $u$ of \eqref{NS} if $u \in L^\infty(Q_R(x_0,t_0))$ for some $R>0$ and $Q_R(x_0,t_0)\subset\Omega \times (0,T)$, and a ``singular point'' otherwise. In a series of papers starting with \cite{Sch76},
Scheffer initiated the study of partial regularity theory of \eqref{NS}. His results were further generalized and strengthened in the celebrated paper by Caffarelli-Kohn-Nirenberg \cite{CKN82}, which showed that the set of possible interior singular points of a suitable
weak solution is of one-dimensional parabolic Hausdorff measure zero. A simplified proof was given by F. Lin  \cite{Lin98}. Since then, many related works emerged and contributed to the partial regularity theory in different directions. For instance, the estimate of the Hausdorff measure was improved by a logarithmic factor in \cite{CL00,LR24}; some new interior regularity criteria were provided in \cite{GKT07,LS99}, and then were extended up to the flat boundary in \cite{GKT06,Ser02b,Ser02,Ser16}, and \cite{SSS04} for a curved boundary.

These results show that a solution is spatially smooth near an interior regular point. One wonders what if the regular point is located on the boundary?  

Scheffer \cite{Sch82} called a boundary point $z_0 \in \pd \Omega$ a regular point to the solution if $\bm{u}$ is bounded in a small neighborhood of $z_0$. 
In the paper of Seregin \cite[Theorem 2.3]{Ser02}, a point $z_0 \in \pd \Omega$ is called a regular point  if $\bm{u}$ is H\"older continuous in $\overline Q_R^+(z_0)$ for some $R>0$. This definition was adapted in \cite{GKT06,SSS04}.
 However, a boundary regular point under these definitions does not guarantee higher spatial regularity as in the interior case, as shown by the counter examples in \cite{CK23,CLT25,Kan04a,KLLT22,KangMin,SS10}, discussed earlier. 
Hence we propose the following alternative definition with higher spatial regularity in mind.

\begin{defn}
We call $z_0=(x_0,t_0)$ with $x_0\in \pd\Omega$ a \emph{boundary regular point} to the solution pair $(\bm{u},p)$ of Navier--Stokes equations \eqref{NS}, if there exists some $R>0$ such that  $\bm{u}\in L^{\infty}(Q^+_{R}(z_0))$ and $p\in L^{q_0,r_0}(Q^+_{R}(z_0))$ for some $1<q_0,r_0<\infty$.
\end{defn}

The definition of interior regular points only concerns the $L^\infty$ norm of $\bm{u}$ and has no assumption on the pressure. Since for suitable weak solutions one often requires that $p\in L^{3/2}$ in space-time, the above condition for $p$ is easy to fulfill. The addition of this condition of pressure enables us to prove the following two theorems for Navier--Stokes equations.

\begin{thm}[NS with no-slip BC]\label{thm4}
\sl{
Suppose that $\pd\Omega\in C^{1,1}$,  $(\bm{u},p)$ is a very weak solution pair  to the Navier--Stokes equations \eqref{NS} with no-slip boundary condition and zero force, $($i.e., it is a very weak solution pair  to the Stokes equations \eqref{Stokes-Eqn} with no-slip boundary condition as in Definition \ref{def:vwsp-0BC}  with force $f=0$ and $\FF = - \bm{u} \otimes \bm{u})$,
%(see Definition \ref{def:vwsp2}), 
 and  $z_0$ is a boundary regular point of $(\bm{u},p)$
with $p\in L^{q_0,r_0}(Q^+_{R}(z_0))$ for some $1<q_0,r_0<\infty$.
 Then there exists some $\rho>0$ such that
$\pd_t\bm{u},\nabla^2\bm{u},\nabla p\in L^{q_1,r_0}(Q^+_{\rho}(z_0))$ for any $q_1<\infty$.
}
\end{thm}

\begin{thm}[NS with Navier BC]\label{thm5}
\sl{
Suppose that $\pd\Omega\in C^{3,1}$,  $(\bm{u},p)$ is a weak solution pair  to the Navier--Stokes equations \eqref{NS} with Navier boundary condition 
and zero force, $($i.e., it is a weak solution pair  to the Stokes equations \eqref{Stokes-Eqn} with Navier boundary condition as in Definition \ref{def:wsp}  with force $f=0$ and $\FF = - \bm{u} \otimes \bm{u})$,
and  $z_0$ is a boundary regular point of $(\bm{u},p)$ with $p\in L^{q_0,r_0}(Q^+_{R}(z_0))$ for some $1<q_0,r_0<\infty$. Then there exists some $\rho>0$ such that
$\nabla\pd_t\bm{u},\nabla^3\bm{u},\nabla^2 p\in L^{q_1,r_0}(Q^+_{\rho}(z_0))$ for any $q_1<\infty$.
}
\end{thm}

\emph{Comments on Theorems \ref{thm4} and \ref{thm5}}
\begin{enumerate}
\item[(i)] 
The proof of Theorems \ref{thm4} and \ref{thm5} follows from Theorems \ref{thm_0BC},  \ref{thm_navierBC} and \ref{thm2_navierBC}:
We first show $\nabla \bm{u}  \in L^{q_0,r_0}$ using \eqref{gra1} and \eqref{main1} with $f=0$ and $\FF =- \bm{u} \otimes \bm{u} \in L^\infty$. We then use 
 \eqref{gra0} and \eqref{second-est} with $f = - \bm{u} \cdot \nabla \bm{u} \in L^{q_0,r_0}$ and $\FF =0$
 to get 
$\pd_t\bm{u},\nabla^2\bm{u},\nabla p\in L^{q_0,r_0}$. It follows by imbedding that $p  \in L^{q_1,r_0}$ for some $q_1>q_0$.
We  iterate this process to get $p \in L^{q_1,r_0}$ for a sequence of $q_1$ increasing to $\infty$. For the Navier BC case, we further use \eqref{third-est} to get 
$\nabla\pd_t\bm{u},\nabla^3\bm{u},\nabla^2 p\in L^{q_1,r_0}$.

\item[(ii)] If we only want second derivative estimates, we can relax to $\pd\Omega\in C^{2,1}$ in Theorem  \ref{thm5}. 

\item[(iii)] As stated in the comments after Theorems \ref{thm_0BC} and  \ref{thm2_navierBC}, we do not expect third derivative estimate for the no-slip BC case in Theorem \ref{thm4}, and we do not expect fourth derivative estimate for the Navier BC case in Theorem \ref{thm5}.

\end{enumerate}

The rest of the paper is organized as follows. In Section \ref{sec2} we introduce notations and some preliminary results, and give the definitions of (very) weak solutions of Stokes systems under no-slip and Navier boundary conditions. In Section \ref{sec3} we prove Theorem \ref{thm_0BC} under no-slip BC: we construct the approximations of velocity first, prove the gradient estimate in the second subsection, and prove the second derivative estimate in the third subsection. Sections \ref{sec4} and  \ref{sec5} are under Navier BC.
In Section \ref{sec4} we prove Theorem \ref{thm_navierBC}: we construct the approximations of vorticity first, prove the gradient estimate with an extra term on the RHS in the second subsection, and remove the extra term in the third subsection. In Section \ref{sec5} we prove Theorem \ref{thm2_navierBC}: we construct the approximations of velocity first, prove the second derivative estimate in the second subsection, and prove the third derivative estimate in the third subsection. In Section \ref{sec6} we give several useful auxiliary lemmas.

\section{Notations and preliminaries }\label{sec2}
\subsection{Notations}\label{sec2.1}
We first clarify the notation system used in this paper.
For two comparable quantities, the inequality $X\lesssim Y$ stands for $X\leqslant C Y$ for some positive constant $C$. The dependence of the constant $C$ on other parameters or constants is usually clear from the context, and we will often suppress this dependence.

For $x\in \R^d$, we write $\pd_{x_i}=\frac{\partial}{\partial x_{i}}$ as the partial derivative with respect to the variable $x_{i}, 1 \leq i \leq d$, and $\nabla_{x}=\left(\pd_{x_{1}},\cdots,\pd_{x_{d}}\right)$ as the gradient. Denote $x^{\prime}=(x_{1},\cdots,x_{d-1})$, and $\nabla'=\nabla_{x'}=\left(\pd_{x_{1}},\cdots,\pd_{x_{d-1}}\right)$ as the horizontal gradient. For a multi-index $\alpha=(\alpha_1,\ldots,\alpha_d)\in \mathbb N^d$, we use the standard notation
$\nabla^\alpha f=\partial_1^{\alpha_1}\cdots \partial_d^{\alpha_d} f$.
On the other hand, for an integer $m\geq 1$, the notation $\nabla^m f$ denotes the collection of all $m$-th order spatial derivatives of $f$. In particular, $\nabla^2 f$ denotes the Hessian matrix of $f$ when $f$ is scalar-valued. If $\bm{f}$ is vector-valued, then $\nabla^m \bm{f}$ is understood componentwise.

For any norm $\norm{\cdot}$, we adopt the convention $\norm{f,g,h}=\norm{f}+\norm{g}+\norm{h}.$

We use the convention that the gradient matrix $(\nabla u)_{ij} = \pd_j u_i$, the curl matrix $\left(\curl \bm{u}\right)_{ij}=\pd_{i}u_j-\pd_{j}u_i$, $
(\div \CF)_i ={\textstyle \sum_j}\, \pd_j \CF_{ij}$, $\left(\bm{u}\otimes\bm{v}\right)_{ij}=u_i\cdot v_j$ and $\CF:\GG={\textstyle \sum_{ij}}\, F_{ij}\cdot G_{ij}$. We have $\int \div \FF \cdot u = 
-\int \FF : \nabla u $ for $u \in C^1_c$.

For an open domain $\Omega\subset \R^d$, $1\leq q,r \leq \infty$ and integer $m\geq 1$, let $L^{q}(\Omega)$ and $W^{m,q}(\Omega)$ be the usual (scalar or vector-valued) Lebesgue and Sobolev spaces, respectively. Denote by $L^{q,r}\left(\Omega\times I\right)=L^{r}\left(I;L^{q}(\Omega)\right)$ a mixed-exponents Lebesgue space for $I\subset \R$ and  $L^{q}\left(\Omega\times I\right)=L^{q,q}\left(\Omega\times I\right)$. 
Let 
\[W_{c,\infty}^{1}\left(\Omega\times I\right)=\{g\, |\,\hbox{supp}(g)\subset\subset\left(\Omega\times I\right);\ g,\,\nabla g, \,\partial_{t}g\in L^{\infty}\left(\Omega\times I\right)\},\] 
\[W_{c,\infty}^{2,1}\left(\Omega\times I\right)=\{g\in W_{c,\infty}^{1}\left(\Omega\times I\right)\,|\,\nabla^2 g\in L^{\infty}\left(\Omega\times I\right)\}.\]

Denote the ball $B_R(x_0)=\{x\in\R^{d}:|x-x_0|<R\}$ for $x_0\in \R^d$; $B^+_R(x_0)=B_R(x_0)\cap \Omega$ for $x_0\in \pd\Omega$ (also denoted as $\Omega_R(x_0)$ in much of the literature), and  $Q^+_R(z_0)=B^+_R(x_0)\times(t_0-R^2,t_0)$ for $z_0=(x_0,t_0)$; $B_R^{\prime}=\{|x'|<R\}$; $\C_R=\{|x'|<R\}\times (-R,R)$ and $\C_R^+=\C_R\cap \R_{+}^d$;   $\Gamma=B_1\cap\pd\Omega$, and $\Sigma=\Gamma\times(-1,0)$.

For abbreviation of notation in later calculation, we denote 
\begin{equation*}
\cW^{q,r}\eqdefa L^{q,r}\left(\R^d\times(-1,0)\right),\quad \cH^{q,r}\eqdefa L^{q,r}\left(\R^d_+\times(-1,0)\right).
\end{equation*}

\subsection{Boundary straightening, cut-off functions and mollification}
Suppose $\pd\Omega\in C^{1,1}$, then by definition for each $x_0\in \pd\Omega$ there exists $\rho>0$ such that upon translating and reorienting the coordinate axes if necessary, we have 
\begin{align}\label{condition-B_R+}
x_0=(0,...,0),\quad B^+_\rho=\{x\in B_\rho \ | \ x_d>\gamma(x')\},
\end{align}
where a function $\gamma \in C^{1,1}(\R^{d-1})$ satisfies that
\begin{equation}\label{gamma}
\gamma|_{x'=0}=0,\quad \nabla\gamma|_{x^{\prime}=0}=0,\quad \|\gamma\|_{C^{1,1}(\R^{d-1})}\lec 1,
\end{equation}
and $\pd\Omega$ does not intersect $\pd B_\rho$ multiple times, i.e. $\pd\Omega \cap \pd B_\rho$ is homeomorphic to $\mathbb{S}^{d-2}$ (two points for $d=2$).

With the set-up above, we define a map $Tx=(x',x_d-\gamma(x'))$ for $x\in B^+_\rho$ and its inverse $T^{-1}x=(x',x_d+\gamma(x'))$ for $x\in T(B^+_\rho)\subset \R^d_+$. Denote $Tf(x)=f(Tx)$,  $T^{-1}f(x)=f(T^{-1} x)$. Notice that, with $T$ and $T^{-1}$ considered as operators on functions, for $1\leq k\leq d-1$
\begin{equation}\label{commu1}
\pd_{k}\, T=T\pd_{k}-\gamma_k^{\prime}\cdot T \pd_{d},\quad \pd_{d}\,T=T  \pd_{d},
\end{equation}
and
\begin{equation}\label{commu2}
\pd_{k}\, T^{-1}=T^{-1}\pd_{k}+\gamma_k^{\prime}\cdot T^{-1}\,\pd_{d},\quad \pd_{d}\,T^{-1}=T^{-1}\,\pd_{d},
\end{equation}
where $\gamma'_k=\pd_k\gamma $.
We agree $\gamma'_d=0$ so that \eqref{commu1}-\eqref{commu2} are valid for $k=d$, and $\nabla \gamma=(\nabla' \gamma,0)$.

\begin{lem}\label{lem21}
\sl{
\begin{equation}\label{change-T-Lap}
T^{-1}\Delta=\Big(\Delta-2\nabla\gamma\cdot \nabla\, \pd_{d}+|\nabla\gamma|^2\cdot \pd_{d}^2-\Delta\gamma\cdot \pd_{d} \Big)T^{-1}.
\end{equation}
 }
\end{lem}
\begin{proof}
In view of \eqref{commu2}, we get that for $1\leq k \leq d-1$,
\begin{align}
T^{-1}\pd_{k}^2 &= \left(\pd_{k}-\gamma_k^{\prime}\cdot \pd_{d}\right)T^{-1}\pd_{k} \notag\\[6pt]
&=\left(\pd_{k}-\gamma_k^{\prime}\cdot \pd_{d}\right)\left(\pd_{k}-\gamma_k^{\prime}\cdot \pd_{d}\right)T^{-1}\notag\\[6pt]
&=\left(\pd_{k}^2 -2\gamma'_k\cdot \pd_{k} \pd_{d}+(\gamma'_k)^2\cdot \pd_{d}^2-\gamma''_{kk}\cdot \pd_{d} \right)T^{-1} ,\label{T1}
\end{align}
and
$
T^{-1}\pd_{d}^2=\pd_{d}^2\,T^{-1}$.
Summing up these estimates, we arrive at \eqref{change-T-Lap}.
\end{proof}

We denote  $\tilde{\bm{u}}=(\tilde{u}_1,\cdots,\tilde{u}_{d})$  with 
\begin{equation}\label{deftilde-u}
\tilde{u}_{k}= \left\{\begin{array}{cc}
u_{k},\quad &1\leq k \leq d-1 \\[6pt]
-\sum_{i=1}^{d-1}\,\gamma'_i\cdot u_i+u_d,\quad&k=d
\end{array}\right..
\end{equation}
In view of \eqref{commu2}, it's easy to see that 
\begin{equation}\label{div2}
\div\, T^{-1}\tilde{\bm{u}}=0,
\end{equation}
if $\div \bm{u}=0$, 
and $T^{-1}\tilde{u}_{d}|_{y_d=0}=0$ if $\bm{u}\cdot \bm{n}|_{\Sigma}=0$.

\smallskip

We now introduce some cut-off functions which will be used throughout this paper. Fix a smooth function $\varphi(s)$ such that $\varphi(s)=1$ for $s<0$, $\varphi(s)=0$ for $s>1$ and $\varphi'(s)\leq 0$. Denote the spatial cut-off function $\zeta_{\rho_1,\rho}(s)=\varphi(\frac{s-\rho_1}{\rho-\rho_1})$ for $s\geq 0$ and  $0<\rho_1<\rho$, which satisfies $\zeta_{\rho_1,\rho}\in C^{\infty}_c[0,\rho)$ with $\zeta_{\rho_1,\rho}=1$ on $[0,\rho_1]$, and its $m$-th order derivative satisfies
\[
|\zeta^{(m)}_{\rho_1,\rho}|\lesssim (\rho-\rho_1)^{-m}.
\] 
Denote the temporal cut-off function $\psi_{\tau}(t)=\varphi(-\frac{t+(1-\tau/4)}{1/4})$ for $t\leq 0$ and $\tau=1,2,3$, which satisfies $\psi_\tau(t)\in C^{\infty}_c(-\frac{5-\tau}{4},0]$ with $\psi_\tau(t)=1$ on $[-\frac{4-\tau}{4},0]$. Let
\begin{equation}\label{def-zeta}
\zeta(x,t)=\zeta_{\rho_1,\rho}(|x'|)\cdot\zeta_{\rho_1,\rho}(|x_d|)\cdot\psi_{\tau}(t).
\end{equation}

For $0<\ep\ll 1$, denote $\eta^{h}_{\ep}(x')=c_1\ep^{-(d-1)}\varphi(\frac{|x'|}{\ep})$, $\eta_{\ep}(t)=c_2 \ep^{-1}\varphi(\frac{|t|}{\ep})$, where $c_1, c_2$ are constants such that
\[\int_{\R^{d-1}}\eta^{h}_{\ep}(x')\dx'=\int_{\R}\eta_{\ep}(t)\dt=1.\]
For $0<\ep_1,\ep_2,\ep_3\ll 1$, let 
\begin{equation}
\begin{aligned}
\eta_{\ep_1,\ep_2,\ep_3}^{\pm}(x,y,t)&=\eta_{\ep_1}^{h}(x'-y')\cdot\Big(\eta_{\ep_2}(x_d-y_d)\pm\eta_{\ep_2}(x_d+y_d)\Big)\cdot\eta_{\ep_3}(t).
\end{aligned}
\end{equation}
We denote
\begin{equation}\label{def-xi}
\xi^{\pm}(x,y,t,s)=\zeta(y,s)\cdot\eta^{\pm}_{\ep_1,\ep_2,\ep_3}(x,y,t-s).
\end{equation}

For $g\in L^{1}_{\text{loc}}(\overline{\R^d_{+}}\times \R)$ and $(x,t)\in\R^d\times \R$, we define the mollification (with even or odd extension),
\begin{align}\label{Eg}
E^{\pm}_{\ep_1,\ep_2,\ep_3}(g)(x,t)&\eqdefa\int_{\R^d_{+}\times\R}g(y,s)\cdot\eta_{\ep_1,\ep_2,\ep_3}^{\pm}(x,y,t-s)\dy\ds\notag\\[6pt]
&=E^h_{\ep_1}\left(E_{\ep_2}^{\pm}\left(E_{\ep_3}\left(g\right)\right)\right),
\end{align}
where the horizontal mollification
\begin{equation*}
E^h_{\ep_1}(g)\eqdefa \int_{\R^{d-1}}g(y',x_d,t)\cdot\eta_{\ep_1}^{h}(x'-y')\dy',
\end{equation*}
the vertical mollification
\begin{equation*}
E_{\ep_2}^{\pm}(g)\eqdefa \int_{\R_+}g(x',y_d,t)\cdot\Big(\eta_{\ep_2}(x_d-y_d)\pm\eta_{\ep_2}(x_d+y_d)\Big)\dy_d,
\end{equation*}
and the temporal mollification
\begin{equation*}
E_{\ep_3}(g)\eqdefa \int_{\R}g(x,s)\cdot\eta_{\ep_3}(t-s)\ds.
\end{equation*}
We also denote $E^{h,\pm}_{\ep_1,\ep_2}(g)\eqdefa E^h_{\ep_1}\left(E^{\pm}_{\ep_2}(g)\right)$, $E^{\pm}_{\ep_2,\ep_3}(g)\eqdefa E^{\pm}_{\ep_2}\left(E_{\ep_3}(g)\right)$, $E^{h}_{\ep_1,\ep_3}(g)\eqdefa E^{h}_{\ep_1}\left(E_{\ep_3}(g)\right)$. The following lemma is frequently  used in this paper to handle the terms after mollification. Recall \eqref{q_*} that $1<q,r<\infty$.
\begin{lem}\label{lem0} 
\sl{Suppose $g\in  L^{q,r}(\R^d_{+}\times\R)$, $1<q,r<\infty$. We have the following properties.
\begin{enumerate}
\item[(\rm  i)] $\pd_{x_{k}}\,\eta^{\pm}_{\ep_1,\ep_2,\ep_3}(x,y,t)=-\pd_{y_{k}}\,\eta_{\ep_1,\ep_2,\ep_3}^{\pm}(x,y,t)$ for $1\leq k\leq d-1$; $\pd_{x_{d}}\eta_{\ep_1,\ep_2,\ep_3}^{\pm}(x,y,t)=-\pd_{y_{d}}\eta_{\ep_1,\ep_2,\ep_3}^{\mp}(x,y,t)$. The second property implies: $\pd_d E^+_{\ep_1,\ep_2,\ep_3}(g)=E^-_{\ep_1,\ep_2,\ep_3}(\pd_d g)$ and $\pd_d E^-_{\ep_1,\ep_2,\ep_3}(g)=E^+_{\ep_1,\ep_2,\ep_3}(\pd_d g)$ if $g|_{y_d=0}=0$.
\vspace{6pt}
\item[(\rm ii)] $E^{+}_{\ep_1,\ep_2,\ep_3}(g),\ E^{-}_{\ep_1,\ep_2,\ep_3}(g)$ are smooth in $\R^d\times \R$ and converge to $g(x,t)$ in $L^{q,r}(\R_{+}^d\times \R)$  as $(\ep_1,\ep_2,\ep_3) \to 0$ for  $1\leq q,r<\infty$.
\vspace{6pt}
\item[(\rm iii)] $\|E^{\pm}_{\ep_1,\ep_2,\ep_3}(g)\|_{L^{q,r}(\R^d\times \R)}\lec\|g\|_{L^{q,r}(\R^d_{+}\times\R)}$; 
\vspace{6pt}
\item[(iv)] If $\|\nabla_{x}E^{\pm}_{\ep_1,\ep_2,\ep_3}(g)\|_{L^{q,r}(\R^d\times \R)}\lec C$ for some constant $C$ independent of $(\ep_1,\ep_2, \ep_3)$, then the weak derivative $\nabla g$ exists and is in $L^{q,r}(\R_{+}^d\times \R)$.
\vspace{6pt}
\item[(v)] Suppose $a\in W^{1,\infty}(\R^{d-1})$, $1\leq k\leq d-1$, we have  the estimate of commutator
\begin{equation}\label{commu3}
\|\pd_{x_k}E^h_{\ep_1}(a\cdot g)-a \cdot \pd_{x_k}E^h_{\ep_1}(g)\|_{L^{q,r}(\R^d_+\times \R)}\lec \|\nabla a\|_{L^{\infty}(\R^{d-1})}\cdot\|g\|_{L^{q,r}(\R^d_{+}\times\R)}.
\end{equation}
\end{enumerate}
}
\end{lem}
\begin{proof}
(i)--(iv) are elementary, so we skip their proof. For (v), we have 
\begin{align*}
&\quad\, \pd_{x_k}E^h_{\ep_1}(a \cdot g)-a(x')\cdot \pd_{x_k}E^h_{\ep_1}(g)\\[6pt]
&= -\int_{\R^{d-1}}\frac{a(x')-a(y')}{|x'-y'|}\cdot g(y',x_d,t)\cdot |x'-y'|\cdot \pd_{x_k}\,\eta_{\ep_1}^{h}(x'-y')\dy'.
\end{align*}
Hence, by using Young's inequality, we have
\begin{align*}
&\quad\ \|\pd_{x_k}E^h_{\ep_1}(ag)-a(x')\pd_{x_k}E^h_{\ep_1}(g)\|_{L^{q,r}(\R^d_+\times \R)} \\[6pt]
&\lec \|\nabla a\|_{L^\infty(\R^{d-1})}\cdot \ep_1\cdot\left\|\int_{\R^{d-1}} |g(y',x_d,t)|\cdot|\pd_{x_k}\eta_{\ep_1}^{h}(x'-y')|\dy'\right\|_{L^{q,r}(\R^d_+\times \R)}\\[6pt]
&\lec  \|\nabla a\|_{L^{\infty}(\R^{d-1})}\cdot\|g\|_{L^{q,r}(\R^d_{+}\times\R)}. \qedhere
\end{align*} 
\end{proof}

\subsection{Various definitions of weak solutions}

We introduce various notions of weak solutions for the non-stationary Stokes equations \eqref{Stokes-Eqn} with no-slip boundary condition \eqref{0-BC0} or Navier boundary condition \eqref{NavierBC0}.
\begin{defn}\label{def:vwsp-0BC}
\sl{Let $\bm{f},\CF\in L^{1}(Q^+_1)$. A pair $(\bm{u},p)$ is called a \emph{very weak solution pair} to the Stokes equations \eqref{Stokes-Eqn} in $Q^+_1$ with \emph{no-slip} boundary condition \eqref{0-BC0} on the boundary $\Sigma$, if
\begin{enumerate}
\item[(i)] $\bm{u},p \in L^{1}(Q^+_1)$ satisfy that for a.e.~$t\in (-1,0)$ and every scalar function $\Psi \in C_{c}^{1}\big(B^+_1\cup\Gamma\big)$, which may not vanish on the boundary $\Gamma$,
\begin{equation}\label{div-free weak}
\int_{B^+_1} \bm{u} \cdot \nabla \Psi\dx=0.
\end{equation}

\item[(ii)]   for each vector $\Phi=\left(\Phi_{1},\cdots,\Phi_{d}\right) \in W_{c,\infty}^{2,1}\left( Q^+_1\cup \Sigma \right)$ with $\Phi|_{\Sigma}=0$, we have
\begin{equation}\label{weakform-0BC}
-\int_{Q^+_1} \big(\bm{u}\cdot \left(\pd_t+\Delta\right)\Phi+p\cdot\div \Phi\big)\dx\dt=\int_{Q^+_1}\left( \bm{f}\cdot\Phi-\CF:\nabla\Phi \right)\dx\dt.
\end{equation}
\end{enumerate}
}
\end{defn}

\begin{defn}\label{def:wsp-0BC}
\sl{Let $\bm{f},\CF\in L^{1}(Q^+_1)$. A pair $(\bm{u},p)$ is called a \emph{weak solution pair} to the Stokes equations \eqref{Stokes-Eqn} in $Q^+_1$ with \emph{no-slip} boundary condition \eqref{0-BC0} on the boundary $\Sigma$, if
\begin{enumerate}
\item[(i)] $\bm{u},\nabla \bm{u}, p \in L^{1}\big(Q^+_1\big)$ with $\div \bm{u}=0$ and $\bm{u}|_{\Sigma}=0$,

\item[(ii)]   for each vector $\Phi=\left(\Phi_{1},\cdots,\Phi_{d}\right) \in W_{c,\infty}^{1}\left( Q^+_1\cup \Sigma \right)$ with $\Phi|_{\Sigma}=0$, we have
\begin{equation}\label{weakform2-0BC}
\int_{Q^+_1} \big(-\bm{u}\cdot \pd_t\Phi+\nabla\bm{u}:\nabla\Phi-p\cdot\div \Phi\big)\dx\dt=\int_{Q^+_1}\left( \bm{f}\cdot\Phi-\CF:\nabla\Phi \right)\dx\dt.
\end{equation}
\end{enumerate}
}
\end{defn}

\begin{defn}\label{def:ws}
\sl{Let $\bm{f},\CF\in L^{1}(Q^+_1)$. A vector field $\bm{u}$ is called a \emph{weak solution} to the Stokes equations \eqref{Stokes-Eqn} in $Q^+_1$ with \emph{Navier} boundary condition \eqref{NavierBC0} on the boundary $\Sigma$, if
\begin{enumerate}
\item[(i)] $\bm{u},\nabla \bm{u}\in L^{1}\big(Q^+_1\big)$ with $\div \bm{u}=0$ and $\bm{u}\cdot\bm{n}|_{\Sigma}=0$;

\item[(ii)]  for each vector $\Phi=\left(\Phi_{1},\cdots,\Phi_{d}\right) \in W_{c,\infty}^{1}\left( Q^+_1 \cup \Sigma \right)$ with $\div \Phi=0$ and  $\Phi\cdot\bm{n}|_{\Sigma}=0$, we have
\begin{align}\label{weakform-NavierBC}
&\ \quad \int_{Q^+_1} \left(-\bm{u} \cdot \partial_{t}\Phi+2\DD \bm{u}: \DD\Phi\right)\dx\dt +\int_{\Sigma}\alpha \bm{u}_{tan}\cdot\Phi_{tan}\dSx\dt\notag\\[6pt]
&=\int_{Q^+_1}\left( \bm{f}\cdot\Phi-\CF:\nabla\Phi \right)\dx\dt.
\end{align}
\end{enumerate}
}
\end{defn}

\begin{defn}\label{def:wsp}
\sl{Let $\bm{f},\CF\in L^{1}(Q^+_1)$. A pair $\left(\bm{u},p\right)$ is called a \emph{weak solution pair} to the Stokes equations \eqref{Stokes-Eqn} in $Q^+_1$ with \emph{Navier} boundary condition \eqref{NavierBC0} on the boundary $\Sigma$, if
\begin{enumerate}
\item[(i)] $\bm{u} ,\nabla \bm{u}, p\in L^{1}(Q^+_1)$ with $\div \bm{u}=0$ and $\bm{u}\cdot\bm{n}|_{\Sigma}=0$;

\item[(ii)]  for each vector $\Phi=\left(\Phi_{1},\cdots,\Phi_{d}\right) \in W_{c,\infty}^{1}\left( Q^+_1 \cup \Sigma \right)$ with $\Phi\cdot\bm{n}|_{\Sigma}=0$, we have
\begin{align}\label{weakform-pair}
&\quad \ \int_{Q^+_1} \left(-\bm{u} \cdot \partial_{t}\Phi+2\DD \bm{u}: \DD \Phi -p\cdot \mathrm{div}\Phi\right)\dx\dt+\int_{\Sigma}\alpha \bm{u}_{tan}\cdot\Phi_{tan}\dSx\dt\notag\\[6pt]
&=\int_{Q^+_1} \left( \bm{f}\cdot\Phi-\CF:\nabla\Phi \right)\dx\dt .
\end{align}
\end{enumerate}
}
\end{defn}

\section{Proof of Theorem \ref{thm_0BC}}\label{sec3}
In this section we give the proof of Theorem \ref{thm_0BC} on derivative estimates of very weak solution pair of \eqref{Stokes-Eqn} under no-slip BC. With the interior estimate \eqref{interior1}, \eqref{interior3} and standard covering argument, we only need to prove the local boundary estimate:  for each $x_0\in B_{1/2}\cap\pd\Omega$, there exists $\rho=\rho(x_0)>0$ such that $\pd_{t} \bm{u},\nabla^2\bm{u},\nabla p \in L^{q,r}\big(B^+_{\rho/10}(x_0)\times(-\frac14,0)\big)$ and
\begin{equation}\label{gra2}
\|\nabla \bm{u}\|_{L^{q,r}\big(B^+_{\rho/{10}}(x_0)\times(-\frac34,0)\big)} \lec\rho^{-1}\|\bm{u}\|_{L^{q,r}(Q^+_1)}+\|p,\CF\|_{L^{q,r}(Q^+_1)}+\|\bm{f}\|_{L^{q_*,r}(Q^+_1)},
\end{equation}
\begin{align}\label{gra3}
&\quad\ \|\pd_{t} \bm{u},\nabla^2\bm{u},\nabla p\|_{L^{q,r}\big(B^+_{\rho/10}(x_0)\times(-\frac14,0)\big)}\lec \rho^{-3-\frac{1}{q-1}}\|\bm{u}, p, \bm{f}\|_{L^{q,r}(Q^+_1)}.
\end{align}
Without loss of generality, we may assume  \eqref{condition-B_R+}-\eqref{gamma}. 
 We require $\rho$ to satisfy some smallness conditions, including
 \begin{equation}\label{rho1}
T(B^+_{\rho/10})\subset \C^+_{\rho/2} ,\quad \C^+_{\rho}\subset T(B^+_1),
\end{equation}
and some others which are to be specified later in subsection \ref{sec3.2}. Note that $\gamma$ defined in \eqref{gamma} always satisfies
\begin{equation}\label{rho2}
\|\nabla\gamma\|_{L^{\infty}(B'_{\rho})}\lec \rho \ll 1.
\end{equation}

\subsection{Heat equations of the approximations of velocity}

The idea of proof is to study the approximation of $(\bm{u},p)$, which is defined in the straightened domain $\R^d_+\times (-1,0)$, with mollification and odd extension in the lower half space. By choosing appropriate test functions, we will derive the desired equation of mollified $(\bm{u},p)$.  To this end, define $\bm{v}=(v_1,\cdots,v_d)$ and $\pi$ with
\begin{equation}\label{def-vpi}
v_k(x,t)=T^{-1}u_k(x,t)\cdot \zeta(x,t),\quad \pi(x,t)=T^{-1}p(x,t)\cdot\zeta(x,t),
\end{equation}
for $1\leq k \leq d$, where $(x,t)\in\R^d_{+}\times\R_-$,
 the cut-off function $\zeta(x,t)$ is defined in \eqref{def-zeta}, and $\rho_1, \tau$ will be specified later. It's easy to get from  \eqref{rho1} and \eqref{def-vpi} that
\begin{equation}\label{add1}
\|v_k\|_{\cH^{q,r}}\lesssim \|u_k\|_{L^{q,r}(Q^+_1)}.
\end{equation}

Changing variables in
the weak form \eqref{weakform-0BC}, using $\det(T^{-1})=1$, \eqref{commu2} and \eqref{change-T-Lap}, we get the equation
\begin{align}
&\,\quad -\int_{\R^d_{+}\times\R} T^{-1}\bm{u}\cdot\Big( \partial_{s}+\Delta -2\nabla\gamma\cdot \nabla\, \pd_{d}+|\nabla\gamma|^2\cdot \pd_{d}^2-\Delta\gamma\cdot \pd_{d} \Big)T^{-1}\Phi\dy\ds  \notag\\[6pt]
&\,\quad -\int_{\R^d_{+}\times\R} T^{-1}p\cdot\Big(\div (T^{-1}\Phi)-\nabla \gamma \cdot \pd_d\,T^{-1}\Phi\Big)\dy\ds \label{ws1}\\[6pt]
&=\int_{\R^d_{+}\times\R}\Big(T^{-1}\bm{f}\cdot T^{-1}\Phi-T^{-1}\CF:\left(\nabla\, T^{-1}\Phi- (\pd_d\,T^{-1}\Phi)\otimes\nabla\gamma\right)\Big)\dy\ds, \notag
\end{align} 
for any $\Phi \in W_{c,\infty}^{2,1}\left( Q^+_1\cup \Sigma \right)$.
For $1\leq k\leq d$, we choose the following test function
\begin{align}\label{testF-0BC}
T^{-1}\Phi(y,s)=\xi^{-}(x,y,t,s)(0,\cdots,0,\underbrace{1}_{k-th},0,\cdots,0).
\end{align}
where $\xi^-$ is defined in \eqref{def-xi}. It satisfies the condition $\Phi|_{\Sigma}=0$. For the other condition $\Phi \in W_{c,\infty}^{2,1}\left( Q^+_1\cup \Sigma \right)$, we have the following remark.

\begin{rem}\label{ep0-argument}
When $s\to 0^-$, the test function may not vanish and this violates Definition \ref{def:vwsp-0BC}. One way to deal with it is to let $t\in(-1,-\ep_0)$ for fixed $0<\ep_3<\ep_0\ll 1$. Then our test function is legitimate since $\eta^{\pm}_{\ep_1,\ep_2,\ep_3}(x,y,t-s)=0$ for $s \in(\ep_3-\ep_0,0]$. We can obtain the desired estimate \eqref{gra3} for $t\in(-\frac14,-\ep_0)$ in the LHS, by 
the rest of the arguments in this section. This estimate is uniform in $\ep_0$ and hence we can let $\ep_0\to 0^+$. 
For simplicity, we shall always skip this step afterwards, at the cost of sacrificing a little bit of rigor in the proof. In what follows we assume the time variable  $t\in (-1,0)$ and $\ep_1,\ep_2\ll\rho$.
\hfill $\square$
\end{rem}

By substituting \eqref{testF-0BC} into the equation \eqref{ws1}, and using Lemma \ref{lem0} (i), we derive the heat equation for $E^{-}_{\ep_1,\ep_2,\ep_3}(v_k)$ defined by \eqref{Eg},
\begin{equation}\label{heatE-OBC}
(\pd_t-\Delta_x)\,E^{-}_{\ep_1,\ep_2,\ep_3}(v_k)=K_1+K_2+K_3+K_4,
\end{equation}
where
\begin{align}
K_1&\eqdefa2\int_{\R^d_{+}\times\R} T^{-1}u_k\cdot\nabla \zeta\cdot\nabla_{y}\eta_{\ep_1,\ep_2,\ep_3}^{-}(x,y,t-s)\dy\ds\notag\\[6pt]
&\quad\,+\int_{\R^d_{+}\times\R} T^{-1}u_k\cdot(\pd_s+\Delta)\,\zeta\cdot\eta_{\ep_1,\ep_2,\ep_3}^-(x,y,t-s)\dy\ds   \notag\\[6pt]
&=-2\sum_{i=1}^{d-1}\pd_{x_i}\,E^{-}_{\ep_1,\ep_2,\ep_3}\left(T^{-1}u_k\cdot \pd_i\,\zeta\right)-2\pd_{x_d}\,E^{+}_{\ep_1,\ep_2,\ep_3}\left(T^{-1}u_k \cdot \pd_d\,\zeta\right)\notag\\[6pt]
&\quad\,+E^{-}_{\ep_1,\ep_2,\ep_3}\left(T^{-1}u_k\cdot(\pd_s+\Delta)\,\zeta\right), \label{K1}
\end{align}
and
\begin{align}
K_2&\eqdefa\int_{\R^d_{+}\times\R}T^{-1}u_k\cdot\Big(-2\nabla\gamma\cdot \nabla_y\, \pd_{y_d}+|\nabla\gamma|^2\cdot \pd_{y_d}^2-\Delta\gamma\cdot \pd_{y_d} \Big)\xi^{-}(x,y,t,s)\dy\ds  \notag\\[6pt]
&=K_2^{high}+K_{2}^{low},\label{K2}
\end{align}
with the higher order term
\begin{align*}
K_2^{high}\eqdefa-2\sum_{i=1}^{d-1}\pd_{x_i}\pd_{x_d}\,E^{+}_{\ep_1,\ep_2,\ep_3}\left(\gamma'_i\cdot v_k\right)+\pd_{x_d}^2\,E^{-}_{\ep_1,\ep_2,\ep_3}\left(|\nabla\gamma|^2\cdot v_k\right),
\end{align*} 
and lower order term
\begin{align*}
K_{2}^{low}&\eqdefa -E^{-}_{\ep_1,\ep_2,\ep_3}\Big(T^{-1}u_k \cdot\left(2\nabla\gamma\cdot \nabla\, \pd_{d}-|\nabla\gamma|^2\cdot \pd_{d}^2+\Delta\gamma\cdot \pd_{d} \right)\zeta\Big)\notag\\[6pt]
&\quad\,+2\sum_{i=1}^{d-1}\pd_{x_i}\,E^{-}_{\ep_1,\ep_2,\ep_3}\left(T^{-1}u_k\cdot \gamma'_i\cdot \pd_d\,\zeta\right)\notag\\[6pt]
&\quad\,+\pd_{x_d}\,E^{+}_{\ep_1,\ep_2,\ep_3}\Big(T^{-1}u_k\cdot\left(2\nabla\gamma\cdot \nabla-2|\nabla \gamma|^2\cdot\pd_{d}+\Delta\gamma\right)\zeta\Big), \label{K2l}
\end{align*}
and
\begin{equation}\label{K3}
K_3\eqdefa\int_{\R^d_{+}\times\R} T^{-1}p\cdot\Big(\pd_{y_k}- \gamma_k' \cdot \pd_{y_d}\Big)\xi^{-}(x,y,t,s)\dy\ds, 
\end{equation}
and
\begin{equation}\label{K4}
K_4\eqdefa\int_{\R^d_+\times \R}\Big(T^{-1}f_k-\sum_{i=1}^{d}T^{-1}\CF_{ki}\cdot\left(\pd_{y_i}-\gamma_i'\cdot \pd_{y_d}\right)\Big)\xi^{-}(x,y,t,s)\dy\ds. 
\end{equation}
In the forms of $K_1$ and $K_2$, we have moved derivatives away from the mollification kernel $\eta_{\ep_1,\ep_2,\ep_3}^\pm$ to avoid negative powers of $\ep_i$ in the estimate.
The forms of $K_3$ and $K_4$ can also be expressed similarly and it's left to the reader.

\subsection{Gradient estimate}
We first give the gradient estimate \eqref{gra1}. The main difficulty here is that we don't know whether the LHS is finite or not, so it takes some effort to deal with the higher order term $K_2^{high}$ and reach \eqref{claim1}. Note that the proof of \eqref{gra1} only relies on the property of heat equation, and does not use the div-free property of $\bm{u}$. 

\begin{proof}[Proof of \eqref{gra1}]
We divide the proof into two steps: in step 1 we raise up the regularity of $\nabla \bm{u}$, and in step 2 we show the estimate.  For $\bm{v}$ defined in \eqref{def-vpi}, let $\rho_1=\rho/2$ and $\tau=1$.

\smallskip

\noindent
\textbf{\textit{Step 1: }} we claim that for $1\leq k \leq d$ and $\ep_3>0$,
\begin{equation}\label{claim1}
\nabla E_{\ep_3}(v_k)\in  \cH^{q,r}.
\end{equation}

\smallskip

For the RHS of \eqref{heatE-OBC}, the lower order terms $K_1$ and $K_2^{low}$, the pressure term $K_3$ and the external force term $K_4$ are easy to handle. Later we will see that the first  term in $K_2^{high}$ can be handled by the commutator estimate \eqref{commu3}. The key difficulty comes from the second term in $K_2^{high}$,
as $\norm{\pd_d E^{-}_{\ep_1,\ep_2,\ep_3} (|\nabla \gamma|^2v)}_{\cW^{q,r}}$ is not bounded by $\norm{\nabla \gamma}_{L^\infty} ^2\norm{\pd_d E^{-}_{\ep_1,\ep_2,\ep_3} (v)}_{\cW^{q,r}}$ plus lower order terms.
We will deal with it by moving this term to the LHS of \eqref{heatE-OBC}, which may be considered a \emph{normal form}. Denote 
\begin{equation}\label{normal-form1}
\hat{v}_k=(1+|\nabla\gamma|^2)v_k,
\end{equation}
and we can reformulate \eqref{heatE-OBC} as 
\begin{multline}\label{heatE-OBC1}
 (\pd_t-\Delta)\,E^{-}_{\ep_1,\ep_2,\ep_3}(\hat{v}_k)(x,t) =
-2\sum_{i=1}^{d-1}\pd_{x_i}\pd_{x_d}\,E^{+}_{\ep_1,\ep_2,\ep_3}\left(\frac{\gamma'_i}{1+|\nabla\gamma|^2}\cdot\hat{v}_k\right) 
\\[3pt]
 +K_1 +K_{2}^{low}+K_3+K_4
+\left(\pd_t-\tsum_{i=1}^{d-1}\pd^2_{x_i}\right)E^{-}_{\ep_1,\ep_2,\ep_3}\left(\frac{|\nabla\gamma|^2}{1+|\nabla\gamma|^2}\cdot\hat{v}_k\right).
\end{multline}
Notice that the support of cut-off function $\zeta$ in spatial direction $\supp_{x} \zeta\subset \C_{\rho}$ and $\ep_1,\ep_2\ll\rho$. By using Lemma \ref{lem0} and Lemma \ref{lem-C-Z est} with $B(R)$ replaced by $\C(2\rho)$, we can get
\begin{align*}
&\|\nabla_x\, E^{-}_{\ep_1,\ep_2,\ep_3}(\hat{v}_k)\|_{\cW^{q,r}}
\lec \left(\ep_3^{-1}+\rho^{-1}\right)\|\bm{u}\|_{L^{q,r}(Q^+_1)}+\|p,\CF\|_{L^{q,r}(Q^+_1)}+\|\bm{f}\|_{L^{q_*,r}(Q^+_1)}\\[6pt]
&\quad\,+\sum_{i=1}^{d-1}\norm{\pd_{x_i}\,E^{+}_{\ep_1,\ep_2,\ep_3}\left(\frac{\gamma'_i}{1+|\nabla\gamma|^2}\cdot\hat{v}_k\right)}_{\cW^{q,r}}
+\norm{\nabla_{x'}\,E^{-}_{\ep_1,\ep_2,\ep_3}\left(\frac{|\nabla\gamma|^2}{1+|\nabla\gamma|^2}\cdot\hat{v}_k\right)}_{\cW^{q,r}}.
\end{align*}
Note that $\ep_3^{-1}$ comes from the term with $\pd_t$. Applying the commutator estimate \eqref{commu3} to the last two terms
in the RHS above and using \eqref{gamma} for the error terms, we can move the coefficients before $\hat{v}_k$ out of the norm,
\begin{align*}
\|\nabla_x\, E_{\ep_1,\ep_2,\ep_3}^{-}(\hat{v}_k)\|_{\cW^{q,r}}&\lec \left(\ep_3^{-1}+\rho^{-1}\right)\|\bm{u}\|_{L^{q,r}(Q^+_1)}+\|p,\CF\|_{L^{q,r}(Q^+_1)}+\|\bm{f}\|_{L^{q_*,r}(Q^+_1)} \notag\\[6pt]
&\quad\,+\|\nabla\gamma\|_{L^{\infty}(B_\rho^{\prime})}\cdot\|\nabla_{x'}\,E_{\ep_1,\ep_2,\ep_3}^{\pm}\left(\hat{v}_k\right)\|_{\cW^{q,r}}. 
\end{align*}
Since the last term in the RHS above has no vertical derivative, we can take $\ep_2 \to 0$ and get 
\begin{align*}
 \|\nabla_x\, E^{h}_{\ep_1,\ep_3}(\hat{v}_k)\|_{\cH^{q,r}}
&\lec \left(\ep_3^{-1}+\rho^{-1}\right)\|\bm{u}\|_{L^{q,r}(Q^+_1)}+\|p,\CF\|_{L^{q,r}(Q^+_1)}+\|\bm{f}\|_{L^{q_*,r}(Q^+_1)}\notag\\[6pt]
&\quad\,+\|\nabla\gamma\|_{L^{\infty}(B_\rho^{\prime})}\cdot\|\nabla_{x'}\,E^{h}_{\ep_1,\ep_3}\left(\hat{v}_k\right)\|_{\cH^{q,r}}.
\end{align*}
And then choose $\rho$ small enough in \eqref{rho2} to absorb the last  term above into the LHS,
\begin{equation*}
\|\nabla_x\, E^{h}_{\ep_1,\ep_3}(\hat{v}_k)\|_{\cH^{q,r}}\lec \left(\ep_3^{-1}+\rho^{-1}\right)\|\bm{u}\|_{L^{q,r}(Q^+_1)}+\|p,\CF\|_{L^{q,r}(Q^+_1)} +\|\bm{f}\|_{L^{q_*,r}(Q^+_1)},
\end{equation*}
which implies that, by taking $\ep_1\rightarrow 0$,
\begin{equation*}
\|\nabla E_{\ep_3}(\hat{v}_k)\|_{\cH^{q,r}}\lec \left(\ep_3^{-1}+\rho^{-1}\right)\|\bm{u}\|_{L^{q,r}(Q^+_1)}+\|p,\CF\|_{L^{q,r}(Q^+_1)} +\|\bm{f}\|_{L^{q_*,r}(Q^+_1)}.
\end{equation*}
Using $\|\nabla E_{\ep_3}({v}_k)\|_{\cH^{q,r}}\lec
\|\nabla E_{\ep_3}(\hat{v}_k)\|_{\cH^{q,r}}   +\norm{\gamma}_{C^{1,1}}\|v\|_{\cH^{q,r}}  $,
we have shown \eqref{claim1}.
\bigskip

\noindent
\textbf{\textit{Step 2: }} the gradient estimate.

\smallskip

To show \eqref{gra2}, we go back to the heat equation \eqref{heatE-OBC}, and use Lemma \ref{lem-C-Z est} to get
\begin{align*}
 &\|\nabla_x\, E^{-}_{\ep_1,\ep_2,\ep_3}(v_k)\|_{\cW^{q,r}}
\lec\rho^{-1}\|\bm{u}\|_{L^{q,r}(Q^+_1)}+\|p,\CF\|_{L^{q,r}(Q^+_1)}+\|\bm{f}\|_{L^{q_*,r}(Q^+_1)} \notag\\[6pt]
&\quad\ +\sum_{i=1}^{d-1}\|\pd_{x_i}\,E^{+}_{\ep_1,\ep_2,\ep_3}(\gamma'_i\cdot v_k)\|_{\cW^{q,r}}+\|\pd_{x_d}\,E^{-}_{\ep_1,\ep_2,\ep_3}(|\nabla\gamma|^2\cdot v_k)\|_{\cW^{q,r}}. 
\end{align*}
By using the commutator estimate \eqref{commu3} to the  second last term in the RHS above and using \eqref{gamma} for the error term, and then taking $\ep_1\rightarrow0$, we get
\begin{align*}
&\|\nabla_x\, E^{-}_{\ep_2,\ep_3}(v_k)\|_{\cW^{q,r}}
\lec\rho^{-1}\|\bm{u}\|_{L^{q,r}(Q^+_1)}+\|p,\CF\|_{L^{q,r}(Q^+_1)} +\|\bm{f}\|_{L^{q_*,r}(Q^+_1)}\notag\\[6pt]
&\quad\,+\|\nabla\gamma\|_{L^{\infty}(B'_{\rho})}\cdot\|\nabla_{x'}\,E^{+}_{\ep_2,\ep_3}(v_k)\|_{\cW^{q,r}}+\|\nabla\gamma\|_{L^{\infty}(B'_{\rho})}\cdot\|\pd_{x_d}\,E^{-}_{\ep_2,\ep_3}(v_k)\|_{\cW^{q,r}}. 
\end{align*}
The second last term in the RHS above does not blow up due to \eqref{claim1}, and the last term does not blow up due to the vertical mollification. By choosing $\rho$ small enough in \eqref{rho2}, the last term in the RHS can be absorbed to the LHS,
\begin{align}\label{gradEs-mol4}
\|\nabla_x\, E^{-}_{\ep_2,\ep_3}(v_k)\|_{\cW^{q,r}}
&\lec\rho^{-1}\|\bm{u}\|_{L^{q,r}(Q^+_1)}+\|p,\CF\|_{L^{q,r}(Q^+_1)} +\|\bm{f}\|_{L^{q_*,r}(Q^+_1)}\notag\\[6pt]
&\quad\,+\|\nabla\gamma\|_{L^{\infty}(B'_{\rho})}\cdot\|\nabla_{x'}\,E^{+}_{\ep_2,\ep_3}(v_k)\|_{\cW^{q,r}}.
\end{align}
Thanks to \eqref{claim1} and \eqref{rho2}, we get that by taking $\ep_2\to 0$, and by absorbing the last term by the LHS,
\begin{equation*}
 \|\nabla E_{\ep_3}(v_k)\|_{\cH^{q,r}}\lec\rho^{-1}\|\bm{u}\|_{L^{q,r}(Q^+_1)}+\|p,\CF\|_{L^{q,r}(Q^+_1)} +\|\bm{f}\|_{L^{q_*,r}(Q^+_1)}. 
\end{equation*}
By taking $\ep_3\to 0$, we get 
\begin{align}\label{0BCunif}
\|\nabla T^{-1}u_k\|_{L^{q,r}\left(\C^+_{\rho/2}\times(-\frac34,0)\right)}\lec\|\nabla v_k\|_{\cH^{q,r}}
\lec\rho^{-1}\|\bm{u}\|_{L^{q,r}(Q^+_1)}+\|p,\CF\|_{L^{q,r}(Q^+_1)}+\|\bm{f}\|_{L^{q_*,r}(Q^+_1)},
\end{align}
which together with \eqref{rho1}  implies that $\nabla\bm{u}\in L^{q,r}\big(B^+_{\rho/10}\times(-\frac34,0)\big)$, and
\begin{equation*}
\|\nabla \bm{u}\|_{L^{q,r}\big(B^+_{\rho/10}\times(-\frac34,0) \big)} \lec\rho^{-1}\|\bm{u}\|_{L^{q,r}(Q^+_1)}+\|p,\CF\|_{L^{q,r}(Q^+_1)}+\|\bm{f}\|_{L^{q_*,r}(Q^+_1)}.
\end{equation*}
Thus we arrive at \eqref{gra2}, which leads to \eqref{gra1}.

\smallskip

Finally, by using Lemma \ref{lem0} (i) and integration by parts, we get
\begin{equation}\label{get0BC}
\pd_{x_{d}}\,E_{\ep_2}^{-}(v_k)-E^{+}_{\ep_2}(\pd_{d} v_k)=2 v_k(x',0,t)\cdot \eta_{\ep_2}(x_d).
\end{equation}
Since the $L^{q,r}$ norm of the LHS of \eqref{get0BC} is uniformly bounded for all $\ep_2>0$ by  \eqref{gradEs-mol4} and \eqref{0BCunif}, we must have $v_k(x',0,t)=0$ and hence recover the no-slip boundary condition
\begin{equation}\label{no-slip1}
\bm{u}|_{\Sigma}=0.
\end{equation}
The divergence-free property $\div\bm{u}=0$  follows from \eqref{div-free weak}. Moreover, it is not difficult to deduce the weak formula \eqref{weakform2-0BC} from \eqref{weakform-0BC}, since $TE^{h,-}_{\ep_1,\ep_2}(T^{-1}\Phi)\in W_{c,\infty}^{2,1}\left(Q_1^+\cup\Sigma\right)$ converges to $\Phi$ for any $\Phi\in W_{c,\infty}^{1}\left(Q_1^+\cup\Sigma\right)$ with no-slip boundary condition $\Phi|_{\Sigma}=0$. Therefore, $(\bm{u},p)$ is a weak solution pair. This completes the proof of the first part of Theorem \ref{thm_0BC}.
\end{proof}

\subsection{Second derivative estimate}\label{sec3.2}

In this subsection we prove \eqref{gra0}. Since $\div\CF$ has the same position as $\bm{f}$ after integration by parts with the test function, for simplicity we assume $\CF=0$. We firstly show the $L^{q,r}$ boundedness of $\nabla^2 \bm{u}, \nabla p$ up to the boundary after time mollification using e.g.~\cite[Theorem IV.5.1 on page 276]{Gal11}. To this end, we
let $\rho_1=\rho/2$ and $\tau=2$ for the cut-off function $\zeta(x,t)$ defined in \eqref{def-zeta} and $v_k,\pi$ defined in \eqref{def-vpi}. By taking $\ep_1,\ep_2\rightarrow0$ in \eqref{heatE-OBC}, we have that for $1\leq k\leq d$,
\begin{equation}\label{heatE-OBC4}
-\Delta\,E_{\ep_3}(v_k)+\pd_{k} E_{\ep_3}\left(\pi\right)=-\pd_{t}E_{\ep_3}(v_k)+\bar{K}_1+\bar{K}_2+\bar{K}_3+\bar{K}_4
\end{equation}
in the distribution $\mathcal{D}'\left(\R_{+}^d\times(-1,0)\right)$, where $\bar{K}_1,\bar{K}_2$ are formulas similar to $K_1,K_2$ in \eqref{K1}-\eqref{K2} with only temporal mollification and no (even or odd) extension, and 
\begin{align*}
\bar{K}_3&=E_{\ep_3}\left(T^{-1}p\cdot\left(\pd_{k}-\gamma_{k}'\cdot\pd_{d}\right)\zeta\right)+\gamma_{k}'\cdot\pd_{d}E_{\ep_3}\left(\pi\right),\\[6pt]
\bar{K}_4&=E_{\ep_3}\left(T^{-1}f_k\cdot\zeta\right).
\end{align*}
Note that for $t$ close to $0^-$, we still use the argument of Remark \ref{ep0-argument}.
Moreover
\begin{equation}\label{div1}
\div E_{\ep_3}(\bm{v})=\bar{K}_5\eqdefa E_{\ep_3}\left(T^{-1}\tilde{\bm{u}}\cdot\nabla\zeta\right)+\sum_{i=1}^{d-1}\,\gamma_{i}'\cdot\pd_{d}E_{\ep_3}\left( v_i\right), 
\end{equation}
where $\tilde{{u}}$ is defined in \eqref{deftilde-u}.

We add a horizontal mollification $E^h_{\ep_1}$ on stationary Stokes equations \eqref{heatE-OBC4}-\eqref{div1}, apply $\nabla_{x'}$, rewrite $\nabla_{x'}E^h_{\ep_1} [\gamma_{k}'\pd_{d}E_{\ep_3}(\pi)]=\pd_{d}\nabla_{x'}E^h_{\ep_1} [\gamma_{k}' E_{\ep_3}(\pi)]$ and similarly for the second term of $\nabla_{x'}E^h_{\ep_1}\bar K_2^{high}$ (we do not rewrite its first term and let the $W^{-1,q}$ estimate absorb $\nabla_{x'}$, so that we do not need $\nabla^3\gamma\in L^\infty$),
and use \cite[Theorem IV.3.3 on page 257]{Gal11}, \eqref{gra2}, and the commutator estimate \eqref{commu3} to get
\begin{align*}
\|\nabla_{x'}\nabla E^{h}_{\ep_1,\ep_3}\left(\bm{v}\right),\nabla_{x'}E^{h}_{\ep_1,\ep_3}\left(\pi\right)\|_{\cH^{q,r}}\lec \left(\ep_3^{-1}+\rho^{-2}\right)\|\bm{u}\|_{L^{q,r}(Q^+_1)}+\rho^{-1}\|p, \bm{f}\|_{L^{q,r}(Q^+_1)}\notag\\[6pt]
+\|\nabla\gamma\|_{L^{\infty}(B^{\prime}_\rho)}\cdot\|\nabla_{x'}\nabla E^{h}_{\ep_1,\ep_3}\left(\bm{v}\right),\nabla_{x'}E^{h}_{\ep_1,\ep_3}\left(\pi\right)\|_{\cH^{q,r}},
\end{align*}
which together with \eqref{rho2} implies that
\begin{align*}
\|\nabla_{x'}\nabla E^{h}_{\ep_1,\ep_3}\left(\bm{v}\right),\nabla_{x'}E^{h}_{\ep_1,\ep_3}\left(\pi\right)\|_{\cH^{q,r}}\lec \left(\ep_3^{-1}+\rho^{-2}\right)\|\bm{u}\|_{L^{q,r}(Q^+_1)}+\rho^{-1}\|p,\bm{f}\|_{L^{q,r}(Q^+_1)}. 
\end{align*}
Then by taking $\ep_1\rightarrow0$, we achieve
\begin{equation}\label{hori}
\nabla'\,\nabla E_{\ep_3}(\bm{v}),\ \nabla'E_{\ep_3}\left(\pi\right)\in \cH^{q,r}. 
\end{equation}
In view of \eqref{heatE-OBC4} and \eqref{hori}, we see that
\begin{equation*}
-\pd_{d}^2E_{\ep_3}\left(T^{-1}\tilde{u}_d\cdot\zeta\right)+\pd_{d}E_{\ep_3}\left(\pi\right)\in \cH^{q,r}, 
\end{equation*}
where we let the last equation of \eqref{heatE-OBC4} minus the first $(d-1)$ equations times $\gamma'_k$.
It follows from the divergence-free property \eqref{div2} that
\begin{equation}\label{hori-1}
\pd_{d}^2E_{\ep_3}\left(T^{-1}\tilde{u}_d\cdot\zeta\right)=\pd_{d}E_{\ep_3}\left(T^{-1}\tilde{\bm{u}}\cdot\nabla \zeta\right)-\sum_{i=1}^{d-1}\,\pd_{d}\pd_i E_{\ep_3}\left(v_i\right)\in \cH^{q,r}. 
\end{equation}
Thus we obtain that
\begin{equation*}
\pd_{d}E_{\ep_3}\left(\pi\right)\in \cH^{q,r},
\end{equation*}
which together with \eqref{heatE-OBC4}, \eqref{hori} and \eqref{hori-1} implies that
\begin{equation}\label{hori2}
\nabla^2 E_{\ep_3}(\bm{v}),\ \nabla E_{\ep_3}\left(\pi\right)\in \cH^{q,r}. 
\end{equation}
Therefore
\begin{equation}\label{hori3}
\nabla^2 E_{\ep_3}\left(T^{-1}\bm{u}\cdot\psi_{2}(t)\right),\nabla E_{\ep_3}\left(T^{-1}p\cdot\psi_{2}(t)\right)\in L^{q,r}\left(\C_{\rho/2}^+\times(-1,0)\right).
\end{equation}

Now the non-stationary Stokes system \eqref{heatE-OBC4}-\eqref{div1} becomes well-defined (with all terms in $L^{q,r}$) instead of in distribution sense, and we will prove \eqref{gra3} 
using a similar proof of \cite{Ser00}.

\begin{rem}\label{rem:seregin-prop1}
We point out that, unlike \eqref{gra1}, the estimate \eqref{gra0} does not seem to follow directly from the standard local estimates for the heat equation and the Laplace equation alone. Let us explain this point heuristically in the flat-boundary case. Under the no-slip boundary condition, the pressure satisfies a Laplace equation, together with the boundary conditions $\pd_d p=\pd_d^2 u_d$ on the boundary. The local elliptic estimate for the pressure would therefore give, schematically,
$\norm{\nabla p}_{L^q(B_1^+)}\lesssim \norm{p,\nabla^2 u_d}_{L^q(B_2^+)}$, 
which still contains the higher-order term \(\nabla^2 u_d\) on the right-hand side. Thus the estimate cannot be closed by using only the local estimates for the heat and Laplace equations. In contrast, under the Navier boundary condition, the boundary condition for the pressure takes the form $\pd_d p= \alpha \pd_d u_d$ on the boundary, using the divergence-free condition. This yields instead $\norm{\nabla p}_{L^q(B_1^+)}\lesssim \norm{p,\nabla u_d}_{L^q(B_2^+)}$, 
which involves only a lower-order term. This is why, as will be seen in the proof of Theorem \ref{thm2_navierBC}, the argument in the Navier BC case can be reduced to the local estimates for the heat and Laplace equations. For the no-slip boundary condition, however, we use a different approach: we transform the local estimate into an initial-boundary value problem and then apply the classical estimate of Giga and Sohr \cite{Giga91}, following Seregin \cite[Proposition 1]{Ser00}. 
\end{rem}

We let $\frac{{\rho}}{2}\leq \rho_1<\rho<\rho_0\ll 1$ and $\tau=2$ for the cut-off function $\zeta(x,t)$ defined in \eqref{def-zeta}. Then  \eqref{hori2} is still valid due to \eqref{hori3}, if we pick a smaller $\rho$.

Consider the following boundary-value problem of stationary Stokes equations in a smooth subdomain $\C_{\omega}$ which satisfies $\C_{\rho}^+\subset \C_{\omega}\subset \C_{\rho_0}^+$,
\begin{equation}\label{div3}
-\Delta \bm{v}^{E}+\nabla \pi^E=0, \quad \div \bm{v}^{E}=\bar{K}_5,\quad\bm{v}^E|_{\pd\C_{\omega}}=0,
\end{equation}
where $\bar{K}_5$ is defined in \eqref{div1}. 
For each $t$, the compatibility condition
$
\int_{C_\omega}\bar K_5(x,t)\,dx=0
$
follows from the definition of $\bar K_5$ in \eqref{div1}.  
Then the stationary Stokes problem \eqref{div3} has a unique solution with $\int_{C_\omega}\pi^E(x,t)\,dx=0$.
By virtue of \cite[Theorem 2.3]{Sol77} (see also \cite[Theorem IV.6.1 on page 283]{Gal11} and \cite[Theorem 3.1]{FS09} combined) and \eqref{gra2}, we get 
\begin{align}\label{err1}
&\quad\ \|\nabla^2 \bm{v}^E,\nabla \pi^E\|_{L^{q,r}\left(\C_{\omega}\times(-1,0)\right)}\lec \|\bar{K}_5,\nabla \bar{K}_5\|_{L^{q,r}\left(\C_{\omega}\times(-1,0)\right)}  \notag\\[6pt]
&\lec \frac{1}{(\rho-\rho_1)^2} \|\bm{u},p,\bm{f}\|_{L^{q,r}(Q^+_1)}+\|\nabla\gamma\|_{L^{\infty}(B^{\prime}_{\rho})}\cdot\|\nabla\pd_{d}E_{\ep_3}\left(\bm{v}\right)\|_{L^{q,r}\left(\C_{\rho}^+\times(-1,0)\right)}. 
\end{align}
The constant in the above inequality may depend on  $\C_\omega$. Nevertheless, we can choose small $\rho< \rho_0\ll1$ to absorb the last term in the end.

We still need the estimate of $\pd_t \bm{v}^E$.
Since $\bar K_5$ is temporally mollified, it is smooth in $t$ for fixed $\varepsilon_3>0$. The same is therefore true for $(\bm{v}^E,\pi^E)$. Indeed, by linearity, for $h\neq0$,
\[
\frac{(\bm{v}^E,\pi^E)(t+h)-(\bm{v}^E,\pi^E)(t)}{h}
\]
solves the same stationary Stokes problem with divergence datum
$(\bar K_5(t+h)-\bar K_5(t))/h$.
Letting $h\to0$ and using the continuous dependence estimate above, we obtain that $\partial_t(\bm{v}^E,\pi^E)$ exists and is the solution of the same stationary Stokes problem with divergence datum $\partial_t\bar K_5$.
Using \cite[Theorem 2.4]{Sol77}, we have that if $\pd_t\bar{K}_5=\div\bm{K}_{51}+K_{52}$, then
\begin{align}\label{Solon77thm24}
\|\pd_t\bm{v}^E(\cdot,t)\|_{L^q(\C_{\omega})}\lesssim \|\bm{K}_{51}\|_{L^q(\C_{\omega})}+\|\bm{K}_{51}\cdot \bm{n}\|_{L^q(\pd\C_{\omega})}+\|K_{52}\|_{L^q(\C_{\omega})}. 
\end{align}
Notice that 
\begin{equation}\label{err3}
\pd_{t}\bar{K}_5=\pd_{t}E_{\ep_3}\left(T^{-1}\tilde{\bm{u}}\cdot\nabla\zeta\right)+\sum_{i=1}^{d-1}\,\gamma_{i}'\cdot\pd_{d}\pd_{t}E_{\ep_3}\left( v_i\right). 
\end{equation}
The second term above will generate the last term in the RHS of \eqref{err2}.
The tricky term is the first term above. To deal with it, we replace the cut-off function $\zeta$ with $\nabla \zeta$ in the definition of $v,\pi$ in \eqref{heatE-OBC4}
and denote
\[
\hat v_k = (T^{-1}u_k) \pd_k \zeta, \quad 
v_k^\sharp = (T^{-1}u_k) \pd_d \zeta, \quad 
\pi _k =  (T^{-1}p )\pd_k \zeta, \quad 
\bar K_{jk}^i=  \bar K_{jk}|_{\zeta \to \pd_i \zeta},
\]
where $\bar K_{jk}$ denotes the $k$-th component of $\bar K_j$ for $j\le4$.
By \eqref{deftilde-u} and \eqref{heatE-OBC4},
\begin{align}
\pd_{t}E_{\ep_3}\left(T^{-1}\tilde{\bm{u}}\cdot\nabla\zeta\right)
&= \pd_{t}E_{\ep_3}\left(\tsum_{k=1}^d \hat v_k - \tsum_{i<d} \gamma_i' v_i^\sharp  \right)\notag
\\
\label{pd_tu_decomp}
&=\tsum_{k=1}^d \bke{ \Delta E_{\ep_3} \hat v_k - \pd _k E_{\ep_3} \pi_k + \tsum_{j=1}^4 \bar K_{jk}^k}
\\
&\quad - \tsum_{i<d}  \gamma_i'\bke{ \Delta E_{\ep_3}  v_i^\sharp - \pd _i E_{\ep_3} \pi_d+ \tsum_{j=1}^4 \bar K_{ji}^d}.\notag
\end{align}
On the right side, the higher order terms are two derivatives of $v$ and  one derivative of $\pi$, coming from $\Delta v$, $\nabla \pi$, $K_2^{high}$ and $K_3$. All the remaining terms are of lower order  and included in $K_{52}$. For the higher order terms with $\pd_j$ where $j<d$, they will generate no boundary term $\bm{K}_{51}\cdot \bm{n}$. Thus (also using $\pd_d \zeta|_{x_d=0}=0$),
\begin{multline}\label{treat-boundary-t}
\pd_{t}E_{\ep_3}\left(T^{-1}\tilde{\bm{u}}\cdot\nabla\zeta\right)
\\
=\sum_{i=1}^{d-1}\left\{(1+|\nabla\gamma|^2)\cdot\pd_{d}^2 E_{\ep_3}\left(T^{-1}u_i\cdot\pd_{i}\zeta\right) +\gamma_{i}'\cdot\pd_{d}E_{\ep_3}\left(T^{-1}p\cdot\pd_{i}\zeta\right)\right\}+\cdots,
\end{multline}
where the first two terms on the right side are the main terms that generate boundary term $\bm{K}_{51}\cdot \bm{n}$ in the RHS of \eqref{Solon77thm24} due to non-zero boundary condition for $\pd_{d} T^{-1} u_i|_{x_d=0}$ and $T^{-1}p|_{x_d=0}$. They come from the second row of \eqref{pd_tu_decomp}. Similar terms in the last row of \eqref{pd_tu_decomp} do not generate boundary terms since $\pd_d\zeta|_{x_d=0}=0$.
The remaining terms only generate the first and third kinds of terms in the RHS of \eqref{Solon77thm24}. For the terms with non-zero boundary condition, we use the following estimate: for $g\in C^1_c(\C_{\rho}^+\cup B'_{\rho})$,  (see the estimate after \cite[(6.7)]{CLT24})
\begin{align}\label{bad1}
\frac{1}{\rho-\rho_1}\|g|_{x_d=0}\|_{L^{q}\left(B_{\rho}'\right)}&\le \frac{c}{\rho-\rho_1} \|g\|_{L^{q}\left(\C_{{\rho}}^+\right)}^{1-\frac{1}{q}}\cdot \|\pd_{d}g\|_{L^{q}\left(\C_{{\rho}}^+\right)}^{\frac{1}{q}}\notag\\[6pt]
&\le \frac{c_\ep}{(\rho-\rho_1)^{\frac{q}{q-1}}} \|g\|_{L^{q}\left(\C_{{\rho}}^+\right)}+\ep\|\pd_{d}g\|_{L^{q}\left(\C_{{\rho}}^+\right)}, 
\end{align}
where $\ep>0$ is a small constant to be chosen. By \eqref{gra2}, \eqref{Solon77thm24}, \eqref{err3} and \eqref{bad1} for the bad terms ($g=\pd_{d} E_{\ep_3}\left(T^{-1}u_i\cdot \psi_2\right)$ and $g=E_{\ep_3}\left(T^{-1}p\cdot \psi_2\right)$) in \eqref{treat-boundary-t}, we obtain that
\begin{align}\label{err2}
&\|\pd_{t} \bm{v}^E\|_{L^{q,r}\left(\C_{\omega}\times(-1,0)\right)}\lec\frac{c_\ep}{(\rho-\rho_1)^{3+\frac{1}{q-1}}}\|\bm{u},p,\bm{f}\|_{L^{q,r}(Q^+_1)}\notag\\[6pt]
&\quad \ + \ep\|\pd_{d}^2E_{\ep_3}\left(T^{-1}\bm{u}\cdot\psi_{2}(t)\right),\pd_{d}E_{\ep_3}\left(T^{-1}p\cdot\psi_{2}(t)\right)\|_{L^{q,r}\left(\C_{\rho}^+\times(-1,0)\right)}\notag\\[6pt]
&\quad\ +\|\nabla\gamma\|_{L^{\infty}(B^{\prime}_{\rho})}\cdot\|\pd_{t}E_{\ep_3}\left(\bm{v}\right)\|_{L^{q,r}\left(\C_{\rho}^+\times(-1,0)\right)}.
\end{align}

In view of \eqref{heatE-OBC4}-\eqref{div1} and \eqref{div3}, we see that the remainder $\bm{v}^R=E_{\ep_3}\left(\bm{v}\right)-\bm{v}^E$ and $\pi^R=E_{\ep_3}\left(\pi\right)-\pi^E$ satisfy the non-stationary Stokes equations in $\C_\omega\times(-1,0)$
\begin{equation}\label{RSE}
\left\{
\begin{array}{l}
\pd_{t}v_k^R-\Delta v_k^R+\pd_{k} \pi^R=-\pd_{t} v_k^E+\bar{K}_1+\bar{K}_2+\bar{K}_3+\bar{K}_4, \quad 1\leq k \leq d,  \\[6pt]
\div \bm{v}^R=0,\\[6pt]
\bm{v}^R|_{\pd \C_{\omega}}=0,\quad\bm{v}^R|_{t=-1}=0.
\end{array}\right.
\end{equation}
Applying \cite[Theorem 2.8]{Giga91} to the non-stationary Stokes equations \eqref{RSE} in  $C_\omega$, using \eqref{gra2},  \eqref{rho2} and the $v^E,\pi^E$ estimates \eqref{err1}, \eqref{err2} with sufficiently small $\ep>0$, we get that for $\rho_1<\rho<\rho_0\ll 1$,
\begin{equation}\label{add4}
\cA(\rho_1)
\leq\frac{C}{(\rho-\rho_1)^{3+\frac{1}{q-1}}} \|\bm{u},p,\bm{f}\|_{L^{q,r}(Q^+_1)}+\frac12\cA(\rho),
\end{equation}
where
\begin{equation*}
\cA(R)
\eqdefa \|\pd_{t}E_{\ep_3}\left(T^{-1}\bm{u}\cdot\psi_{2}(t)\right),\nabla^2 E_{\ep_3}\left(T^{-1}\bm{u}\cdot\psi_{2}(t)\right),\nabla E_{\ep_3}\left(T^{-1}p\cdot\psi_{2}(t)\right)\|_{L^{q,r}\left(\C_{R}^{+}\times(-1,0)\right)}.
\end{equation*}
By using Lemma \ref{lemA2}, we achieve
\begin{equation}\label{err5}
\cA(\rho/2)\lec \rho^{-3-\frac{1}{q-1}}\|\bm{u},p,\bm{f}\|_{L^{q,r}(Q^+_1)}.
\end{equation}
Then, by taking $\ep_3\rightarrow0$ in \eqref{err5}, the three terms admit weak limits in $L^{q,r}$ (up to subsequence). The limit of the first term satisfies the weak formulation for the time derivative, and
\begin{align}\label{add5}
&\quad\ \|\pd_{t}T^{-1}\bm{u},\nabla^2 T^{-1}\bm{u},\nabla T^{-1}p\|_{L^{q,r}\left(\C_{\rho/2}^{+}\times(-\frac12,0)\right)}\lec\rho^{-3-\frac{1}{q-1}}\|\bm{u},p,\bm{f}\|_{L^{q,r}(Q^+_1)}.
\end{align}
This gives rise to \eqref{gra3} and completes the proof of Theorem \ref{thm_0BC}.

\section{Proof of Theorem \ref{thm_navierBC}}\label{sec4}

In this section we give the proof of Theorem \ref{thm_navierBC}.
Analogous to the proof of Theorem \ref{thm_0BC}, 
with interior estimate \eqref{interior1} and covering lemma, we only need to prove the local boundary estimate, i.e., \eqref{limit6} in Lemma \ref{lem4-NavierBC}. The idea is to choose suitable test function $\Phi$ in the definition of weak solution \eqref{weakform-NavierBC} to get the heat equation of the approximation of the vorticity (denoted as $E^+_{\ep_1,\ep_2,\ep_3}(\omega_{ij})$ and $E^-_{\ep_1,\ep_2,\ep_3}(\omega_{id})$). Similar to Step 1 in the proof of \eqref{gra1}, in this section, lots of work needs to be done to raise up the regularity of $\bm{u}$.

Recall $\tilde{\bm{u}}, \zeta$ defined in \eqref{deftilde-u} and \eqref{def-zeta}.  Denote the vector $\bm{v}=(v_1,\cdots,v_d)$ as (different from \eqref{def-vpi})
\begin{equation}\label{defv1}
\bm{v}(x,t)\eqdefa T^{-1}\tilde{\bm{u}}(x,t)\cdot \zeta(x,t),
\end{equation}
where we choose $\tau=1$ in $\zeta$, and $\rho_1,\rho$ will be specified later. Denote vector
\begin{equation}\label{def-v}
E^{\#}_{\ep_1,\ep_2,\ep_3}(\bm{v})\eqdefa \Big(E^+_{\ep_1,\ep_2,\ep_3}(v_1),\cdots,E^+_{\ep_1,\ep_2,\ep_3}(v_{d-1}),E^-_{\ep_1,\ep_2,\ep_3}(v_d)\Big),
\end{equation}
which spatially is compactly supported in $\R^d$. In our previous work \cite[Theorem 1.1]{CLT24} where we prove the gradient estimate near a flat boundary, we use a similar extension while there is some weight function $e^{\alpha y_3/2}$, which makes the proof easier and the proof could work for very weak solutions. Here near curved boundary, it's hard to find such a weight function, so our proof becomes more complicated and it only works for weak solutions.

We define the curl matrix $\bm{\omega}=\curl\bm{v}$ of $\bm{v}$ by
\begin{equation}
\omega_{ij}\eqdefa \pd_i v_j-\pd_j v_i,\quad\omega_{ij}=-\omega_{ji}.
\end{equation}
Note that $\bm{\omega}$ is well-defined since we assume $\nabla \bm{u}\in L^1(Q^+_1)$ in the definition of weak solutions. 
The curl matrix of the vector $E^{\#}_{\ep_1,\ep_2,\ep_3}(\bm{v})$ is
\begin{equation*}
\Big(\curl E^{\#}_{\ep_1,\ep_2,\ep_3}(\bm{v})\Big)_{ij} = \pd_i\left(E^{\#}_{\ep_1,\ep_2,\ep_3}(\bm{v})\right)_j-\pd_j\left(E^{\#}_{\ep_1,\ep_2,\ep_3}(\bm{v})\right)_i
\end{equation*}
for $1\leq i,j\leq d$. Moreover, we have that
\begin{equation}\label{E-sharp-v}
\Big(\curl E^{\#}_{\ep_1,\ep_2,\ep_3}(\bm{v})\Big)_{ij}=
\left \{
\begin{split}
E^+_{\ep_1,\ep_2,\ep_3}(\omega_{ij}), \qquad &  i<j\leq d-1,
\\[4pt]
E^-_{\ep_1,\ep_2,\ep_3}(\omega_{id}), \qquad & 1\leq i<j=d.
\end{split}
\right.
\end{equation}

By virtue of  \eqref{div2}, we have
\begin{equation*}
\div E^{\#}_{\ep_1,\ep_2,\ep_3}(\bm{v})=E^+_{\ep_1,\ep_2,\ep_3}\left(\div \bm{v}\right)=E^+_{\ep_1,\ep_2,\ep_3}\left(\,  T^{-1}\tilde{\bm{u}}\cdot \nabla\,\zeta\right).
\end{equation*}
Notice that
\begin{align}\label{lap-Ev}
\Delta \left(E^{\#}_{\ep_1,\ep_2,\ep_3}(\bm{v})\right)_{i}=\pd_i\left(\div E^{\#}_{\ep_1,\ep_2,\ep_3}(\bm{v})\right)-\sum_{j=1}^{d}\pd_j\left(\curl E^{\#}_{\ep_1,\ep_2,\ep_3}(\bm{v})\right)_{ij}.
\end{align}
Hence using Lemma \ref{lem-C-Z est}, we have 
\begin{equation}\label{divcurl1}
\|\nabla E^{\#}_{\ep_1,\ep_2,\ep_3}(\bm{v})\|_{\cW^{q,r}}\lec \|\div E^{\#}_{\ep_1,\ep_2,\ep_3}(\bm{v})\|_{\cW^{q,r}}+\|\curl E^{\#}_{\ep_1,\ep_2,\ep_3}(\bm{v})\|_{\cW^{q,r}},
\end{equation}
while it's easy to get
\begin{equation}\label{divcontrol}
\|\div E^{\#}_{\ep_1,\ep_2,\ep_3}(\bm{v})\|_{\cW^{q,r}}\lec \frac{1}{\rho-\rho_1}\|\bm{u}\|_{L^{q,r}\left(Q^+_1\right)}.
\end{equation}

In what follows we shall deduce the heat equation of  $\curl\left(E^{\#}_{\ep_1,\ep_2,\ep_3}(\bm{v})\right)_{ij}$. There will be some inherent  cancellations in the calculation.

From the form of weak solution \eqref{weakform-NavierBC}, by using $\det{(T^{-1})}=1$, $\DD\bm{u}:\DD\Phi=\DD\bm{u}:\nabla\Phi$, and $\Phi=\Phi_{tan}$ on the boundary, we get
\begin{equation}\label{Stokes-curved}
\begin{split}
&\int_{\R^d_{+}\times\R} T^{-1}\left(-\bm{u} \cdot \partial_s\Phi+2\DD\bm{u}:\nabla\Phi \right)\dy\ds +\int_{\R^{d-1}\times \R}\alpha T^{-1} (\bm{u}\cdot\Phi) \sqrt{1+|\nabla \gamma|^2}\dy'\ds\\[6pt]
&=\int_{\R^d_{+}\times\R}T^{-1}\left( \bm{f}\cdot\Phi-\CF:\nabla\Phi \right)\dy\ds,
\end{split}
\end{equation}
while using \eqref{commu2} the second term in the first integral above can be rewritten as 
\begin{equation}\label{Stokes-curved-term2}
\begin{split}
&\,\quad T^{-1}(2\DD\bm{u}:\nabla\Phi)\\
&=\nabla_y\, T^{-1}\bm{u}: \nabla_y\, T^{-1} \Phi-\sum_{k=1}^{d-1}\gamma_{k}'\cdot\left(\pd_{y_k}\,T^{-1}\bm{u}\cdot \pd_{y_d}\,T^{-1}\Phi+\pd_{y_{d}}\,T^{-1}\bm{u}\cdot\pd_{y_{k}}\,T^{-1}\Phi\right)\\
&\,\quad+|\nabla\gamma|^2\cdot\pd_{y_d}\,T^{-1}\bm{u}\cdot \pd_{y_d}\,T^{-1}\Phi+\sum_{1\leq k,\ell\leq d}T^{-1}\pd_{y_\ell}u_{k}\cdot T^{-1}\pd_{y_k}\Phi_{\ell}.
\end{split}
\end{equation}
Above, the last sum is $T^{-1}((\nabla\bm{u})^T:\nabla\Phi)$ and was not expanded using \eqref{commu2}.

\subsection{Heat equations of the approximations of vorticity}

We now rewrite the heat equations satisfied by $\curl E^\sharp_{\ep_1,\ep_2,\ep_3}(\bm{v})$, the approximations of the vorticity matrix.

\smallskip

\noindent$\bullet$ Case 1: the heat equations of $E^+_{\ep_1,\ep_2,\ep_3}(\omega_{ij})$ for $1\leq i<j\leq d-1$.

Notice that
\begin{align}\label{heatE-I_1}
&\ \quad(\pd_t-\Delta)E^+_{\ep_1,\ep_2,\ep_3}(\omega_{ij})\notag\\[6pt]
&=-\int_{\R^d_{+}\times\R}\left(T^{-1}u_i\cdot\zeta\cdot \pd_{y_j}-T^{-1}u_j\cdot\zeta\cdot \pd_{y_i}\right)(\pd_s+\Delta_y)\,\eta^+_{\ep_1,\ep_2,\ep_3}(x,y,t-s)\dy\ds\notag\\[6pt]
&\eqdefa I_1+I_2, 
\end{align}
where
\begin{equation}\label{I1}
I_1\eqdefa-\int_{\R^d_{+}\times\R}\left(T^{-1}u_i\cdot \pd_{y_j}-T^{-1}u_j\cdot \pd_{y_i}\right)(\pd_s+\Delta_y)\,\xi^+(x,y,t,s)\dy\ds,
\end{equation}
and
\begin{align}\label{I2}
I_2&\eqdefa\int_{\R^d_{+}\times\R}\left(T^{-1}u_i\cdot \pd_{y_j}-T^{-1}u_j\cdot \pd_{y_i}\right)\Big((\pd_s+\Delta_y)\,\zeta\cdot\eta^+_{\ep_1,\ep_2,\ep_3}(x,y,t-s)\Big)\dy\ds\notag\\[6pt]
&\,\quad+2\int_{\R^d_{+}\times\R}\left(T^{-1}u_i\cdot \pd_{y_j}-T^{-1}u_j\cdot \pd_{y_i}\right)\Big(\nabla_y\zeta\cdot\nabla_y\eta^+_{\ep_1,\ep_2,\ep_3}(x,y,t-s)\Big)\dy\ds\notag\\[6pt]
&\quad\ -(\pd_{t}-\Delta_x)\int_{\R^d_{+}\times\R}\bke{\left(T^{-1}u_i\cdot \pd_{y_j}-T^{-1}u_j\cdot \pd_{y_i}\right)\zeta}\cdot\eta^+_{\ep_1,\ep_2,\ep_3}(x,y,t-s)\dy\ds. 
\end{align}
Now we choose the following test function $\Phi$ (with $(x,t)$ as parameters), which satisfies the div-free property and boundary condition in the definition of weak solution \eqref{weakform-NavierBC} in $Q_1^+$ with normal direction $(\nabla' \gamma,-1)$, given by
\begin{equation*}
T^{-1}\Phi(y,s)=(0,\cdots,0,\underbrace{\pd_{y_j}\xi^+}_{i-th},0,\cdots,0,\underbrace{-\pd_{y_i}\xi^+}_{j-th},0,\cdots,0,\underbrace{-\gamma_{j}'\cdot \pd_{y_{i}}\xi^++\gamma_{i}'\cdot \pd_{y_{j}}\xi^+}_{d-th}).
\end{equation*}
Note that for $t$ close to $0^-$ we use the argument in Remark \ref{ep0-argument}. 
Denote the perturbation part
\[T^{-1}\Phi^{pert}\eqdefa(0,\cdots,0,\underbrace{-\gamma_{j}'\cdot \pd_{y_{i}}\xi^++\gamma_{i}'\cdot \pd_{y_{j}}\xi^+}_{d-th}),\]
then we see that 
\begin{align}\label{I_1-expression}
I_1=\int_{\R^d_{+}\times\R} \left(-T^{-1}\bm{u} \cdot \partial_sT^{-1}(\Phi-\Phi^{pert})+\nabla_y\, T^{-1}\bm{u}: \nabla_y\, T^{-1} (\Phi-\Phi^{pert})\right)\dy\ds, 
\end{align}
since $\pd_{d}\xi^+|_{y_d=0}=0$ (cf.~$\pd_{d}\xi^-|_{y_d=0}\not=0$) and $i,j<d$. Then, since $\Phi$ satisfies \eqref{Stokes-curved} and \eqref{Stokes-curved-term2}, we get
\begin{equation*}
I_1= \tsum_{k=1}^{8}I_{1k},
\end{equation*}
where 
\begin{equation*}
I_{11}\eqdefa\int_{\R^d_{+}\times\R}T^{-1}u_d\cdot\pd_s \Big(-\gamma_{j}'\cdot \pd_{y_{i}}\,\xi^++\gamma_{i}'\cdot \pd_{y_{j}}\,\xi^+\Big)\dy\ds,
\end{equation*}
\begin{align*}
I_{12}&\eqdefa-\int_{\R^d_{+}\times\R}\nabla_y\, T^{-1}u_d\cdot\nabla_y\, \Big(-\gamma_{j}'\cdot \pd_{y_{i}}\,\xi^++\gamma_{i}'\cdot \pd_{y_{j}}\,\xi^+\Big)\dy\ds\\[6pt]
&=-\sum_{k=1}^{d-1}\int_{\R^d_{+}\times\R}\pd_{y_k}\, T^{-1}u_d\cdot \Big(-\gamma_{jk}''\cdot \pd_{y_{i}}\,\xi^++\gamma_{ik}''\cdot \pd_{y_{j}}\,\xi^+\Big)\dy\ds\\[6pt]
&\quad\ +\sum_{k=1}^{d-1}\int_{\R^d_{+}\times\R} T^{-1}u_d\cdot \pd_{y_k}\,\Big(-\gamma_{j}'\cdot\pd_{y_k} \pd_{y_{i}}\,\xi^++\gamma_{i}'\cdot\pd_{y_k} \pd_{y_{j}}\,\xi^+\Big)\dy\ds\\[6pt]
&\quad\ +\int_{\R^d_{+}\times\R} T^{-1}u_d\cdot \Big(-\gamma_{j}'\cdot \pd_{y_d}^2\pd_{y_{i}}\,\xi^++\gamma_{i}'\cdot \pd_{y_d}^2 \pd_{y_{j}}\,\xi^+\Big)\dy\ds,
\end{align*}
We rewrite $I_{12}$ to avoid higher order terms with $\nabla'  T^{-1}\bm{u}$ in the integral (since we 
can only treat higher order terms with $\pd_d  T^{-1}\bm{u}$ in the integral; see $H_5$ in step 1 of Lemma \ref{lem1-navierBC}), and to avoid third derivative of $\gamma$ (since we only assume $\pd \Omega\in C^{1,1}$).

\begin{align*}
I_{13}&\eqdefa\,\sum_{k=1}^{d-1}\int_{\R^d_{+}\times\R}\gamma_{k}'\cdot\Big(\pd_{y_k}\,T^{-1}u_i\cdot \pd_{y_{j}}\,\pd_{y_{d}}\xi^{+}-\pd_{y_k}T^{-1}u_j\cdot \pd_{y_{i}}\pd_{y_{d}}\xi^{+}\Big)\dy\ds \\[6pt]
&\quad\ +\sum_{k=1}^{d-1}\int_{\R^d_{+}\times\R}\gamma_{k}'\cdot\Big(\pd_{y_k}\,T^{-1}u_d\cdot\pd_{y_{d}} \left(-\gamma_{j}'\cdot\pd_{y_{i}}\xi^{+}+\gamma_{i}'\cdot\pd_{y_{j}}\xi^{+}\right)\Big)\dy\ds\\[6pt]
&=-\sum_{k=1}^{d-1}\int_{\R^d_{+}\times\R}T^{-1}u_i\cdot \pd_{y_k}\left(\gamma_{k}'\cdot\pd_{y_{j}}\,\pd_{y_{d}}\xi^{+}\right)-T^{-1}u_j\cdot\pd_{y_k} \left(\gamma_{k}'\cdot\pd_{y_{i}}\pd_{y_{d}}\xi^{+}\right)\ds \\[6pt]
&\quad\ -\sum_{k=1}^{d-1}\int_{\R^d_{+}\times\R}T^{-1}u_d\cdot\pd_{y_{d}} \pd_{y_k}\left(-\gamma_{j}'\gamma_{k}'\cdot\pd_{y_{i}}\xi^{+}+\gamma_{i}'\gamma_{k}'\cdot\pd_{y_{j}}\xi^{+}\right)\dy\ds,
\end{align*}
\begin{align*}
I_{14}&\eqdefa\,\sum_{k=1}^{d-1}\int_{\R^d_{+}\times\R}\gamma_{k}'\cdot\Big(\pd_{y_d}\,T^{-1}u_i\cdot \pd_{y_{j}}\,\pd_{y_{k}}\xi^{+}-\pd_{y_d}T^{-1}u_j\cdot \pd_{y_{i}}\pd_{y_{k}}\xi^{+}\Big)\dy\ds \\[6pt]
&\quad\ +\sum_{k=1}^{d-1}\int_{\R^d_{+}\times\R}\gamma_{k}'\cdot\Big(\pd_{y_d}\,T^{-1}u_d\cdot\pd_{y_{k}} \left(-\gamma_{j}'\cdot\pd_{y_{i}}\xi^{+}+\gamma_{i}'\cdot\pd_{y_{j}}\xi^{+}\right)\Big)\dy\ds,
\end{align*}
\begin{align*}
I_{15}&\eqdefa\,-\int_{\R^d_{+}\times\R}|\nabla\gamma|^2\cdot\Big(\pd_{y_d}\,T^{-1}u_i\cdot \pd_{y_{j}}\,\pd_{y_{d}}\xi^{+}-\pd_{y_d}T^{-1}u_j\cdot \pd_{y_{i}}\pd_{y_{d}}\xi^{+}\Big)\dy\ds \\[6pt]
&\quad\ -\int_{\R^d_{+}\times\R}|\nabla\gamma|^2\cdot\Big(\pd_{y_d}\,T^{-1}u_d\cdot\pd_{y_{d}} \left(-\gamma_{j}'\cdot\pd_{y_{i}}\xi^{+}+\gamma_{i}'\cdot\pd_{y_{j}}\xi^{+}\right)\Big)\dy\ds,
\end{align*}
\begin{align*}
I_{16}&\eqdefa-\sum_{1\leq k,\ell\leq d}\int_{\R^d_{+}\times\R}T^{-1}\pd_{y_\ell}u_{k}\cdot T^{-1}\pd_{y_k}\Phi_{\ell}\dy\ds   \\[6pt]
&=-\sum_{\ell=i,j,d}\sum_{1\leq k\leq d}\int_{\R^d_{+}\times\R}\left(\pd_{y_\ell}-\gamma_{\ell}'\cdot\pd_{y_d}\right)T^{-1}u_{k}\cdot \left(\pd_{y_k}-\gamma_{k}'\cdot\pd_{y_d}\right)T^{-1}\Phi_{\ell}\dy\ds\\[6pt]
&=-\sum_{1\leq k\leq d-1}\int_{\R^d_{+}\times\R}\pd_{y_d}\Big(T^{-1}u_{k}\cdot \left(-\gamma_{jk}''\cdot \pd_{y_{i}}+\gamma_{ik}''\cdot \pd_{y_{j}}\right)\xi^+\Big)\dy\ds, 
\end{align*}
\begin{align*}
I_{17}&\eqdefa-\int_{\pd \R^d_{+}\times\R}T^{-1}\left(\alpha \bm{u}\cdot\Phi\right)\cdot \sqrt{1+|\nabla \gamma|^2}\dy'\ds  \\[6pt]
&=\alpha \int_{\R^d_{+}\times\R}\pd_{y_d}\Big(T^{-1}u_{i}\cdot \pd_{y_j}\xi^+-T^{-1}u_{j}\cdot\pd_{y_i}\xi^+\Big)\cdot \sqrt{1+|\nabla \gamma|^2}\dy\ds\\[6pt]
&\quad\ +\alpha\int_{\R^d_{+}\times\R}\pd_{y_d}\Big(T^{-1}u_d\cdot\left(-\gamma_{j}'\cdot \pd_{y_{i}}+\gamma_{i}'\cdot \pd_{y_{j}}\right)\xi^+\Big)\cdot \sqrt{1+|\nabla \gamma|^2}\dy\ds, 
\end{align*}
\begin{align*}
I_{18}&\eqdefa\int_{\R^d_{+}\times\R}T^{-1}\left( \bm{f}\cdot\Phi-\CF:\nabla\Phi \right)\dy\ds  \\[6pt]
&=\int_{\R^d_{+}\times\R}T^{-1} f_i\cdot\pd_{y_j}\xi^+-T^{-1}f_j\cdot\pd_{y_i}\xi^++T^{-1}f_d\cdot(-\gamma_{j}'\cdot \pd_{y_{i}}+\gamma_{i}'\cdot \pd_{y_{j}})\xi^+ \dy\ds\\[6pt]
&\quad\ - \sum_{1\leq\ell\leq d}\,\int_{\R^d_{+}\times\R}T^{-1}\CF_{i\ell}\cdot\pd_{y_\ell}\pd_{y_j}\xi^+-T^{-1}\CF_{j\ell}\cdot\pd_{y_\ell}\pd_{y_i}\xi^+\dy\ds\\[6pt]
&\quad\ -\sum_{1\leq\ell\leq d}\,\int_{\R^d_{+}\times\R}T^{-1}\CF_{d\ell}\cdot\pd_{y_\ell}(-\gamma_{j}'\cdot \pd_{y_{i}}+\gamma_{i}'\cdot \pd_{y_{j}})\xi^+\dy\ds.
\end{align*}

\smallskip

\noindent$\bullet$ Case 2: the heat equations of $E^-_{\ep_1,\ep_2,\ep_3}(\omega_{id})$ for $1\leq i\leq d-1$.

Analogously, we have that for $1\leq i<j=d$,
\begin{equation}\label{heatE-I_3}
\begin{split}
&\ \quad(\pd_t-\Delta)E^-_{\ep_1,\ep_2,\ep_3}(\omega_{id})\\[6pt]
&=-\int_{\R^d_{+}\times\R}\left(T^{-1}u_i\cdot\zeta\cdot \pd_{y_d}-T^{-1}\tilde{u}_d\cdot\zeta\cdot\pd_{y_i}\right)(\pd_s+\Delta_y)\,\eta^-_{\ep_1,\ep_2,\ep_3}(x,y,t-s)\dy\ds\\[6pt]
&\eqdefa I_3+I_4,
\end{split}
\end{equation}
where
\begin{equation*}
I_3\eqdefa-\int_{\R^d_{+}\times\R}\left(T^{-1}u_i\, \pd_{y_d}-T^{-1}\tilde{u}_d\, \pd_{y_i}\right)(\pd_s+\Delta_y)\,\xi^-(x,y,t,s)\dy\ds,
\end{equation*}
\begin{align*}
I_4&\eqdefa\int_{\R^d_{+}\times\R}\left(T^{-1}u_i\, \pd_{y_d}-T^{-1}\tilde{u}_d \,\pd_{y_i}\right)\Big((\pd_s+\Delta_y)\,\zeta\cdot\eta^-_{\ep_1,\ep_2,\ep_3}(x,y,t-s)\Big)\dy\ds\\[6pt]
&\,\quad+2\int_{\R^d_{+}\times\R}\left(T^{-1}u_i\, \pd_{y_d}-T^{-1}\tilde{u}_d \,\pd_{y_i}\right)\Big(\nabla_y\zeta\cdot\nabla_y\eta^-_{\ep_1,\ep_2,\ep_3}(x,y,t-s)\Big)\dy\ds\\[6pt]
&\,\quad-(\pd_{t}-\Delta)\int_{\R^d_{+}\times\R}\bke{\left(T^{-1}u_i\, \pd_{y_d}-T^{-1}\tilde{u}_d \,\pd_{y_i}\right)\zeta}\cdot\eta^-_{\ep_1,\ep_2,\ep_3}(x,y,t-s)\dy\ds. 
\end{align*}

\noindent To express $I_3$ in \eqref{heatE-I_3}, we choose the following test function $\Phi$ satisfying the condition for \eqref{weakform-NavierBC}, (div-free and orthogonal to $(\nabla' \gamma,-1)$ since $\pd_i \xi^-|_{y_d=0}=0$),
\begin{equation}\label{testFw_id}
T^{-1}\Phi(y,s)=(0,\cdots,0,\underbrace{\pd_{y_d}\xi^-}_{i-th},0,\cdots,0,\underbrace{-\pd_{y_{i}}\xi^-+\gamma_{i}'\cdot \pd_{y_{d}}\xi^-}_{d-th}),
\end{equation}
and denote
\[T^{-1} \Phi^{pert}\eqdefa(0,\cdots,0,\underbrace{\gamma_{i}'\cdot \pd_{y_{d}}\xi^-}_{d-th}).\]
Then we notice that (no boundary term using $\pd_d^2\xi^-=0=T^{-1}\tilde u_d$ when $y_d=0$)
\begin{align}\label{I_3-Phi-p}
I_3&=\int_{\R^d_{+}\times\R} \left(-T^{-1}\tilde{\bm{u}} \cdot \partial_sT^{-1}(\Phi-\Phi^{pert})+\nabla_y\, T^{-1}\tilde{\bm{u}}: \nabla_y\, T^{-1}(\Phi-\Phi^{pert})\right)\dy\ds\notag  \\[6pt]
&=\int_{\R^d_{+}\times\R} -T^{-1}\bm{u} \cdot \partial_sT^{-1}(\Phi-\Phi^{pert})+\nabla_y\, T^{-1}\bm{u}: \nabla_y\, T^{-1}(\Phi-\Phi^{pert})\dy\ds\notag\\[6pt]
&\,\quad -\sum_{k=1}^{d-1}\int_{\R^d_{+}\times\R}\gamma'_k\cdot T^{-1}u_k \cdot \partial_s\pd_{y_i}\xi^- -\nabla_y\,\left(\gamma'_k\cdot T^{-1}u_k\right)\cdot \nabla_y\, (\pd_{y_i}\xi^-)\dy\ds. 
\end{align}

\noindent By plugging \eqref{testFw_id} into \eqref{Stokes-curved}, and using \eqref{Stokes-curved-term2}, \eqref{I_3-Phi-p}, we see that $I_3$ can be expressed as 
\[I_3=\tsum_{k=1}^{8}I_{3k},\]
where 
\begin{align*}
I_{31}&\eqdefa\int_{\R^d_{+}\times\R}T^{-1}u_d\cdot\pd_{s}(\gamma_{i}' \cdot \pd_{y_d}\xi^-)-\nabla_y\, T^{-1}u_d\cdot\nabla_y\, (\gamma_{i}'\cdot \pd_{y_{d}}\,\xi^-)\dy\ds  \\[6pt]
&\eqdefa I_{311}+I_{312}
\end{align*}
with the main term
\begin{equation*}
I_{311}\eqdefa-\int_{\R^d_{+}\times\R}\pd_{y_d}\, T^{-1}u_d\cdot \gamma_{i}'\cdot \pd_{y_{d}}^2\,\xi^-\dy\ds,
\end{equation*}
\begin{align*}
I_{312} &\eqdefa \int_{\R^d_{+}\times\R}T^{-1}u_d\cdot \bkt{ \pd_{s}(\gamma_{i}' \cdot \pd_{y_d}\xi^-)+ \tsum_{k=1}^{d-1} \pd_{y_k}  (\gamma_{i}'\cdot \pd_{y_k} \pd_{y_{d}}\,\xi^-)}\dy\ds\\
&\,\quad -\sum_{k=1}^{d-1}\int_{\R^d_{+}\times\R} \pd_{y_k}T^{-1}u_d\cdot \gamma_{ik}''\cdot \pd_{y_{d}}\,\xi^- \dy\ds, 
\end{align*}
and
\begin{align*}
I_{32}&\eqdefa  -\sum_{k=1}^{d-1}\int_{\R^d_{+}\times\R}\gamma'_k\cdot T^{-1}u_k \cdot \partial_s\pd_{y_i}\xi^- -\nabla_y\,\left(\gamma'_k\cdot T^{-1}u_k\right)\cdot \nabla_y\, (\pd_{y_i}\xi^-)\dy\ds
\\
&=-\sum_{k=1}^{d-1}\int_{\R^d_{+}\times\R}\gamma'_k\cdot T^{-1}u_k \cdot (\partial_s + \tsum_{k=1}^{d-1} \pd_{y_k}^2)\pd_{y_i}\xi^-\dy\ds
\\
&\quad\qquad +  \int_{\R^d_{+}\times\R}\pd_{y_d} \,\left(\gamma'_k\cdot T^{-1}u_k\right)\cdot \pd_{y_d}\, (\pd_{y_i}\xi^-)\dy\ds,
\end{align*}
\begin{align*}
I_{33}&\eqdefa\,\sum_{k=1}^{d-1}\int_{\R^d_{+}\times\R}\gamma_{k}'\cdot\Big(\pd_{y_k}\,T^{-1}u_i\cdot \pd_{y_{d}}^2\,\xi^{-}-\pd_{y_k}T^{-1}u_d\cdot \pd_{y_{d}}\left(\pd_{y_{i}}-\gamma_{i}'\cdot\pd_{y_d}\right)\xi^{-}\Big)\dy\ds  \\[6pt]
&=-\sum_{k=1}^{d-1}\int_{\R^d_{+}\times\R}T^{-1}u_i\cdot \pd_{y_k}\left(\gamma_{k}'\cdot\pd_{y_{d}}^2\xi^{-}\right)-T^{-1}u_d\cdot\pd_{y_k}\pd_{y_d} \left(\gamma_{k}'\cdot\left(\pd_{y_{i}}-\gamma_{i}'\cdot\pd_{y_d}\right)\xi^{-}\right)\dy\ds, 
\end{align*}
\begin{align*}
I_{34}\eqdefa\,\sum_{k=1}^{d-1}\int_{\R^d_{+}\times\R}\gamma_{k}'\cdot\Big(\pd_{y_d}\,T^{-1}u_i\cdot \pd_{y_{d}}\,\pd_{y_{k}}\xi^{-}-\pd_{y_d}T^{-1}u_d\cdot \pd_{y_{k}} \left(\pd_{y_{i}}-\gamma_{i}'\cdot\pd_{y_d}\right)\xi^{-}\Big)\dy\ds, 
\end{align*}
\begin{align*}
I_{35}&\eqdefa\,-\int_{\R^d_{+}\times\R}|\nabla\gamma|^2\cdot\Big(\pd_{y_d}\,T^{-1}u_i\cdot \pd_{y_{d}}^2\xi^{-}-\pd_{y_d}T^{-1}u_d\cdot \pd_{y_{d}}\left(\pd_{y_{i}}-\gamma_{i}'\cdot\pd_{y_d}\right)\xi^{-}\Big)\dy\ds, 
\end{align*}
\begin{align*}
I_{36}&\eqdefa-\sum_{1\leq k,\ell\leq d}\int_{\R^d_{+}\times\R}T^{-1}\pd_{y_\ell}u_{k}\cdot T^{-1}\pd_{y_k}\Phi_{\ell}\dy\ds   \\[6pt]
&=-\sum_{\ell=i,d}\sum_{1\leq k\leq d}\int_{\R^d_{+}\times\R}\left(\pd_{y_\ell}-\gamma_{\ell}'\cdot\pd_{y_d}\right)T^{-1}u_{k}\cdot \left(\pd_{y_k}-\gamma_{k}'\cdot\pd_{y_d}\right)T^{-1}\Phi_{\ell}\dy\ds\\[6pt]
&=-\sum_{1\leq k\leq d-1}\int_{\R^d_{+}\times\R}\pd_{y_d}\left(T^{-1}u_{k}\cdot \gamma_{ik}''\cdot \pd_{y_{d}}\xi^-\right)\dy\ds, 
\end{align*}
\begin{align*}
I_{37}&\eqdefa-\int_{\R^{d-1}\times \R}T^{-1}\left(\alpha \bm{u}\cdot\Phi\right)\cdot \sqrt{1+|\nabla \gamma|^2}\dy'\ds  \\[6pt]
&=\alpha \int_{\R^d_{+}\times\R}\pd_{y_d}\Big(T^{-1}u_{i}\cdot \pd_{y_d}\xi^--T^{-1}u_{d}\cdot\left(\pd_{y_i}-\gamma_{i}'\cdot\pd_{y_d}\right)\xi^-\Big)\cdot \sqrt{1+|\nabla \gamma|^2}\dy\ds, 
\end{align*}
\begin{align*}
I_{38}&\eqdefa\int_{\R^d_{+}\times\R}T^{-1}\left( \bm{f}\cdot\Phi-\CF:\nabla\Phi \right)\dy\ds  \\[6pt]
&=\int_{\R^d_{+}\times\R}T^{-1} f_i\cdot\pd_{y_d}\xi^--T^{-1}f_d\cdot\left(\pd_{y_i}-\gamma_{i}'\cdot\pd_{y_d}\right)\xi^- \dy\ds\\[6pt]
&\quad\ - \sum_{1\leq\ell\leq d}\,\int_{\R^d_{+}\times\R}T^{-1}\CF_{i\ell}\cdot\pd_{y_\ell}\pd_{y_d}\xi^--T^{-1}\CF_{d\ell}\cdot\pd_{y_\ell}\left(\pd_{y_i}-\gamma_{i}'\cdot\pd_{y_d}\right)\xi^-\dy\ds. 
\end{align*}

\subsection{Gradient estimate with an extra term}

Our proof of Theorem \ref{thm_navierBC} is divided into several steps. In this subsection we will prove a local estimate \eqref{est-lem-navier2} which differs from Theorem \ref{thm_navierBC} with an extra term $\norm{\nabla u}_{L^1}$ on the right side.

\begin{prop}\label{lem2-NavierBC}
\sl{
 Under assumptions in Theorem \ref{thm_navierBC}, there exists some $\rho_0\ll 1$ such that for $\rho< \rho_0$, we have $\psi_1(t)\cdot\nabla T^{-1} \bm{u}\in L^{q,r}\big(\C^+_{\rho}\times(-1,0)\big)$ and 
\begin{align}\label{est-lem-navier2}
\| \psi_1(t)\cdot\nabla T^{-1}\bm{u}\|_{L^{q,r}\left(\C^+_{\rho}\times(-1,0)\right)}&\lec \rho^{-1}\|\bm{u}\|_{L^{q,r}\left(Q^+_1\right)}
+\|\bm{f}\|_{L^{q_*,r}\left(Q^+_1\right)}+\|\CF\|_{L^{q,r}\left(Q^+_1\right)}
\notag\\[6pt] 
&\,\quad +\|\nabla \bm{u}\|_{L^{1}\left(Q^+_1\right)}.
\end{align}
}
\end{prop}

Proposition \ref{lem2-NavierBC} follows from bootstrapping in the indexes $q_1,q_2,r_1,r_2$ in the following Lemma \ref{lem1-navierBC} and the fact that $\nabla\bm{u}\in L^{1}(Q^+_1)$.

\begin{lem} \label{lem1-navierBC}
\sl{
Let $1\leq q_1\leq  q_2\leq q$, $1\leq r_1\leq r_2\leq r$,  $1<q_2,r_2$, and
 $d\left(\frac{1}{q_1}-\frac{1}{q_2}\right)+2\left(\frac{1}{r_1}-\frac{1}{r_2}\right)<1$. Under assumptions in Theorem \ref{thm_navierBC},  there exists some $0<\rho_0< 1$ such that for  $0<\rho< \rho_0$, we have $\psi_1(t)\cdot\nabla T^{-1} \bm{u}\in L^{q_2,r_2}\big(\C^+_{\rho/2}\times(-1,0)\big)$ and
\begin{align}\label{est-lem-navier}
 \| \psi_1(t)\cdot\nabla T^{-1}\bm{u}\|_{L^{q_2,r_2}\left(\C^+_{\rho/2}\times(-1,0)\right)}&\lec \rho^{-1}\|\bm{u}\|_{L^{q,r}\left(Q^+_1\right)}
 +\|\bm{f}\|_{L^{q_*,r}\left(Q^+_1\right)}+\|\CF\|_{L^{q,r}\left(Q^+_1\right)}
 \notag\\[6pt] 
 &\,\quad +\|\psi_1(t)\cdot\nabla T^{-1} \bm{u}\|_{L^{q_1,r_1}\left(\C^+_{\rho}\times(-1,0)\right)},
\end{align}
provided that $\psi_1(t)\cdot\nabla T^{-1}\bm{u}\in L^{q_1,r_1}(\C_{\rho}^+\times(-1,0))$}.
\end{lem}

Note that a pressure bound is not needed on the right side of \eqref{est-lem-navier}, as shown for the flat boundary case in \cite{CLT24}.

\begin{proof}
Let $\rho_1=\frac{\rho}{2}$ in the cut-off function $\zeta$ in \eqref{def-zeta}.  As in the proof of \eqref{gra1}, we adopt 3 steps to obtain the desired results.

\smallskip

\noindent \textbf{\textit{Step 1: }} 
Classification of terms and the normal form
\smallskip

We will prove Lemma \ref{lem1-navierBC} using the heat equations  \eqref{heatE-I_1} and \eqref{heatE-I_3} of
$E^+_{\ep_1,\ep_2,\ep_3}(\omega_{ij})$ and $E^-_{\ep_1,\ep_2,\ep_3}(\omega_{id})$ for $1\leq i,j\leq d-1$.
However, \eqref{heatE-I_3} contains bad terms with three vertical derivative $\pd_{y_d}^3$ on the velocity in $I_{311}$ and $I_{35}$. We will first rewrite \eqref{heatE-I_3} in its \emph{normal form} \eqref{heatE-hatw}, similar to Step 1 of the proof of \eqref{gra1}.
Notice that by divergence-free property and \eqref{commu2},
\begin{align}
\pd_{y_d} T^{-1}u_d &= T^{-1}\pd_{y_d} u_d
=\sum_{k=1}^{d-1} T^{-1}\pd_{y_k}u_k =\sum_{k=1}^{d-1} \bke{\pd_{y_k}T^{-1}u_k -\gamma_{k}'\pd_{y_d}T^{-1}u_k} \nonumber\\
&=\sum_{k=1}^{d-1}\bke{ \pd_{y_k}T^{-1}u_k
-\gamma_{k}'\pd_{y_k}T^{-1}\tilde{u}_d
 +\gamma_{k}' \left(\pd_{y_k}T^{-1}\tilde{u}_d-\pd_{y_d}T^{-1}u_k\right) }.\label{eq4.17}
\end{align}
Hence
\begin{align*}
I_{311}&=\int_{\R^d_{+}\times\R} \pd_{y_d} T^{-1}u_d  \cdot \gamma_{i}'\cdot \pd_{y_{d}}^2\,\xi^-\dy\ds 
\\
&= \sum_{k=1}^{d-1} \int_{\R^d_{+}\times\R} \bigg(
-T^{-1}u_k \cdot\pd_{y_k} \bke{ \gamma_{i}'\cdot \pd_{y_{d}}^2\,\xi^-}
+ T^{-1}u_d \cdot\pd_{y_k}  \bke{ \gamma_{k}' \gamma_{i}'\cdot \pd_{y_{d}}^2\,\xi^-}
\\
&\qquad + \left(\pd_{y_k}T^{-1}\tilde{u}_d-\pd_{y_d}T^{-1}u_k\right)\cdot\gamma_{k}' \gamma_{i}'\cdot \pd_{y_{d}}^2\,\xi^-\bigg) \dy\ds = \pd_{x_d}^2\,I_{w1}+I_{311}^* 
\end{align*}
where $I_{311}^*$ is the error term and 
\begin{equation*}
I_{w1}\eqdefa\sum_{k=1}^{d-1} E^{-}_{\ep_1,\ep_2,\ep_3}(\omega_{kd}\cdot\gamma_{k}'\cdot\gamma_{i}'). 
\end{equation*}
Analogously, also using \eqref{eq4.17}, we decompose $I_{35}=\pd_{x_d}^2\,I_{w2}+I_{35}^*$ with
\begin{equation*}
I_{w2}\eqdefa  E^{-}_{\ep_1,\ep_2,\ep_3}(\omega_{id}\cdot|\nabla\gamma|^2)+\sum_{k=1}^{d-1} E^{-}_{\ep_1,\ep_2,\ep_3}(\omega_{kd}\cdot\gamma_{k}'\cdot\gamma_{i}'\cdot|\nabla\gamma|^2). 
\end{equation*}

For the terms $I_{w1}$ and $I_{w2}$ above, we want to move them to the LHS of equation \eqref{heatE-I_3}. 
This is another key step in which a \emph{normal form} is used, compare \eqref{heatE-OBC1} in Step 1 of the proof of \eqref{gra1}.
Hence we define for $1\leq i\leq d-1$,
\begin{equation}\label{omegaid1}
\hat{\omega}_{id}\eqdefa(1+|\nabla\gamma|^2)\cdot\omega_{id}+\sum_{k=1}^{d-1}(1+|\nabla\gamma|^2)\cdot\gamma'_k\cdot\gamma'_i\cdot\omega_{kd},
\end{equation}
which satisfies the heat equation
\begin{align}\label{heatE-hatw}
(\pd_t-\Delta)E^-_{\ep_1,\ep_2,\ep_3}(\hat{\omega}_{id})&=I_{311}^*+I_{312}+I_{32}+I_{33}+I_{34}+I_{35}^*+I_{36}+I_{37}+I_{38}+I_4\notag\\[6pt]
&\quad \ +\left(\pd_t-\tsum_{\ell=1}^{d-1}\pd_{x_\ell}^2\right)\left(I_{w1}+I_{w2}\right). 
\end{align}

For later use, we note that from \eqref{omegaid1}, we have
\begin{equation}\label{omegaid2}
\omega_{id}=\sum_{k=1}^{d-1}b_{ik}(\nabla \gamma)\cdot\hat{\omega}_{kd},
\end{equation}
where $b_{ik}(\nabla \gamma)$ is a matrix function close to identity.
Notice that 
\begin{align}
 E^-_{\ep_1,\ep_2,\ep_3}(\omega_{id})&=\sum_{k=1}^{d-1}b_{ik}(\nabla \gamma(x'))\cdot E^-_{\ep_1,\ep_2,\ep_3}(\hat{\omega}_{kd})+\sum_{k=1}^{d-1}\int_{\R^d_{+}\times\R}\left(b_{ik}(\nabla \gamma(y'))-b_{ik}(\nabla \gamma(x'))\right)\cdot\notag\\[6pt]
 &\,\quad  \cdot(1+|\nabla\gamma|^2)\left(\omega_{kd}+\sum_{j=1}^{d-1}\left(\gamma_{k}'\cdot\gamma_{j}'\cdot\omega_{jd}\right)\right)\cdot\eta_{\ep_1,\ep_2,\ep_3}^{-}(x,y,t-s)\dy\ds.  \label{0904a}
\end{align}
From here we claim
\begin{equation}\label{omegaid3}
  \|(-\Delta)^{-1}\nabla\nabla_{x'} E^-_{\ep_1,\ep_2,\ep_3}(\omega_{id})\|_{\cW^{q_2,r_2}}
\lec \sum_{k=1}^{d-1}\|E^-_{\ep_1,\ep_2,\ep_3}(\hat{\omega}_{kd})\|_{\cW^{q_2,r_2}}+\|\bm{u}\|_{L^{q,r}\left(Q^+_1\right)}.
\end{equation}
Indeed, the first term on the RHS of \eqref{omegaid3} comes from the first term on the RHS of \eqref{0904a} and Lemma \ref{lem-C-Z est}. The integral of \eqref{0904a} is bounded by 
$\|\bm{u}\|_{L^{q,r}(Q^+_1)}$, similar to the proof of Lemma \ref{lem0} (v).  Specifically, we decompose $\omega_{kd}= \pd_k v_d-\pd_d v_k$ and integrate by parts to move these derivatives $\pd_k$ and $\pd_d$ to other terms.  For $ \pd_k v_d$,  the integral becomes $v_d$ multiplied by $\eta_{\ep_1,\ep_2,\ep_3}^{-}$ or $(b_{ik}(\nabla \gamma(y'))-b_{ik}(\nabla \gamma(x')))\pd_{y_k}\eta_{\ep_1,\ep_2,\ep_3}^{-}$ (which is bounded by $\ep_1|\pd_{y_k}\eta_{\ep_1,\ep_2,\ep_3}^{-}|$) and some bounded functions. Its estimates follow from Lemma \ref{lem-C-Z est}. For $ \pd_d v_k$, using $\eta_{\ep_1,\ep_2,\ep_3}^{-}|_{y_d=0}=0$,
the integral becomes $v_k$ multiplied by some function and $\pd_{y_d}\eta_{\ep_1,\ep_2,\ep_3}^{-}=-\pd_{x_d}\eta_{\ep_1,\ep_2,\ep_3}^{+}$.
We pull $\partial_{x_d}$ outside the integral and move $\nabla_{x'}$ inside the integral. In this way, we are reduced to the previous case, where the term $(b_{ik}(\nabla \gamma(y'))-b_{ik}(\nabla \gamma(x')))\pd_{y_k}\eta_{\ep_1,\ep_2,\ep_3}^{-}$ is controlled by $\ep_1|\pd_{y_k}\eta_{\ep_1,\ep_2,\ep_3}^{-}|$. The terms involving $\omega_{jd}$ are treated in the same way.

We now estimate $E^+_{\ep_1,\ep_2,\ep_3}(\omega_{ij})$ and $E^-_{\ep_1,\ep_2,\ep_3}(\hat \omega_{id})$ using their equations \eqref{heatE-I_1} and \eqref{heatE-hatw}.
Note that all terms in the RHS of \eqref{heatE-I_1} and \eqref{heatE-hatw} either involve the velocity $\bm{u}$ or the external forces $\bm{f},\CF$, and can be classified into the following six types. The terms involving the velocity can be written in the form of
\begin{equation}
\pd_{t}^{a_1}\nabla_{x}^{a_2} E_{\ep_1,\ep_2,\ep_3}^{\pm}\left(\nabla^{a_3} T^{-1}\bm{u}\cdot \left(\nabla^2\gamma\right)^{a_4}\cdot g_1(\nabla\gamma)\cdot\pd_{t}^{a_5}\nabla^{a_6}\zeta\right)
\end{equation} 
with non-negative integers  $2a_1+a_2+a_3+a_4+2a_5+a_6\leq 3$, $a_3,a_4\leq 1$ and some smooth bounded function $g_1$. For these terms we do not integrate by parts. For derivatives $\nabla_y$ and $ \pd_s$ acting on $\eta^\pm_{\ep_1,\ep_2,\ep_3}$, we use
Lemma \ref{lem0}(i) to change $\nabla_y$ to $\nabla_x$, $\pd_s \eta^\pm_{\ep_1,\ep_2,\ep_3}=-\pd_t\eta^\pm_{\ep_1,\ep_2,\ep_3}$, and pull $\nabla_{x}$ and $\pd_t$ outside of the integral. They are called higher order terms if $2a_1+a_2+a_3=3$ (all derivatives hit on $T^{-1}\bm{u} $), otherwise lower order terms. Higher order terms always come with $g_1(0)=0$, hence $|g_1(\nabla\gamma)|\lec |\nabla\gamma|$, and can be considered as perturbation terms.
\begin{enumerate}
\renewcommand{\theenumi}{\roman{enumi}}
\item Lower order terms $H_1$ (e.g.~$\nabla^2 E^\pm_{\ep_1,\ep_2,\ep_3}(T^{-1}\bm{u} \cdot  \nabla\zeta)$)  for $a_3=0$.
 By using Lemma \ref{lem-C-Z est}, we have 
\begin{equation}\label{est-H1}
\|\int_{-1}^{t}e^{(t-s)\Delta} H_1\ds\|_{L^{q_2,r_2}\left(\C_{2\rho}\times(-1,0)\right)}\lec \rho^{-1} \|\bm{u}\|_{L^{q,r}\left(Q^+_1\right)}.
\end{equation}

\item Lower order terms $H_2$ (e.g.~$\nabla E^\pm_{\ep_1,\ep_2,\ep_3}(\nabla T^{-1}\bm{u} \cdot \nabla \zeta)$)  for $a_3=1$ and $a_1=a_5=0$. If $a_4=0$, either $g_1(0)=0$ or $a_2+a_6=1$ (hence there will be no coefficient $1/\rho$ in the RHS of the estimate below).
For these terms, by Lemma \ref{lem-C-Z est} and \eqref{rho2}, we have
\begin{equation}\label{H_2-est}
\|\int_{-1}^{t}e^{(t-s)\Delta} H_2\ds\|_{L^{q_2,r_2}\left(\C_{2\rho}\times(-1,0)\right)}\lec \|\psi_1(t)\cdot\nabla T^{-1} \bm{u}\|_{L^{q_1,r_1}\left(\C^+_{\rho}\times(-1,0)\right)}.
\end{equation}

\item Higher order term $H_3$ (e.g.~$\pd_t \nabla E^\pm_{\ep_1,\ep_2,\ep_3}(T^{-1}\bm{u} \cdot g_1(\nabla \gamma) \zeta)$) for $a_1=a_2=1$. Moreover, if $\nabla^{a_2}$  refers to $\pd_{x_{d}}$, $H_3$ contains only the even extension
$\pd_{t}\pd_{x_{d}}E_{\ep_1,\ep_2,\ep_3}^{+}\left(T^{-1}\bm{u}\cdot g_1(\nabla\gamma)\cdot\zeta\right)$.
This will be crucial in the later estimate \eqref{H7-est}.
By Lemma \ref{lem-C-Z est} and ``bounding'' $\pd_{t} $ by $\ep_3^{-1}$, we have 
\begin{equation}\label{eL2}
\|\int_{-1}^{t}e^{(t-s)\Delta} H_3\ds\|_{L^{q_2,r_2}\left(\C_{2\rho}\times(-1,0)\right)}\lec  \ep_3^{-1}  \|\bm{u}\|_{L^{q,r}\left(Q^+_1\right)}.
\end{equation}

\item Higher order terms $H_4$ (e.g.~$\nabla^2 \nabla' E^\pm_{\ep_1,\ep_2,\ep_3}(T^{-1}\bm{u} \cdot g_1(\nabla \gamma) \zeta)$) for $a_1=a_3=0$ and $\nabla_{x}^{a_2}$ involving at least one horizontal derivative $\nabla_{x'}$. 
By  the estimate of commutator \eqref{commu3} and Lemma \ref{lem-C-Z est}, we have
\begin{align*}
 \|\int_{-1}^{t}e^{(t-s)\Delta} H_4\ds\|_{\cW^{q_2,r_2}}
\lec \|\bm{u}\|_{L^{q,r}\left(Q^+_1\right)}
+\|\nabla\gamma\|_{L^{\infty}\left(B_{\rho}'\right)}\cdot \|\nabla_{x'}E^{\pm}_{\ep_1,\ep_2,\ep_3}(\bm{v})\|_{\cW^{q_2,r_2}}. 
\end{align*}

\item Higher order terms $H_5$ (e.g.~$\nabla \nabla' E^\pm_{\ep_1,\ep_2,\ep_3}(\pd_d T^{-1}\bm{u} \cdot g_1(\nabla \gamma) \zeta)$) for $a_1=0$, $a_3=1$ and $\nabla_{x}^{a_2}$ involving at least one horizontal derivative $\nabla_{x'}$. Actually, we only consider the case that $\nabla^{a_3}$ refers to the vertical derivative $\pd_{y_d}$, because other $H_5$ with $\nabla^{a_3}=\nabla_{x'}$ can be transferred to $H_1+H_4$ using integration by parts as we did for $I_{12}$, $I_{312}$ and $I_{32}$.
For these terms, using Lemma \ref{lem-C-Z est}(v) and ``bounding'' $\nabla' $ by $\ep_1^{-1}$, we have
\begin{align}\label{estH3-1}
 \|\int_{-1}^{t}e^{(t-s)\Delta} H_5\ds\|_{L^{q_2,r_2}\left(\C_{2\rho}\times(-1,0)\right)}
\lec \ep_1^{-1}\|\psi_1(t)\cdot\pd_d T^{-1}\bm{u}\|_{L^{q_1,r_1}\big(\C^+_{\rho}\times(-1,0)\big)}.
\end{align}

\item Terms $H_6$ involving $\bm{f}$ and $\CF$. For these terms,  we have
 \begin{equation}\label{est-H6}
\|\int_{-1}^{t}e^{(t-s)\Delta}  H_6\ds\|_{L^{q,r}\left(\C_{2\rho}\times(-1,0)\right)}\lec \|\bm{f}\|_{L^{q_*,r}\left(Q^+_1\right)}+\|\CF\|_{L^{q,r}\left(Q^+_1\right)}.
\end{equation}
\end{enumerate}

Accordingly, we obtain that for $1\leq i\leq j\leq d-1$,
\begin{align}\label{limit1}
&\|E^+_{\ep_1,\ep_2,\ep_3}(\omega_{ij}),\, E^-_{\ep_1,\ep_2,\ep_3}(\hat{\omega}_{id})\|_{\cW^{q_2,r_2}}\lec  \left(\ep_3^{-1}+\rho^{-1}\right)\|\bm{u}\|_{L^{q,r}\left(Q^+_1\right)}
 +\|\bm{f}\|_{L^{q_*,r}\left(Q^+_1\right)}\notag\\[6pt]
&\,\qquad +\|\CF\|_{L^{q,r}\left(Q^+_1\right)}+\ep_1^{-1}\|\psi_1(t)\cdot\nabla T^{-1} \bm{u}\|_{L^{q_1,r_1}\left(\C^+_{\rho}\times(-1,0)\right)}
\notag\\[6pt] 
&\,\qquad +\|\nabla\gamma\|_{L^{\infty}\left(B_{\rho}'\right)}\cdot \|\nabla_{x'}E^{\pm}_{\ep_1,\ep_2,\ep_3}(\bm{v})\|_{\cW^{q_2,r_2}}. 
\end{align}
Notice that $\|\psi_1(t)\cdot\nabla T^{-1} \bm{u}\|_{L^{q_1,r_1}\left(\C^+_{\rho}\times(-1,0)\right)}<\infty$ by assumption of Lemma \ref{lem1-navierBC}, and  the last term $\|\nabla_{x'}E^{\pm}_{\ep_1,\ep_2,\ep_3}(\bm{v})\|_{\cW^{q_2,r_2}}$ 
 has no vertical derivative and remains bounded as $\ep_2\to 0$. By taking $\ep_2\to 0$, we get that for $1\leq i< j\leq d-1$,
\begin{equation}\label{boundedness of ep2}
\|E^{h}_{\ep_1,\ep_3}(\omega_{ij})\|_{\cH^{q_2,r_2}}+\|E^{h}_{\ep_1,\ep_3}(\hat{\omega}_{id})\|_{\cH^{q_2,r_2}}<\infty. 
\end{equation}

\bigskip

\noindent \textbf{\textit{Step 2: }} We will show
\begin{equation}\label{claim2}
\nabla E_{\ep_3}(\bm{v})\in  \cH^{q_2,r_2}.
\end{equation}

\smallskip

We give a different estimate of the higher order term $H_5$, to remove the $\ep_1^{-1}$ in the RHS of \eqref{limit1}. Using divergence-free property \eqref{div2} to change $\pd_dv_d$, $\pd_dv_i=-\omega_{id}+\pd_iv_d$ for $1\leq i\leq d-1$, and \eqref{omegaid2}, $H_5$ can be decomposed into
\begin{equation}\label{H5-2}
H_5=\nabla_{x'}\nabla E^{\pm}_{\ep_1,\ep_2,\ep_3}\left( \hat{\omega}_{kd}\cdot b_{ik}(\nabla\gamma) \cdot g_{2}(\nabla\gamma)\right)+H_1+H_4, 
\end{equation}
for some smooth function $g_2$ with $g_2(0)=0$. 
We have replaced  ${\omega}_{id}$ by $\hat{\omega}_{kd}$, which is to be absorbed into the LHS. Rewrite $E^{\pm}_{\ep_1,\ep_2,\ep_3}\left( \hat{\omega}_{kd}\cdot b_{ik}(\nabla\gamma) \cdot g_{2}(\nabla\gamma)\right)$ as in \eqref{0904a}.
By a similar proof of \eqref{omegaid3}, we get the estimate
\begin{align}
& \|\int_{-1}^{t}e^{(t-s)\Delta}\, \nabla_{x'} \nabla\,E^{\pm}_{\ep_1,\ep_2,\ep_3}\left( \hat{\omega}_{kd}\cdot b_{ik}(\nabla\gamma) \cdot g_{2}(\nabla\gamma)\right)\ds\|_{L^{q_2,r_2}\left(\C_{2\rho}\times(-1,0)\right)} \notag\\[0pt]
&\lec \sum_{k=1}^{d-1}\|\nabla\gamma\|_{L^{\infty}\left(B_{\rho}'\right)}\cdot\|E^{\pm}_{\ep_1,\ep_2,\ep_3}(\hat{\omega}_{kd})\|_{\cW^{q_2,r_2}}+\|\bm{u}\|_{L^{q,r}\left(Q^+_1\right)}\label{0911a}\\
&\,\quad +\|\psi_1(t)\cdot\nabla T^{-1} \bm{u}\|_{L^{q_1,r_1}\left(\C^+_{\rho}\times(-1,0)\right)}.\notag
\end{align}
Comparing with \eqref{omegaid3}, we gain a small factor $\|\nabla\gamma\|_{L^{\infty}}$ thanks to $g_2$, and 
the extra last term is because $\pd_d$ cannot be moved to $\eta^+_{\ep_1,\ep_2,\ep_3}$ by integration by parts in the case of $E^+_{\ep_1,\ep_2,\ep_3}$.

Then, by replacing \eqref{estH3-1} with \eqref{0911a}, we get instead of \eqref{limit1} another inequality,
\begin{align*}
&\|E^+_{\ep_1,\ep_2,\ep_3}(\omega_{ij}), \, E^-_{\ep_1,\ep_2,\ep_3}(\hat{\omega}_{id})\|_{\cW^{q_2,r_2}}\lec  \mathcal M
\notag\\[6pt] 
&\,\quad +\|\nabla\gamma\|_{L^{\infty}\left(B_{\rho}'\right)}\cdot\left( \|\nabla_{x'}E^{\pm}_{\ep_1,\ep_2,\ep_3}(\bm{v})\|_{\cW^{q_2,r_2}}+\sum_{k=1}^{d-1}\|E^{\pm}_{\ep_1,\ep_2,\ep_3}\left(\hat{\omega}_{kd}\right)\|_{\cW^{q_2,r_2}}\right),
\end{align*}
where
\[
\mathcal M=\left(\frac1{\ep_3}+\frac1\rho\right)\|\bm{u}\|_{L^{q,r}\left(Q^+_1\right)}
+\|\bm{f}\|_{L^{q_*,r}\left(Q^+_1\right)}+\|\CF\|_{L^{q,r}\left(Q^+_1\right)}+\|\psi_1(t)\cdot\nabla T^{-1} \bm{u}\|_{L^{q_1,r_1}\left(\C^+_{\rho}\times(-1,0)\right)}.
\]
We now take $\ep_2\to 0$.
In view of \eqref{boundedness of ep2} the conclusion of step 1, the last term remains finite in the limit and can be absorbed into the LHS for $\rho<\rho_0\ll 1$. We get
\begin{align}\label{limit2}
\|E^{h}_{\ep_1,\ep_3}(\omega_{ij}),\,E^{h}_{\ep_1,\ep_3}(\hat{\omega}_{id})\|_{\cH^{q_2,r_2}}\lec \mathcal M +\|\nabla\gamma\|_{L^{\infty}\left(B_{\rho}'\right)}\cdot \|\nabla_{x'}E^{h}_{\ep_1,\ep_3}(\bm{v})\|_{\cH^{q_2,r_2}}. 
\end{align}
Note that we have removed the $\ep_1^{-1}$ term of \eqref{limit1}.

\smallskip

We get a formula for $\nabla_{x'}\,E^{\#}_{\ep_1,\ep_2,\ep_3}(\bm{v})$
by applying $\nabla_{x'}(-\Delta)^{-1}$ to 
\eqref{lap-Ev}. 
By this formula, \eqref{E-sharp-v},
\eqref{divcontrol},  \eqref{omegaid3}, \eqref{limit2}, and Lemma \ref{lem-C-Z est}, we get
\begin{align*} \|\nabla_{x'}\,E^{\#}_{\ep_1,\ep_2,\ep_3}(\bm{v})\|_{\cW^{q_2,r_2}} \lec 
\mathcal M +\|\nabla\gamma\|_{L^{\infty}\left(B_{\rho}'\right)}\cdot \|\nabla_{x'}E^{h}_{\ep_1,\ep_3}(\bm{v})\|_{\cH^{q_2,r_2}}.
\end{align*}
This way of estimating does not extend to $\|\pd_{x_d}\,E^{\#}_{\ep_1,\ep_2,\ep_3}(\bm{v})\|_{\cW^{q_2,r_2}} $ because we used \eqref{omegaid3}.

By taking $\ep_2\to 0$ in the above inequality, we get  that for $\rho<\rho_0\ll 1$,
\begin{align*}
 \|\nabla_{x'}\,E^{h}_{\ep_1,\ep_3}(\bm{v})\|_{\cH^{q_2,r_2}} &\lec \mathcal M.
\end{align*}
By substituting the above inequality back into \eqref{limit2} and taking $\ep_1\to 0$, we get that for $1\leq i,j\leq d-1$
\begin{align*}
\|E_{\ep_3}(\omega_{ij}),\,E_{\ep_3}(\hat{\omega}_{id})\|_{\cH^{q_2,r_2}}&\lec  \mathcal M,
\end{align*}
which together with \eqref{omegaid2} implies that 
\begin{align}\label{limt5}
\|E_{\ep_3}(\omega_{ij}),\,E_{\ep_3}(\omega_{id})\|_{\cH^{q_2,r_2}}
&\lec  \mathcal M.
\end{align}

Apply $\nabla_{x}(-\Delta)^{-1}$ to \eqref{lap-Ev} to get a formula for $\nabla_{x}\,E^{\#}_{\ep_1,\ep_2,\ep_3}(\bm{v})$.
By virtue of \eqref{E-sharp-v}, \eqref{divcontrol}, \eqref{limt5}, and Lemma \ref{lem-C-Z est}, we achieve
\begin{align*}
\|\nabla\,E^{\#}_{\ep_1,\ep_2,\ep_3}(\bm{v})\|_{\cW^{q_2,r_2}} 
&\lec \mathcal M,
\end{align*}
which leads to $\|\nabla\,E_{\ep_3}(\bm{v})\|_{\cH^{q_2,r_2}} \lec \mathcal M$ by taking $\ep_1,\ep_2\rightarrow0$, establishing \eqref{claim2}.

\bigskip

\noindent \textbf{\textit{Step 3: }} The $L^{q_2,r_2}$ boundedness of $\nabla \bm{v}$ in $\C^+_{\rho}\times(-1,0)$.

\smallskip
Now we revisit the heat equations \eqref{heatE-I_1} and \eqref{heatE-I_3}. The latter is equal to
\begin{align}\label{heatE-wid}
(\pd_t-\Delta)E^-_{\ep_1,\ep_2,\ep_3}(\omega_{id})&=I_{311}^*+I_{312}+I_{32}+I_{33}+I_{34}+I_{35}^*+I_{36}+I_{37}+I_{38}+I_4\notag\\[6pt]
&\quad \ +\pd_{x_d}^2\left(I_{w1}+I_{w2}\right).
\end{align}
We have new higher order terms in $\pd_{x_d}^2\left(I_{w1}+I_{w2}\right)$. Denote by $H_7$ 
the sum of all higher order terms (terms with 3 derivatives on $\bm{u}$)
in the RHS of either  \eqref{heatE-I_1} or \eqref{heatE-wid}. 
With \eqref{claim2}, $H_4$ type terms are estimated as before, $H_5$ terms and $\pd_{x_d}^2\left(I_{w1}+I_{w2}\right)$ are estimated directly. For $H_3$ type terms of the form $\pd_t \nabla E^\pm_{\ep_1,\ep_2,\ep_3}(T^{-1}\bm{u} \cdot g_1(\nabla \gamma) \zeta)$, we move $\nabla$ inside the integral and integrate by parts.  In this process there is no boundary term because there is no $H_3$ term like $\pd_t \pd_{x_d} E^-_\ep[T^{-1}\bm{u} \cdot g_1(\nabla \gamma) \zeta]$. We get the following estimate
\begin{align}\label{H7-est}
\|\int_{-1}^{t}e^{(t-s)\Delta} H_7\ds\|_{\cW^{q_{2},r_2}}\lec  \|\nabla\gamma\|_{L^{\infty}\left(B_{\rho}'\right)}\cdot\|\nabla E_{\ep_3}(\bm{v})\|_{\cH^{q_2,r_2}}+\rho^{-1}\|\bm{u}\|_{L^{q,r}\left(Q^+_1\right)}.
\end{align}
Estimating \eqref{heatE-I_1} and \eqref{heatE-I_3} using
\eqref{est-H1}, \eqref{H_2-est}, \eqref{est-H6} and \eqref{H7-est}, we get
\begin{align}\label{limit8}
\|E^+_{\ep_1,\ep_2,\ep_3}(\omega_{ij}),\,E^-_{\ep_1,\ep_2,\ep_3}(\omega_{id})\|_{\cW^{q_2,r_2}}\lec  \rho^{-1}\|\bm{u}\|_{L^{q,r}\left(Q^+_1\right)}
+\|\bm{f}\|_{L^{q_*,r}\left(Q^+_1\right)}+\|\CF\|_{L^{q,r}\left(Q^+_1\right)}
\notag\\[6pt] 
+\|\psi_1(t)\cdot\nabla T^{-1} \bm{u}\|_{L^{q_1,r_1}\left(\C^+_{\rho}\times(-1,0)\right)}+\|\nabla\gamma\|_{L^{\infty}\left(B_{\rho}'\right)}\cdot\|\nabla E_{\ep_3}(\bm{v})\|_{\cH^{q_2,r_2}}.
\end{align}
By virtue of \eqref{lap-Ev}, \eqref{divcontrol} and \eqref{limit8}, we achieve that for $\rho<\rho_0\ll 1$,
\begin{align*}
\|\nabla\,E_{\ep_3}(\bm{v})\|_{\cH^{q_2,r_2}} 
&\lec \rho^{-1}\|\bm{u}\|_{L^{q,r}\left(Q^+_1\right)}
+\|\bm{f}\|_{L^{q_*,r}\left(Q^+_1\right)}+\|\CF\|_{L^{q,r}\left(Q^+_1\right)}
\notag\\[6pt] 
&\,\quad +\|\psi_1(t)\cdot\nabla T^{-1} \bm{u}\|_{L^{q_1,r_1}\left(\C^+_{\rho}\times(-1,0)\right)},
\end{align*}
which leads to \eqref{est-lem-navier} by taking $\ep_3\rightarrow0$, showing Lemma \ref{lem1-navierBC}.
\end{proof}

\subsection{Removing the extra term}

Now we begin to remove the term $\|\nabla\bm{u}\|_{L^{1}\left(Q^+_1\right)}$ in the RHS of \eqref{est-lem-navier2}. Notice that it comes from the estimate of $H_2$ in \eqref{H_2-est}. More precisely, $H_2$ with $a_4=1$ and $\nabla^{a_3}$ referring to the horizontal derivative $\nabla'$ is what the difficulty lies in 
(see, for instance, the first term in $I_{12}$ and the last term in $I_{312}$). For such terms, we can not move the horizontal derivative $\nabla'$ away since $\gamma\in C^{1,1}$ only.  A new strategy is employed to overcome this difficulty.
\begin{lem}\label{lem3-NavierBC}
\sl{
Let $q_3=\min(q,r)$. Under assumptions in Theorem \ref{thm_navierBC}, there exists some $\rho_0\ll 1$ such that for $\rho< \rho_0$, we have
\begin{align}\label{limit3}
\|\psi_1(t)\cdot\nabla T^{-1}\bm{u}\|_{L^{q_3}\left(\C^+_{\rho}\times(-1,0)\right)}\lec \rho^{-1-\frac{d}{q_3}}\left(\|\bm{u}\|_{L^{q,r}\left(Q^+_1\right)}
+\|\bm{f}\|_{L^{q_*,r}\left(Q^+_1\right)}+\|\CF\|_{L^{q,r}\left(Q^+_1\right)}\right).
\end{align}
}
\end{lem}

\begin{proof}
For some $\rho_0$ whose value will be specified later (and smaller than the $\rho_0$ in Proposition \ref{lem2-NavierBC}), we will show that for any  $\rho_1<\rho<\rho_0/2$,
\begin{align}\label{limit4}
&\quad\ \|\psi_1(t)\cdot\nabla T^{-1} \bm{u}\|_{L^{q_3}\left(\C^+_{\rho_1}\times(-1,0)\right)}^{q_3}\notag\\[6pt]
&\leq C(\rho-\rho_1)^{-q_3-d}\left(\|\bm{u}\|_{L^{q,r}\left(Q^+_1\right)}
+\|\bm{f}\|_{L^{q_*,r}\left(Q^+_1\right)}+\|\CF\|_{L^{q,r}\left(Q^+_1\right)}\right)^{q_3}
\notag\\[6pt] 
&\,\quad +\frac{1}{2}\|\psi_1(t)\cdot\nabla T^{-1} \bm{u}\|_{L^{q_3}\left(\C^+_{\rho}\times(-1,0)\right)}^{q_3}.
\end{align}
Then by applying Lemma \ref{lemA2}, we get that for $\rho<\rho_0/4$,
\begin{equation*}
 \|\psi_1(t)\cdot\nabla T^{-1} \bm{u}\|_{L^{q_3}\left(\C^+_{\rho}\times(-1,0)\right)}^{q_3}
\lec\rho^{-q_3-d}\left(\|\bm{u}\|_{L^{q,r}\left(Q^+_1\right)}
+\|\bm{f}\|_{L^{q_*,r}\left(Q^+_1\right)}+\|\CF\|_{L^{q,r}\left(Q^+_1\right)}\right)^{q_3},
\end{equation*}
which gives rise to \eqref{limit3}.

To show \eqref{limit4},
we choose a new cut-off function $\zeta= \zeta_{\rho_2/2,\rho_2}(x-x_0)\cdot\psi_1(t)$ in \eqref{defv1}, where $x_0$ is any point on $\C_{\rho_0/2}\cap \pd\R^{d}_{+}$, and $\rho_2<\frac{\rho_0}{2}$  (we will let $\rho_2=(\rho-\rho_1)/4$). Denote $\C_R^+(x_0)=x_0+\C_{R}^+$ and $\C_R(x_0)=x_0+\C_R$. All estimates in the proof of Lemma \ref{lem1-navierBC} carry over to $\C_{2\rho_2}(x_0)$ (with smaller $\rho_0$). For instance, lower order terms
\begin{equation*}
\|\int_{-1}^{t}e^{(t-s)\Delta} H_1\ds\|_{L^{q_3}\left(\C_{2\rho_2}(x_0)\times(-1,0)\right)}\lec \rho_2^{-1} \|\bm{u}\|_{L^{q,r}\left(Q^+_1\right)},
\end{equation*}
all higher order terms
\begin{align*}
\|\int_{-1}^{t}e^{(t-s)\Delta} H_7\ds\|_{\cW^{q_3,q_3}}\lec  \|\nabla\gamma\|_{L^{\infty}\left(B_{\rho_0}'\right)}\cdot\|\nabla \bm{v}\|_{\cH^{q_3,q_3}}+\rho_2^{-1}\|\bm{u}\|_{L^{q,r}\left(Q^+_1\right)},
\end{align*}
and external-force terms
\begin{equation*}
\|\int_{-1}^{t}e^{(t-s)\Delta}  H_6\ds\|_{L^{q,r}\left(\C_{2\rho_2}(x_0)\times(-1,0)\right)}\lec \|\bm{f}\|_{L^{q_*,r}\left(Q^+_1\right)}+\|\CF\|_{L^{q,r}\left(Q^+_1\right)}.
\end{equation*}

We adopt a different estimate of $H_2$. By Lemma \ref{lem-C-Z est} and \eqref{rho2}, we get
\begin{equation*}
\|\int_{-1}^{t}e^{(t-s)\Delta} H_2\ds\|_{L^{q_3}\left(\C_{2\rho_2}(x_0)\times(-1,0)\right)}\lec \rho_0\|\psi_1(t)\cdot\nabla T^{-1} \bm{u}\|_{L^{q_3}\left(\C^+_{\rho_2}(x_0)\times(-1,0)\right)}.
\end{equation*}
It is similar to \eqref{H_2-est} but we gain a small factor because the same exponent $q_3$ in both sides. Analogous to the estimates in Step 3 of the proof for Lemma \ref{lem1-navierBC}, we get
\begin{align*}
\|\psi_1(t)\cdot\nabla T^{-1} \bm{u}\|_{L^{q_3}\left(\C^+_{\rho_2/2}(x_0)\times(-1,0)\right)} 
&\lec \rho_2^{-1}\|\bm{u}\|_{L^{q,r}\left(Q^+_1\right)}
+\|\bm{f}\|_{L^{q_*,r}\left(Q^+_1\right)}+\|\CF\|_{L^{q,r}\left(Q^+_1\right)}
\notag\\[6pt] 
&\,\quad +\rho_0\|\psi_1(t)\cdot\nabla T^{-1} \bm{u}\|_{L^{q_3}\left(\C^+_{\rho_2}(x_0)\times(-1,0)\right)},
\end{align*}
which implies that for sufficiently small $\rho_0$,
\begin{align}\label{limit9}
\|\psi_1(t)\cdot\nabla T^{-1} \bm{u}\|_{L^{q_3}\left(\C^+_{\rho_2/2}(x_0)\times(-1,0)\right)}^{q_3}
\leq C\rho_2^{-q_3}\left(\|\bm{u}, \CF\|_{L^{q,r}\left(Q^+_1\right)}
+\|\bm{f}\|_{L^{q_*,r}\left(Q^+_1\right)}\right)^{q_3}
\notag\\[6pt] 
+\frac{1}{2*100^{d}}\|\psi_1(t)\cdot\nabla T^{-1} \bm{u}\|_{L^{q_3}\left(\C^+_{\rho_2}(x_0)\times(-1,0)\right)}^{q_3}.
\end{align}
The constant $C$ above depends only on $q,r$ and $\|\gamma\|_{C^{1,1}}$, and is independent of $x_0$. Thus, by using \eqref{limit9}, interior estimate \eqref{interior1} and finite overlapping covering lemma (covering $\C^+_{\rho_1}$ with $\C_{(\rho-\rho_1)/4}^+(x_0)$ for boundary point $x_0\in \C_{\rho_1}\cap\pd\R_{+}^d$, and $B_{(\rho-\rho_1)/{10}}(x_0)$ for interior point $x_0\in \C^+_{\rho_1}$ satisfying $\hbox{dist}(x_0,\pd\R_{+}^d)>(\rho-\rho_1)/5$), and noting that the number of overlapping is less than $100^{d}$ if we double their radius (this kind of covering is easy to build), we achieve \eqref{limit4}.
\end{proof}

Now we are in position to prove Theorem \ref{thm_navierBC}, which follows from the lemma below.
\begin{lem}\label{lem4-NavierBC}
\sl{
Under assumptions in Lemma \ref{lem1-navierBC},  there exists some $\rho_0\ll 1$ such that for $\rho< \rho_0$, we have $\psi_1(t)\cdot\nabla T^{-1} \bm{u}\in L^{q,r}\big(\C^+_{\rho}\times(-1,0)\big)$ and 
\begin{align}\label{limit6}
&\quad\| \psi_1(t)\cdot\nabla T^{-1}\bm{u}\|_{L^{q,r}\left(\C^+_{\rho}\times(-1,0)\right)}\lec \rho^{-1-2d-\frac{d}{\min(q,r)}} \left(\|\bm{u},\CF\|_{L^{q,r}\left(Q^+_1\right)}+\|\bm{f}\|_{L^{q_*,r}\left(Q^+_1\right)}\right).
\end{align}
}
\end{lem}
\begin{proof}
Without loss of generality, we  may assume $q\leq r$ (the  case $q> r$ is similar). For some $\rho_0$ whose value will be specified later (and smaller than the $\rho_0$ in Lemma \ref{lem3-NavierBC}), we let $\rho_1<\rho<\rho_0$ in the cut-off function $\zeta$ of \eqref{defv1}. Analogously, we get lower order terms
\begin{equation*}
\|\int_{-1}^{t}e^{(t-s)\Delta} H_1\ds\|_{L^{q,r}\left(\C_{2\rho}\times(-1,0)\right)}\lec (\rho-\rho_1)^{-3} \|\bm{u}\|_{L^{q,r}\left(Q^+_1\right)},
\end{equation*}
 all higher order terms
\begin{align*}
\|\int_{-1}^{t}e^{(t-s)\Delta} H_7\ds\|_{\cW^{q,r}}\lec  \|\nabla\gamma\|_{L^{\infty}\left(B_{\rho}'\right)}\cdot\|\nabla \bm{v}\|_{\cH^{q,r}}+(\rho-\rho_1)^{-1}\|\bm{u}\|_{L^{q,r}\left(Q^+_1\right)},
\end{align*}
and external-force terms
\begin{equation*}
\|\int_{-1}^{t}e^{(t-s)\Delta}  H_6\ds\|_{L^{q,r}\left(\C_{2\rho}\times(-1,0)\right)}\lec (\rho-\rho_1)^{-2} \left(\|\bm{f}\|_{L^{q_*,r}\left(Q^+_1\right)}+\|\CF\|_{L^{q,r}\left(Q^+_1\right)}\right).
\end{equation*}
Whereas for $H_2$, by using Lemma \ref{lem-C-Z est}, we get
\begin{align*}
&\ \quad \|\int_{-1}^{t}e^{(t-s)\Delta} H_2\ds\|_{L^{q,r}\left(\C_{2\rho}\times(-1,0)\right)}\lec (\rho-\rho_1)^{-2}\|\psi_1(t)\cdot T^{-1}\nabla\bm{u}\|_{L^{q,\frac{dqr}{(d-1)q+r}}\left(\C^+_{\rho}\times (-1,0)\right)}\notag\\[6pt]
&\lec (\rho-\rho_1)^{-2}\|\psi_1(t)\cdot T^{-1}\nabla\bm{u}\|_{L^{q,q}\left(\C^+_{\rho}\times (-1,0)\right)}^{\frac1d}\cdot\|\psi_1(t)\cdot T^{-1}\nabla\bm{u}\|_{L^{q,r}\left(\C^+_{\rho}\times (-1,0)\right)}^{1-\frac1d},
\end{align*}
where it's easy to check that $\frac{(d-1)q+r}{dqr}-\frac{1}{r}<\frac{1}{2}$. Analogous to the estimates in  Step 2 of the proof for Lemma \ref{lem1-navierBC}, we get that, by using \eqref{limit3},
\begin{align*}
\|\psi_1(t)\cdot\nabla T^{-1} \bm{u}\|_{L^{q,r}\left(\C^+_{\rho_1}\times(-1,0)\right)}&\leq C\rho^{-1-\frac{d}{q}} (\rho-\rho_1)^{-2d} \left(\|\bm{u},\CF\|_{L^{q,r}\left(Q^+_1\right)}+\|\bm{f}\|_{L^{q_*,r}\left(Q^+_1\right)}\right)\\[6pt]
&\quad\ +\frac{1}{2}\|\psi_1(t)\cdot T^{-1}\nabla\bm{u}\|_{L^{q,r}\left(\C^+_{\rho}\times (-1,0)\right)}.
\end{align*}
Then applying Lemma \ref{lemA2}, we get that for $\rho<\frac{\rho_0}{2}$,
\begin{equation*}
 \|\psi_1(t)\cdot\nabla T^{-1} \bm{u}\|_{L^{q,r}\left(\C^+_{\rho}\times(-1,0)\right)}\lec\rho^{-1-2d-\frac{d}{q}} \left(\|\bm{u},\CF\|_{L^{q,r}\left(Q^+_1\right)}+\|\bm{f}\|_{L^{q_*,r}\left(Q^+_1\right)}\right),
\end{equation*}
which gives rise to \eqref{limit6}.
\end{proof}

\section{Proof of Theorem \ref{thm2_navierBC}}\label{sec5}

In this section we give the proof of Theorem \ref{thm2_navierBC}. Since any $\bm{u}$ in a weak solution pair is a weak solution, using the gradient estimate in Section \ref{sec4}, and setting $\CF=0$, we have $\nabla \bm{u}\in L^{q,r}\left(B^+_{1/2}\times(-\frac{3}{4},0)\right)$ and
\begin{equation}\label{gra2-1}
\| \nabla\bm{u}\|_{L^{q,r}\left(B^+_{1/2}\times(-\frac{3}{4},0)\right)}\lec \|\bm{u}\|_{L^{q,r}\left(Q^+_1\right)}+\|\bm{f}\|_{L^{q,r}\left(Q^+_1\right)}.
\end{equation}

\noindent Therefore, we focus on the estimates of higher order derivatives, which are divided into two parts. Readers can find the second derivative estimate in Lemma \ref{lem51} and third derivative estimate in Lemma \ref{lem52}. They imply Theorem \ref{thm2_navierBC}.

\subsection{Stokes equations of the approximations of velocity}

Let $\bm{v}$ be defined as in \eqref{defv1} and 
\begin{equation*}
\pi(x,t)=T^{-1}p(x,t)\cdot\zeta(x,t),
\end{equation*}
where $\tau\ge2,\rho,\rho_1$ in the definition \eqref{def-zeta} of $\zeta$ will be specified later.
At first, we derive the Stokes equation of $E^{\#}_{\ep_1,\ep_2,\ep_3}(\bm{v})$ (defined in \eqref{def-v}) and $E^{+}_{\ep_1,\ep_2,\ep_3}\left(\pi\right)$.

\smallskip

\noindent \textbf{\textit{Case 1}}:\quad The equation of $E^{+}_{\ep_1,\ep_2,\ep_3}(v_k)$ for $1\leq k\leq d-1$.

\smallskip

Substitute the test functions (which are tangential to boundary but not divergence free)
\[
T^{-1}\Phi(y,s)=\xi^+(x,y,t,s)\cdot(0,\cdots,0,\underbrace{1}_{k-th},0,\cdots,0,\underbrace{\gamma_{k}'}_{d-th})
\]
to \eqref{weakform-pair}. 
Recall that for $t$ close to $0^-$ we use the argument in Remark \ref{ep0-argument}.
Using \eqref{Stokes-curved-term2} and divergence-free property, we get
\begin{equation}\label{navier1}
(\pd_t-\Delta)E^{+}_{\ep_1,\ep_2,\ep_3}(v_k)+\pd_{x_{k}}E^{+}_{\ep_1,\ep_2,\ep_3}\left(\pi\right) =J_k\eqdefa J_{k1}+J_{k2}+J_{k3}+J_{k4},
\end{equation}
with 
\begin{align*}
J_{k1}&\eqdefa2\int_{\R^d_{+}\times\R} T^{-1}u_k\cdot\nabla \zeta\cdot\nabla_{y}\eta_{\ep_1,\ep_2,\ep_3}^{+}(x,y,t-s)\dy\ds\notag\\[6pt]
&\quad\,+\int_{\R^d_{+}\times\R} T^{-1}u_k\cdot(\pd_s+\Delta)\,\zeta\cdot\eta_{\ep_1,\ep_2,\ep_3}^+(x,y,t-s)\dy\ds\notag \\[6pt]
&=-2\int_{\R^d_{+}\times\R} \nabla_{y}T^{-1}u_k\cdot\nabla \zeta\cdot\eta_{\ep_1,\ep_2,\ep_3}^{+}(x,y,t-s)\dy\ds\notag\\[6pt]
&\quad\,+\int_{\R^d_{+}\times\R} T^{-1}u_k\cdot(\pd_s-\Delta)\,\zeta\cdot\eta_{\ep_1,\ep_2,\ep_3}^+(x,y,t-s)\dy\ds,
\end{align*}
(The second equality above transfers $M_2$ term in $J_{k1}$ to $M_1$ term, see \eqref{M1}-\eqref{M2}. It is desired although not necessary.)
\begin{align*}
J_{k2}&\eqdefa \int_{\R^d_{+}\times\R}T^{-1}u_d\cdot\gamma_{k}'\cdot\pd_{s}\xi^{+}-\nabla_{y}T^{-1}u_d\cdot\nabla_{y}\left(\gamma_{k}'\cdot\xi^{+}\right)\dy\ds,
\end{align*}
\begin{align*}
J_{k3}&\eqdefa\sum_{i=1}^{d-1}\int_{\R^d_{+}\times\R}\gamma_{i}'\cdot\left(\pd_{y_{i}}T^{-1}u_k\cdot\pd_{y_{d}}\xi^{+}+\pd_{y_{i}}T^{-1}u_d\cdot\pd_{y_{d}}\left(\gamma_k'\cdot\xi^{+}\right)\right)\dy\ds\notag\\[6pt]
&\quad\ +\sum_{i=1}^{d-1}\int_{\R^d_{+}\times\R}\gamma_{i}'\cdot\left(\pd_{y_{d}}T^{-1}u_k\cdot\pd_{y_{i}}\xi^{+}+\pd_{y_{d}}T^{-1}u_d\cdot\pd_{y_{i}}\left(\gamma_k'\cdot\xi^{+}\right)\right)\dy\ds\notag\\[6pt]
&\quad\ -\int_{\R^d_{+}\times\R}|\nabla\gamma|^2\cdot\left(\pd_{y_{d}}T^{-1}u_k\cdot\pd_{y_{d}}\xi^{+}+\pd_{y_{d}}T^{-1}u_d\cdot\pd_{y_{d}}\left(\gamma_k'\cdot\xi^{+}\right)\right)\dy\ds\notag\\[6pt]
&\quad\ -\sum_{i=1}^{d-1}\int_{\R^d_{+}\times\R}\pd_{y_{d}}\left(T^{-1}u_i\cdot\gamma_{ik}^{\prime\prime}\cdot\xi^{+}\right)\dy\ds,
\end{align*}
\begin{align*}
J_{k4}&\eqdefa\int_{\R^d_{+}\times\R}T^{-1}p\cdot\pd_{k}\zeta\cdot\eta_{\ep_1,\ep_2,\ep_3}^+(x,y,t-s)\dy\ds\notag\\[6pt]
&\quad\ +\alpha\int_{\R^d_{+}\times\R}\pd_{y_{d}}\left(\left(T^{-1}u_k+T^{-1}u_d\cdot\gamma_{k}\right)\cdot\sqrt{1+|\nabla\gamma|^2}\cdot\xi^{+}\right)\dy\ds\notag\\[6pt]
&\quad\ +\int_{\R^d_{+}\times\R}\left(T^{-1}f_{k}+T^{-1}f_{d}\cdot\gamma_{k}'\right)\cdot\xi^{+}\dy\ds.
\end{align*}
The last term of $J_{k3}$ comes from the last term of \eqref{Stokes-curved-term2} 
using
\begin{align*}
\sum_{1\leq i,\ell\leq d}T^{-1}(\pd_{y_\ell}u_{i})\cdot T^{-1}(\pd_{y_i}\Phi_{\ell})=
\sum_{1\leq i\leq d}
\bket{\pd_{y_k}T^{-1}\tilde{u}_i\cdot \pd_{y_i}\xi^++\pd_{y_{d}}\left(T^{-1}u_i\cdot\gamma_{ik}^{\prime\prime}\cdot\xi^{+}\right)}, 
\end{align*}
and that $T^{-1}\tilde{\bm{u}}$ is solenoidal in $\R^d_+$.
\smallskip

\noindent \textbf{\textit{Case 2}}:\quad The equation of $E^{-}_{\ep_1,\ep_2,\ep_3}(v_d)$.

\smallskip

By substituting the test functions (zero on boundary but not divergence free)
\[
T^{-1}\Phi(y,s)=\xi^-(x,y,t,s)\cdot(-\nabla\gamma,1)
\]
to \eqref{weakform-pair}, and using \eqref{Stokes-curved-term2} and divergence-free property, we get
\begin{equation}\label{navier2}
(\pd_t-\Delta)E^{-}_{\ep_1,\ep_2,\ep_3}(v_d)+\pd_{x_{d}}E^{+}_{\ep_1,\ep_2,\ep_3}\left(\pi\right)=J_d\eqdefa J_{d1}+J_{d2}+J_{d3},
\end{equation}
with
\begin{align*}
J_{d1}&\eqdefa2\int_{\R^d_{+}\times\R} T^{-1}\tilde{u}_d\cdot\nabla \zeta\cdot\nabla_{y}\eta_{\ep_1,\ep_2,\ep_3}^{-}(x,y,t-s)\dy\ds\notag\\[6pt]
&\quad\,+\int_{\R^d_{+}\times\R} T^{-1}\tilde{u}_d\cdot(\pd_s+\Delta)\,\zeta\cdot\eta_{\ep_1,\ep_2,\ep_3}^-(x,y,t-s)\dy\ds\notag\\[6pt]
&=-2\int_{\R^d_{+}\times\R} \nabla_{y}T^{-1}\tilde{u}_d\cdot\nabla \zeta\cdot\eta_{\ep_1,\ep_2,\ep_3}^{-}(x,y,t-s)\dy\ds\notag\\[6pt]
&\quad\,+\int_{\R^d_{+}\times\R} T^{-1}\tilde{u}_d\cdot(\pd_s-\Delta)\,\zeta\cdot\eta_{\ep_1,\ep_2,\ep_3}^-(x,y,t-s)\dy\ds,
\end{align*}
(The second equality in $J_{d1}$ has the same comment for $J_{k1}$.)
\begin{align*}
J_{d2}&\eqdefa \sum_{1\leq i,j\leq d-1}\int_{\R^d_{+}\times\R}\pd_{y_{i}}T^{-1}u_j\cdot\gamma_{ij}''\cdot\xi^{-}-T^{-1}u_j\cdot\gamma_{ij}''\cdot\pd_{y_{i}}\xi^{-}\dy\ds\notag\\[6pt]
&\quad\ +\sum_{i=1}^{d-1}\int_{\R^d_{+}\times\R}\gamma_{i}'\cdot\left(\pd_{y_{i}}T^{-1}\tilde{u}_d\cdot\pd_{y_{d}}\xi^{-}+\pd_{y_{d}}T^{-1}\tilde{u}_d\cdot\pd_{y_{i}}\xi^{-}\right)\dy\ds\notag\\[6pt]
&\quad\ -2\sum_{1\leq i,j\leq d-1}\int_{\R^d_{+}\times\R}\gamma_{i}'\cdot\gamma_{ij}''\cdot \pd_{y_{d}}T^{-1}u_j\cdot \xi^{-}\dy\ds\notag\\[6pt]
&\quad\ -\int_{\R^d_{+}\times\R}|\nabla\gamma|^2\cdot\pd_{y_{d}}T^{-1}\tilde{u}_d\cdot\pd_{y_{d}}\xi^{-}\dy\ds,
\end{align*}
\begin{align*}
J_{d3}&\eqdefa\int_{\R^d_{+}\times\R}T^{-1}p\cdot\pd_{d}\zeta\cdot\eta_{\ep_1,\ep_2,\ep_3}^-(x,y,t-s)\dy\ds\notag\\[6pt]
&\quad\ -\sum_{i=1}^{d-1}\int_{\R^d_{+}\times\R}T^{-1}p\cdot(\pd_{y_{i}}-\gamma_{i}'\cdot\pd_{y_{d}})\left(\gamma_{i}'\cdot\xi^{-}\right)\dy\ds \notag\\[6pt]
&\quad\ +\int_{\R^d_{+}\times\R}\left(T^{-1}f_{d}-\sum_{i=1}^{d-1}T^{-1}f_i\cdot\gamma_{i}'\right)\cdot\xi^{-}\dy\ds.
\end{align*}
The last term in \eqref{Stokes-curved-term2} disappears since 
\begin{align*}
&\,\quad\int_{\R^d_{+}\times\R}\sum_{1\leq k,\ell\leq d}T^{-1}\pd_{y_\ell}u_{k}\cdot T^{-1}\pd_{y_k}\Phi_{\ell}\dy\ds=\int_{\Omega\times \R}\sum_{1\leq k,\ell\leq d}\pd_{y_\ell}u_{k}\cdot \pd_{y_k}\Phi_{\ell}\dy \ds\\
&=\int_{\Omega\times \R}\sum_{1\leq \ell\leq d}\pd_{y_\ell}\div \bm{u}\cdot\Phi_{\ell}\dy ds=0, 
\end{align*}
where we have used $\Phi|_{\pd\Omega}=0$ and $\div \bm{u}=0$.

Notice that 
\begin{equation}\label{div4}
\div E^{\#}_{\ep_1,\ep_2,\ep_3}(\bm{v})=J_{0}\eqdefa E^{+}_{\ep_1,\ep_2,\ep_3}\left(T^{-1}\tilde{\bm{u}}\cdot\nabla\zeta\right). 
\end{equation}
From the Stokes equations \eqref{navier1}, \eqref{navier2} and \eqref{div4} in $\R^d$, we have
\begin{equation}\label{pressure1}
-\Delta E_{\ep_1,\ep_2,\ep_3}^{+}\left(\pi\right)=\pd_{t} J_{0}-\Delta J_{0}-\div \bm{J}, 
\end{equation}
where $\bm{J}=(J_1,\ldots,J_d)$ is the RHS of \eqref{navier1} and \eqref{navier2},
and hence
\begin{equation}\label{velocity1}
(\pd_{t}-\Delta)E^{\#}_{\ep_1,\ep_2,\ep_3}(\bm{v})=(1+\nabla (-\Delta)^{-1}\div)\bm{J}-(\pd_t-\Delta)\nabla \left(-\Delta\right)^{-1} J_{0}. 
\end{equation}
The functions in \eqref{pressure1} and \eqref{velocity1} are  smooth in $\mathbb R^d\times(-1,0)$ due to mollification. Moreover, the cut-off function $\zeta$ makes them compactly supported in the spatial variables and vanish for $t$ close to $-1$. Hence we may apply the whole-space Newtonian potential and Duhamel's formula to obtain
\begin{equation}\label{vp}
\left\{
\begin{split}
E_{\ep_1,\ep_2,\ep_3}^{+}\left(\pi\right)&= \left(-\Delta\right)^{-1}\left(\pd_{t} J_{0}-\Delta J_{0}-\div \bm{J}\right),\\[6pt]
E^{\#}_{\ep_1,\ep_2,\ep_3}(\bm{v})&=\int_{-1}^{t}e^{(t-s)\Delta}\left(1+\nabla (-\Delta)^{-1}\div\right)\bm{J}\ds-\nabla \left(-\Delta\right)^{-1} J_{0}.
\end{split}
\right.
\end{equation}

\smallskip

Now we give some estimates of the terms in $\bm{J}=(J_1,\cdots,J_d)$.
All the terms in $\bm{J}$ can be classified into the following six types. The terms involving the velocity can be written in the form of
\begin{equation}
\pd_{t}^{b_1}\nabla_{x}^{b_2} E_{\ep_1,\ep_2,\ep_3}^{\pm}\left(\nabla^{b_3} T^{-1}\tilde{\bm{u}}\cdot \left(\nabla^2\gamma\right)^{b_4}\cdot g_3(\nabla\gamma)\cdot\pd_{t}^{b_5}\nabla^{b_6}\zeta\right)
\end{equation} 
with non-negative integers  $2b_1+b_2+b_3+b_4+2b_5+b_6\leq 2$, $b_3,b_4\leq 1$ and some smooth bounded function $g_3$. They are called higher order terms if $2b_1+b_2+b_3=2$, otherwise lower order terms. Higher order terms are always the perturbation terms such that $g_3(0)=0$, i.e., $|g_3(\nabla\gamma)|\lec |\nabla\gamma|$. 
\begin{enumerate}
\renewcommand{\theenumi}{\roman{enumi}}
\item Lower order terms $M_1$ for $b_1=b_2=0$.
By using \eqref{gra2-1}, we have 
\begin{equation} \label{M1}
\|M_1\|_{\cW^{q,r}}\lec \frac{1}{(\rho-\rho_1)^{2}} \|\bm{u}\|_{L^{q,r}\left(Q^+_1\right)}+\frac{1}{\rho-\rho_1}\|\bm{f}\|_{L^{q,r}\left(Q^+_1\right)}. 
\end{equation}
\item Lower order terms $M_2$ for $b_2=1$, e.g.~$\nabla E_{\ep_1,\ep_2,\ep_3}^{\pm}\left({\bm{v}}\cdot \nabla^2\gamma\right)
$ from $J_{k3}$.
By using Lemma \ref{lem-C-Z est} and \eqref{gra2-1}, and moving $\nabla'$ inside integral and integrating by parts if necessary, we have 
\begin{align}\label{M2}
\Big\|\nabla'\nabla(-\Delta)^{-1} M_2, \int_{-1}^{t}\nabla'\nabla e^{(t-s)\Delta} M_2\ds\Big\|_{\cW^{q,r}}\lec \frac{1}{(\rho-\rho_1)^{2}} \|\bm{u}, \bm{f}\|_{L^{q,r}\left(Q^+_1\right)}. 
\end{align}
Note that when $b_4=1$, the $\nabla'$ above may act on $\nabla^2\gamma$ and hence generate $\nabla^3\gamma$ in the RHS of the estimate. This is the reason we require $\Gamma \in C^{2,1}$ in Theorem \ref{thm2_navierBC}.
\item Higher order terms $M_3$  for $\pd_{t}E_{\ep_1,\ep_2,\ep_3}^{+}\left(\bm{v}\cdot g_{3}(\nabla\gamma)\right)$.
For these terms (only from $J_{k2}$), we have
\begin{equation}\label{M3}
\|M_3\|_{\cW^{q,r}}\lec  \|\nabla\gamma\|_{L^{\infty}\left(B_{\rho}'\right)}\cdot  \|\pd_{t}E_{\ep_3}\left(\bm{v}\right)\|_{\cH^{q,r}}. 
\end{equation}

\item Higher order terms $M_4$ for $b_2=b_3=1$, e.g.~$\nabla E_{\ep_1,\ep_2,\ep_3}^{+}\left(\nabla\bm{v}\cdot g_{3}(\nabla\gamma)\right)$.
By  Lemma \ref{lem-C-Z est}, the estimate of commutator \eqref{commu3} and \eqref{gra2-1}, we have 
\begin{align}\label{M4}
\Big\|\nabla'\nabla(-\Delta)^{-1}M_4, \int_{-1}^{t}\nabla'\nabla e^{(t-s)\Delta} M_4\ds\Big\|_{\cW^{q,r}}\lec \frac{1}{(\rho-\rho_1)^{2}}\|\bm{u},\bm{f}\|_{L^{q,r}\left(Q^+_1\right)}
\notag\\[6pt]
+\|\nabla\gamma\|_{L^{\infty}\left(B_{\rho}'\right)}\cdot \|\nabla_{x'}E^{\pm}_{\ep_1,\ep_2,\ep_3}\left(\nabla\bm{v}\right)\|_{\cW^{q,r}}. 
\end{align}

\item Terms $M_5$ involving pressure, e.g.~$\pd_d E^+_{\ep_1,\ep_2,\ep_3}( \pi \cdot g_{3}(\nabla\gamma))$. For these terms, by  Lemma \ref{lem-C-Z est} and the estimate of commutator \eqref{commu3}, we get
\begin{align}\label{M5}
\Big\|\nabla'\nabla(-\Delta)^{-1}M_5, \int_{-1}^{t}\nabla'\nabla e^{(t-s)\Delta} M_5\ds\Big\|_{\cW^{q,r}}\lec \frac{1}{\rho-\rho_1}\|p\|_{L^{q,r}\big(Q^+_1\big)}\notag\\[6pt]
+\|\nabla\gamma\|_{L^{\infty}\left(B_{\rho}'\right)}\cdot\|\nabla_{x'}E_{\ep_1,\ep_2,\ep_3}^{\pm}\left(\pi\right)\|_{\cW^{q,r}}. 
\end{align}

\item Terms $M_6$ involving $\bm{f}$. For these terms,  we have
\begin{equation}\label{M6}
\|M_6\|_{\cW^{q,r}}\lec \|\bm{f}\|_{L^{q,r}\left(Q^+_1\right)}. 
\end{equation}
\end{enumerate}

\subsection{Second derivative estimate}

This subsection gives the proof of second derivative estimates, similar to those in Section \ref{sec3.2}, \cite[Theorem 1.2]{CLT24}, and \cite[\S IV.5]{Gal11}.

\begin{lem}\label{lem51}
\sl{
Suppose that the  boundary portion $\Gamma \in C^{2,1}$,  $\bm{f}\in L^{q,r}(Q^+_1)$, $\CF=0$, and $(\bm{u},p)\in L^{q,r}(Q^+_1)$ is a weak solution pair  with Navier boundary condition (see Definition \ref{def:wsp}).
Then we have $\pd_t\bm{u},\nabla^2 \bm{u}, \nabla p\in L^{q,r}\left(B^+_{1/2}\times(-\frac12,0)\right)$, and
\begin{align}\label{eq5.15}
&\quad\ \|\pd_t\bm{u}, \nabla^2\bm{u},\nabla p\|_{L^{q,r}\left(B^+_{1/2}\times(-\frac12,0)\right)}\lec \|\bm{u}\|_{L^{q,r}\left(Q^+_1\right)}+\|p\|_{L^{q,r}\left(Q^+_1\right)}+\|\bm{f}\|_{L^{q,r}\left(Q^+_1\right)}
\end{align}
}
\end{lem}
\begin{proof}
Similar to the argument at the beginning of Section \ref{sec3}, we only prove the corresponding local estimate near origin. Let $\tau=2$ in $\zeta$.  The proof is divided into two steps: in step 1 we get \eqref{claim3} and \eqref{use1'}; in step 2 we deal with the estimate of $\|\pd_{t}E_{\ep_3}\left(\bm{v}\right)\|_{\cH^{q,r}}$.

\smallskip

\noindent \textbf{\textit{Step 1: }} We claim that
\begin{equation}\label{claim3}
\nabla^2 E_{\ep_3}(\bm{v}),\nabla E_{\ep_3}\left(\pi\right)\in  \cH^{q,r}.
\end{equation}

By using Lemma \ref{lem-C-Z est}, \eqref{gra2-1} and \eqref{div4}, we have
\begin{align}\label{2nd-est1}
 \|\nabla J_{0}\|_{\cW^{q,r}} 
\lec \frac{1}{(\rho-\rho_1)^2}\|\bm{u}\|_{L^{q,r}\left(Q^+_1\right)}+\frac{1}{\rho-\rho_1}\|\bm{f}\|_{L^{q,r}\left(Q^+_1\right)},
\end{align}
and by first equality of  \eqref{div4}
\begin{align}\label{2nd-est4}
\Big\|\nabla(-\Delta)^{-1}\pd_{t}J_{0}\|_{\cW^{q,r}}\lec  \|\pd_{t}E_{\ep_3}\left(\bm{v}\right)\|_{\cH^{q,r}}.
\end{align}
Moreover, it is easy to verify that, by using the estimates of $M_1$-$M_6$ in \eqref{M1}-\eqref{M6},
\begin{align} \label{2nd-est2}
&\Big\|\nabla'\nabla(-\Delta)^{-1} \bm{J},\, \int_{-1}^{t}\nabla'\nabla e^{(t-s)\Delta} \bm{J}\ds\Big\|_{\cW^{q,r}}\lec \frac{1}{(\rho-\rho_1)^2} \|\bm{u},p, \bm{f}\|_{L^{q,r}\left(Q^+_1\right)}\notag\\[6pt]
&\,\quad +\|\nabla\gamma\|_{L^{\infty}\left(B_{\rho}'\right)}\cdot  \left(\|\pd_{t}E_{\ep_3}\left(\bm{v}\right)\|_{\cH^{q,r}}+\|\nabla_{x'}E^{\pm}_{\ep_1,\ep_2,\ep_3}\left(\nabla\bm{v}\right),\nabla_{x'}E_{\ep_1,\ep_2,\ep_3}^{\pm}\left(\pi\right)\|_{\cW^{q,r}}\right).
\end{align}
Hence, in view of \eqref{vp}, and  using \eqref{2nd-est1}, \eqref{2nd-est4},  \eqref{2nd-est2}, we achieve
\begin{align*}
&\|\nabla_{x'}\nabla_{x}E^{\#}_{\ep_1,\ep_2,\ep_3}(\bm{v}),\, \nabla_{x'}E^{+}_{\ep_1,\ep_2,\ep_3}(\pi)\|_{\cW^{q,r}}\lec \frac{1}{(\rho-\rho_1)^2} \|\bm{u}, p, \bm{f}\|_{L^{q,r}\left(Q^+_1\right)}\notag\\[6pt]
&\quad\ +\|\pd_{t}E_{\ep_3}\left(\bm{v}\right)\|_{\cH^{q,r}}+\|\nabla\gamma\|_{L^{\infty}\left(B_{\rho}'\right)}\cdot  \|\nabla_{x'}\nabla_{x} E^{h}_{\ep_1,\ep_3}\left(\bm{v}\right),\nabla_{x'}E^{h}_{\ep_1,\ep_3}\left(\pi\right)\|_{\cH^{q,r}} ,
\end{align*}
which implies that, by taking $\ep_2\rightarrow0$, using \eqref{rho2}, and then taking $\ep_1\rightarrow0$, 
\begin{equation}\label{hori4}
\nabla'\nabla E_{\ep_3}(\bm{v}),\nabla' E_{\ep_3}\left(\pi\right)\in  \cH^{q,r},
\end{equation}
and
\begin{align}\label{hori5}
&\|\nabla'\nabla E_{\ep_3}(\bm{v}),\, \nabla' E_{\ep_3}(\pi)\|_{\cH^{q,r}}\lec \frac{1}{(\rho-\rho_1)^2} \|\bm{u},p, \bm{f}\|_{L^{q,r}\left(Q^+_1\right)}+\|\pd_{t}E_{\ep_3}\left(\bm{v}\right)\|_{\cH^{q,r}}. 
\end{align}

By taking $\ep_1,\ep_2\rightarrow0$ in \eqref{navier1} and \eqref{navier2} for $(x,t)\in \R_{+}^d\times(-1,0)$, we get that 
\begin{equation}\label{navier3}
(\pd_t-\Delta)E_{\ep_3}(\bm{v})+\nabla E_{\ep_3}\left(\pi\right) =\bar{\bm{J}},
\end{equation}
in the distribution $\mathcal{D}'\left(\R_{+}^d\times(-1,0)\right)$, where $\bar{\bm{J}}$ are formulas that $\bm{J}$ becomes after taking $\ep_1,\ep_2\rightarrow0$.
Note that we do not have $\pd_{d}^2E_{\ep_3}(\bm{v}),\pd_{d}E_{\ep_3}(\pi),  \bar{\bm{J}}$ in $\cH^{q,r}$ yet. Collecting $\pd_{d}^2E_{\ep_3}(\bm{v}),\pd_{d}E_{\ep_3}(\pi)$ terms in $\bar{\bm{J}}$ (from $\bar J_{k2}$, $\bar J_{k3,3}$, $\bar J_{d2,4}$ and $\bar J_{d3,2}$) and moving them to the LHS of \eqref{navier3},
we see that 
\begin{align}\label{galdi-trick1}
(1+|\nabla\gamma|^2)\cdot\left(\pd_d^2E_{\ep_3}(v_k)+\gamma'_k\pd_d^2E_{\ep_3}(T^{-1}u_d\cdot\zeta)\right)\in \cH^{q,r}, \quad (k\leq d-1), 
\end{align}
and
\begin{align}\label{galdi-trick2}
(1+|\nabla\gamma|^2)\cdot\left(-\pd_{d}^2E_{\ep_3}\left(v_d\right)+\pd_{d}E_{\ep_3}\left(\pi\right)\right)\in \cH^{q,r}. 
\end{align}
We give some details on how to obtain \eqref{galdi-trick1} and \eqref{galdi-trick2} from \eqref{navier3}. Let us illustrate the point by considering the term $\bar J_{k2}$. After letting $\varepsilon_1,\varepsilon_2\to0$, we have
\begin{align*}
\bar J_{k2}
&=
E_{\varepsilon_3}\left(T^{-1}u_d\cdot\gamma'_k\partial_s\zeta\right)
-\partial_tE_{\varepsilon_3}\left(T^{-1}u_d\cdot \gamma'_k\zeta\right)  \\
&\quad
-E_{\varepsilon_3}\left(\nabla_{y} T^{-1}u_d\cdot \nabla_{y}(\gamma'_k\zeta)\right)
+\sum_{i=1}^{d-1}\pd_{x_i}E_{\varepsilon_3}\left(\pd_{y_i} T^{-1}u_d\cdot \gamma'_k\zeta\right) \\
&\quad
+\gamma'_k\partial_d^2E_{\varepsilon_3}\left(T^{-1}u_d\cdot\zeta\right)
-\gamma'_k\partial_dE_{\varepsilon_3}\left(T^{-1}u_d\cdot\partial_d\zeta\right).
\end{align*}
All terms on the right-hand side belong to $\mathcal H^{q,r}$, except possibly $\gamma'_k\partial_d^2E_{\varepsilon_3}\left(T^{-1}u_d\cdot\zeta\right)$.
Thus this term has to be moved to the left-hand side when deriving \eqref{galdi-trick1}. The same inspection applies to the other terms in $\bar{\bm{J}}$ (from $\bar J_{k2}$, $\bar J_{k3,3}$, $\bar J_{d2,4}$ and $\bar J_{d3,2}$).

In view of \eqref{div4}, \eqref{2nd-est1} and \eqref{hori4}, we see that $\pd_{d}^2E_{\ep_3}\left(v_d\right)\in \cH^{q,r}$,  which together with
\eqref{galdi-trick2} leads to $\pd_{d}E_{\ep_3}\left(\pi\right)\in \cH^{q,r}$.  Multiplying \eqref{galdi-trick1} by $\gamma'_k$ and summing in $k$, we have
\begin{align*}
\sum_{k=1}^{d-1}\gamma'_k\pd_d^2E_{\ep_3}(v_k)+|\nabla\gamma|^2\pd_d^2E_{\ep_3}(T^{-1}u_d\cdot\zeta)\in \cH^{q,r},
\end{align*}
which, when added to $\pd_{d}^2E_{\ep_3}\left(v_d\right)\in \cH^{q,r}$, leads to $\pd_d^2E_{\ep_3}(T^{-1}u_d\cdot\zeta)\in \cH^{q,r}$, and this together with \eqref{galdi-trick1} leads to $\pd_d^2E_{\ep_3}(v_k)\in \cH^{q,r}$. Hence, we reach \eqref{claim3}. Together with \eqref{hori5}, we achieve
\begin{align}\label{use1'}
\|\nabla^2 E_{\ep_3}(\bm{v}),\,\nabla E_{\ep_3}(\pi)\|_{\cH^{q,r}}\lec \frac{1}{(\rho-\rho_1)^2} \|\bm{u},p,\bm{f}\|_{L^{q,r}\left(Q^+_1\right)}+\|\pd_{t}E_{\ep_3}\left(\bm{v}\right)\|_{\cH^{q,r}}. 
\end{align}

\bigskip

\noindent \textbf{\textit{Step 2: }}  The estimates of second order derivative.

\smallskip

By \eqref{use1'}, it suffices to derive a
bound of  $\|\pd_{t}E_{\ep_3}\left(\bm{v}\right)\|_{\cH^{q,r}}$. 
From \eqref{vp} we have
\begin{align}\label{2nd-est3}
 \pd_{t}E^{\#}_{\ep_1,\ep_2,\ep_3}(\bm{v})=\int_{-1}^{t}e^{(t-s)\Delta}\pd_{s}\left(1+\nabla (-\Delta)^{-1}\div\right)\bm{J}\ds -\nabla (-\Delta)^{-1}\pd_{t} J_{0}. 
\end{align}

We now estimate the first term in the RHS of \eqref{2nd-est3}. By virtue of $\eqref{M1}$, \eqref{M3}, \eqref{M6}, and Calderón--Zygmund estimate for the heat equation in the whole space, we get
\begin{align}\label{M136}
\Big\|\int_{-1}^{t}e^{(t-s)\Delta}\pd_{s}\left(M_1,\, M_3,\,M_6\right)\ds \Big\|_{\cW^{q,r}} &\lec \frac{1}{(\rho-\rho_1)^2} \|\bm{u}, \bm{f}\|_{L^{q,r}\left(Q^+_1\right)}\notag \\[6pt]
&\,\quad +\|\nabla\gamma\|_{L^{\infty}(B^{\prime}_\rho)}\cdot \|\pd_{t}E_{\ep_3}\left(\bm{v}\right)\|_{\cH^{q,r}}. 
\end{align}
For $M_5$, by using Lemma \ref{lem-C-Z est}, we get
\begin{equation}\label{M5'}
 \Big\|\int_{-1}^{t}e^{(t-s)\Delta}\pd_{s}M_5\ds\Big\|_{\cW^{q,r}}\lec \frac{1}{\rho-\rho_1}\|p\|_{L^{q,r}\big(Q^+_1\big)}+\|\nabla\gamma\|_{L^{\infty}(B^{\prime}_\rho)}\cdot\|\nabla E_{\ep_3}\left(\pi\right)\|_{\cH^{q,r}}. 
\end{equation}
Since the operator $\nabla (-\Delta)^{-1}\div$ is bounded in $L^q(\R^d)$, it can be included in the estimates \eqref{M136} and \eqref{M5'} for $M_1,M_3,M_5$ and $M_6$.
 
For $M_2$, since horizontal derivatives can be moved inside the mollification (as well as $\pd_d$ before $E^{+}_{\ep_1,\ep_2,\ep_3}$), some terms in $M_2$ can be transferred to $M_1$ (though $\nabla^3\gamma$ might be generated, it will not affect any  previous estimates). The only kind of term that requires special treatment is
\begin{equation}\label{M2''}
\pd_{t}M_2= \pd_{x_{d}}E^{h,-}_{\ep_1,\ep_2}\left(\pd_{t}E_{\ep_3}\left(T^{-1}\tilde{\bm{u}}\cdot\nabla^{b_6}\zeta\right)\cdot\left(\nabla^2\gamma\right)^{b_4}\cdot g_3(\nabla\gamma)\right).
\end{equation} 
Replacing the cut-off function $\zeta$ with $ \nabla^{b_6}\zeta$ in \eqref{navier3}, we have
\begin{align}\label{pdt1}
 \pd_{t}E_{\ep_3}\left(T^{-1}\tilde{\bm{u}}\cdot\nabla^{b_6}\zeta\right)=g_{4}(\nabla\gamma)\cdot\left(\Delta E_{\ep_3}\left(T^{-1}\tilde{\bm{u}}\cdot\nabla^{b_6}\zeta\right)-\nabla E_{\ep_3}\left(T^{-1}p\cdot\nabla^{b_6}\zeta\right)+\cdots\right).
\end{align}
The above $g_4(\nabla \gamma)$ is 
originally 1. Denote by $\bar{M}_j$ the limit of those new $M_j$ with $\nabla^{b_6}\zeta$ after taking $\ep_1,\ep_2\rightarrow0$, $j=1,\ldots, 6$. We move $\bar{M}_3$ to the LHS and can solve \eqref{pdt1} with $g_4(\nabla \gamma)$ a matrix close to identity, and RHS has no $\bar{M}_3$ anymore.  
By substituting \eqref{pdt1} into \eqref{M2''} and moving one derivative out of $E^{h,-}_{\ep_1,\ep_2}$ if necessary (the above two terms in the RHS of \eqref{pdt1}, and $\bar{M}_2,\bar{M}_4,\bar{M}_5$),  we get
by Lemma \ref{lem-C-Z est},
\begin{align}\label{M2'}
&\quad\ \|\int_{-1}^{t}e^{(t-s)\Delta}\pd_{s}M_2\ds,\int_{-1}^{t}e^{(t-s)\Delta}\nabla^2(-\Delta)^{-1}\pd_{s}M_2\ds\|_{L^{q,r}\left(\C(2\rho)\times(-1,0)\right)}\notag\\[6pt]
&\lec \frac{1}{(\rho-\rho_1)^2}\|\bm{u},p, \bm{f}\|_{L^{q,r}\left(Q^+_1\right)}.
\end{align}

For $M_4$, the terms inside $M_4$ can be divided into three categories. The first and second types contain at least one horizontal derivative,
\begin{align*}
M_{41}&=\nabla'E_{\ep_1,\ep_2,\ep_3}^{\pm}\left(\nabla T^{-1}\tilde{\bm{u}} \cdot g_3(\nabla\gamma)\cdot\zeta\right),\\
M_{42}&=\nabla E_{\ep_1,\ep_2,\ep_3}^{\pm}\left(\nabla'T^{-1}\tilde{\bm{u}} \cdot g_3(\nabla\gamma)\cdot\zeta\right).
\end{align*}
The third type has two vertical derivatives,
\begin{equation*}
M_{43}=\pd_d E_{\ep_1,\ep_2,\ep_3}^{-}\left(\pd_d T^{-1}\tilde{\bm{u}}\cdot g_3(\nabla\gamma)\cdot \zeta\right) \ \text{ or } \  \pd_d E_{\ep_1,\ep_2,\ep_3}^{+}\left(\pd_dT^{-1}\tilde{u}_d\cdot|\nabla\gamma|^2\cdot \zeta\right).
\end{equation*}
(There are no terms of the type $\pd_d E^{+}\left(\pd_d T^{-1}{u_k}\cdot g_3(\nabla\gamma)\cdot \zeta\right)$, $k<d$.)
For $M_{41}$, we can move the horizontal derivative inside the mollification; for
$M_{42}$ and $M_{43}$, we can move the inside derivative out of the mollification (as $\eta^-_{\ep_1,\ep_2,\ep_3}=T^{-1}\tilde u_d=0$ when $y_d=0$). There will be some lower order terms $M_1,M_2$ generated by this process.
Hence
\begin{align}\label{M4'}
&\quad\ \|\int_{-1}^{t}e^{(t-s)\Delta}\pd_{s}M_4\ds,\int_{-1}^{t}e^{(t-s)\Delta}\nabla^2(-\Delta)^{-1}\pd_{s}M_4\ds\|_{L^{q,r}\left(\C(2\rho)\times(-1,0)\right)}\notag\\[6pt]
&\lec\frac{1}{(\rho-\rho_1)^2}\|\bm{u},p, \bm{f}\|_{L^{q,r}\left(Q^+_1\right)}+\|\nabla\gamma\|_{L^{\infty}(B^{\prime}_\rho)}\cdot\|\nabla'\nabla E_{\ep_3}\left(\bm{v}\right),\pd_{t}E_{\ep_3}\left(\bm{v}\right)\|_{\cH^{q,r}},
\end{align}
where $\|\nabla'\nabla E_{\ep_3}\left(\bm{v}\right)\|_{\cH^{q,r}}$ is generated by $M_{41}$, and $\|\pd_{t}E_{\ep_3}\left(\bm{v}\right)\|_{\cH^{q,r}}$ by $M_{42}$ and $M_{43}$.

By combining \eqref{M136}, \eqref{M5'}, \eqref{M2'}, and  \eqref{M4'}, we get
\begin{align}\label{add2}
\Big\|\int_{-1}^{t}e^{(t-s)\Delta}\pd_{s}\left(1+\nabla (-\Delta)^{-1}\div\right)\bm{J}\ds
\Big\|_{L^{q,r}\left(\C(2\rho)\times(-1,0)\right)}\lec \frac{1}{(\rho-\rho_1)^2} \|\bm{u},p,\bm{f}\|_{L^{q,r}\left(Q^+_1\right)} \notag\\[6pt]
+\|\nabla\gamma\|_{L^{\infty}(B^{\prime}_\rho)}\cdot\|\pd_{t}E_{\ep_3}\left(\bm{v}\right),\nabla'\nabla E_{\ep_3}\left(\bm{v}\right),\nabla E_{\ep_3}\left(\pi\right)\|_{\cH^{q,r}}. 
\end{align}

Now we estimate the last term $\nabla (-\Delta)^{-1}\pd_{t} J_{0}$ in the RHS of \eqref{2nd-est3}. In view of \eqref{div4} and \eqref{pdt1}, we have
\begin{align}\label{pd_t-to-pdx}
\begin{split}
\pd_tJ_{0}&= \sum_{k=1}^d E^{h,+}_{\ep_1,\ep_2}\bke{\pd_t E_{\ep_3} \left(T^{-1}\tilde u_k\cdot\pd_k\zeta\right)}
\\
&=\sum_{k,j=1}^d E^{h,+}_{\ep_1,\ep_2}\left(g_{4}^{kj}(\nabla\gamma)\left(\Delta E_{\ep_3}\left(T^{-1}\tilde{u}_j\cdot\pd_k\zeta\right)-\pd_j E_{\ep_3}\left(T^{-1}p\cdot\pd_k\zeta\right)+\sum _{3\not=\ell\le 6} \bar M_{\ell,j}^{(k)}\right)\right).  
\end{split}
\end{align}
Here $g_{4}^{kj}$ means the $kj$ component of the matrix $g_4$, and
$\bar M_{\ell,j}^{(k)}$ means the $j$-th component of $\bar M_{\ell}$ with $\zeta$ replaced by $\pd_k \zeta$.
For any term in $\pd_tJ_{0}$ with a tangential derivative, we can pull the tangential derivative out of convolution and bound the resulting term with Lemma \ref{lem-C-Z est}. Hence we only need to bound $\nabla (-\Delta)^{-1}N$ in $L^{q,r}$ for
terms $N$ of the form
\[
N=E^{h,+}_{\ep_1,\ep_2}\left(g(\nabla\gamma)\pd_d W\right), \quad W=\pd_d E_{\ep_3}\left(T^{-1}\tilde{u}_j\cdot\pd_k\zeta\right) \ \text{or } \ W=E_{\ep_3}\left(T^{-1}p\cdot\pd_k\zeta\right).
\]
By integration by parts,
\begin{align*}
N&=\pd_d E^{h,-}_{\ep_1,\ep_2}\left(g(\nabla\gamma) W\right)
- \int_{\R^{d-1}} g(\nabla\gamma) W(y',0) \eta_{\ep_1}^h(x'-y') 2\eta_{\ep_2}(x_d)\,dy'=N_1+N_2,
\end{align*}
where $N_1$ can be estimated as other terms, and we only need to estimate the boundary term $N_2$. By the inequality
\begin{equation}\label{eq5.34}
\|\nabla(-\Delta)^{-1}g\|_{L^{q}\left(\C_{2\rho}\right)} 
\lesssim \norm{ (|x|^{1-d})*|g|}_{L^{q}\left(\C_{2\rho}\right)} \lec A\|g\|_{L^{1}\left((-2\rho,2\rho);L^{q}\left(B_{2\rho}'\right)\right)}
\end{equation}
for $g\in L^{q}_{c}\left(\C_{2\rho}\right)$, where $A=  \norm{|x|^{1-d}}_{L^{q}\left((-8\rho,8\rho);L^{1}\left(B_{8\rho}'\right)\right)} \lec 1$,
we get
\begin{align*}
\norm{\nabla (-\Delta)^{-1}N_2}_{L^{q,r}(B_{2\rho}'\times(-1,0))} & \lec \norm{\eta_{\ep_2}}_{L^1(\R)} \norm{ E_{\ep_1}^h ( g(\nabla\gamma) W|_{y_d=0}) }_{L^{q,r}(B_{2\rho}'\times(-1,0))} 
\\
&\lec  \norm{ W|_{y_d=0} }_{L^{q,r}(B_{2\rho}'\times(-1,0))} 
\\
& \lec  \frac{c_\ep}{(\rho-\rho_1)^{\frac{1}{q-1}}} \|W\|_{L^{q,r}\left(Q^+_1\right)}+ \ep(\rho-
\rho_1)\|\pd_{d}W\|_{L^{q,r}\left(Q^+_1\right)},
\end{align*}
where the last inequality is 
by \eqref{bad1}, and $\ep>0$ is a small constant to be chosen. We conclude
\begin{align}\label{add3}
&\Big\|\nabla(-\Delta)^{-1}\pd_{t}J_0\Big\|_{L^{q,r}\left(\C_{2\rho}\times(-1,0)\right)}\lec \frac{c_\ep}{(\rho-\rho_1)^{3+\frac{1}{q-1}}} \|\bm{u},p,\bm{f}\|_{L^{q,r}\left(Q^+_1\right)}+
 \ep{\mathcal A}_\rho,
\end{align}
where
\[
{\mathcal A}_\rho = \left\|\pd_{d}^2E_{\ep_3}\left(T^{-1}\bm{u}\cdot\psi_{2}(t)\right),\pd_{d}E_{\ep_3}\left(T^{-1}p\cdot\psi_{2}(t)\right)\right\|_{L^{q,r}\left(\C_{\rho}^+\times(-1,0)\right)}.
\]

By combining  \eqref{2nd-est3} and the estimates
 \eqref{add2}, \eqref{add3}, we have
\begin{align*}
&\|\pd_{t}E^{\#}_{\ep_1,\ep_2,\ep_3}(\bm{v})\|_{L^{q,r}\left(\C(2\rho)\times(-1,0)\right)}\lec\frac{c_\ep}{(\rho-\rho_1)^{3+\frac{1}{q-1}}}\|\bm{u},p, \bm{f}\|_{L^{q,r}\left(Q^+_1\right)}\notag\\[6pt]
&\,\quad +\|\nabla\gamma\|_{L^{\infty}(B^{\prime}_\rho)}\cdot\|\pd_{t}E_{\ep_3}\left(\bm{v}\right),\nabla'\nabla E_{\ep_3}\left(\bm{v}\right),\nabla E_{\ep_3}\left(\pi\right)\|_{\cH^{q,r}} +\ep{\mathcal A}_\rho.
\end{align*}
Hence, by taking small $\rho$, we achieve
\begin{align}\label{use2}
 &\|\pd_{t}E_{\ep_3}(\bm{v})\|_{\cH^{q,r}} \lec\frac{c_\ep}{(\rho-\rho_1)^{3+\frac{1}{q-1}}}\|\bm{u},p,\bm{f}\|_{L^{q,r}\left(Q^+_1\right)}\notag\\[6pt]
 &\,\quad +\|\nabla\gamma\|_{L^{\infty}(B^{\prime}_\rho)}\cdot\|\nabla'\nabla E_{\ep_3}\left(\bm{v}\right),\nabla E_{\ep_3}\left(\pi\right)\|_{\cH^{q,r}}+ \ep{\mathcal A}_\rho.
\end{align}
It together with \eqref{use1'} implies that, by taking small $\rho$ and $\ep$,
\begin{align*}
&\|\pd_{t}E_{\ep_3}(\bm{v}), \nabla^2 E_{\ep_3}(\bm{v}),\nabla E_{\ep_3}(\pi)\|_{\cH^{q,r}}\leq \frac{C}{(\rho-\rho_1)^{3+\frac{1}{q-1}}}\|\bm{u},p,\bm{f}\|_{L^{q,r}\left(Q^+_1\right)} + \frac{1}{100}{\mathcal A}_\rho. 
\end{align*}
Hence we get \eqref{add4} for any $\rho_1<\rho<\rho_0\ll1$. Finally, \eqref{add5} is  derived  by  the same arguments.  This completes the proof of Lemma \ref{lem51}.
\end{proof}

\subsection{Third derivative estimate}
\begin{lem}\label{lem52}
\sl{
Under assumptions in Lemma \ref{lem51}, and assuming additionally $\Gamma\in C^{3,1}$ and $\nabla\bm{f}\in L^{q,r}\left(Q^+_1\right)$,  we have $\nabla\pd_t\bm{u}, \nabla^3 \bm{u}, \nabla^2 p\in L^{q,r}\left(Q^+_{1/2}\right)$, and
\begin{align}
\| \nabla\pd_t\bm{u}, \nabla^3\bm{u},  \nabla^2 p\|_{L^{q,r}\left(Q^+_{1/2}\right)} \lec \|\bm{u}\|_{L^{q,r}\left(Q^+_1\right)}+\|p\|_{L^{q,r}\left(Q^+_1\right)}+\|\bm{f}\|_{L^r\left(-1,0;W^{1,q}(B^+_1)\right)}.
\end{align}
}
\end{lem}
\begin{proof} As in the proof for Lemma \ref{lem51},
we only present the local estimate near origin. Let $\zeta$ be cut-off function defined in \eqref{def-zeta} with
$\rho_1=\frac{\rho}{2}$ and $\tau=3$. The proof has two steps.

\smallskip
 
\noindent \textbf{\textit{Step 1: }} We first show spatial regularity of the  time-mollified solution. We claim that
\begin{equation}\label{claim4}
 \nabla^3 E_{\ep_3}(\bm{v}), \nabla^2 E_{\ep_3}\left(\pi\right)\in  \cH^{q,r}.
\end{equation}
In view of \eqref{vp}, and using Lemma \ref{lem-C-Z est}, Lemma \ref{lem51} and \eqref{gra2-1}, we get
\begin{align}\label{use4}
&\quad\,\|\nabla_{x'}^2\nabla_{x}E^{\#}_{\ep_1,\ep_2,\ep_3}(\bm{v})\|_{\cW^{q,r}}+\|\nabla_{x'}^2E^{+}_{\ep_1,\ep_2,\ep_3}(\pi)\|_{\cW^{q,r}}\notag\\[6pt]
&\lec \|\pd_{t}J_0,\nabla_{x'}^2 J_0, \nabla_{x'}^2\nabla_{x}(-\Delta)^{-1}\bm{J}\|_{\cW^{q,r}}\notag\\[6pt]
&\lec {\mathcal B}_\rho+\|\nabla\gamma\|_{L^{\infty}\left(B_{\rho}'\right)}\cdot  \|\nabla'\pd_{t}E_{\ep_3}\left(\bm{v}\right),\nabla_{x'}^2\nabla_{x} E^{h}_{\ep_1,\ep_3}\left(\bm{v}\right),\nabla_{x'}^2 E^{h}_{\ep_1,\ep_3}\left(\pi\right)\|_{\cH^{q,r}},
\end{align}
where the last line has no $E^{\pm}_{\ep_2}$ mollification and
\begin{equation}
{\mathcal B}_\rho= \rho^{-3} \|\bm{u}\|_{L^{q,r}\left(Q^+_1\right)}+\rho^{-2}\|p, \bm{f}, \nabla \bm{f}\|_{L^{q,r}\left(Q^+_1\right)}.
\end{equation}
The estimates of terms involving $\bm{J}$ are quite similar to $M_1$-$M_6$ in \eqref{M1}-\eqref{M6}, so we omit details here. By taking $\ep_2\rightarrow0$ in \eqref{use4} and using \eqref{rho2}, then taking $\ep_1\rightarrow0$, we achieve
\begin{equation}\label{eq5.41}
(\nabla')^2\nabla E_{\ep_3}(\bm{v}),\ (\nabla')^2 E_{\ep_3}(\pi)\in \cH^{q,r},
\end{equation}
and
\begin{equation}\label{use5}
\|(\nabla')^2\nabla  E_{\ep_3}(\bm{v})\|_{\cH^{q,r}}+\|(\nabla')^2 E_{\ep_3}(\pi)\|_{\cH^{q,r}}\lec {\mathcal B}_\rho +\|\nabla\gamma\|_{L^{\infty}\left(B_{\rho}'\right)}\cdot  \|\nabla'\pd_{t}E_{\ep_3}\left(\bm{v}\right)\|_{\cH^{q,r}}. 
\end{equation}

Now we go back to \eqref{navier3}. The arguments \eqref{navier3}--\eqref{use1'} in Step 1 of the proof of Lemma \ref{lem51} give
\begin{equation}\label{navier4}
\pd_{d}^2 E_{\ep_3}(\bm{v}),\ \pd_{d} E_{\ep_3}(\pi)=\hbox{span}\Big\{\nabla'\nabla E_{\ep_3}(\bm{v}),\ \nabla' E_{\ep_3}(\pi), \pd_t E_{\ep_3}(\bm{v})\Big\} + l.o.t.,
\end{equation}
where ``span'' means all linear combinations of elements in the set with coefficients $g(\nabla \gamma)$ being polynomials of $\nabla '\gamma$ and  $(1+|\nabla \gamma|^2)^{-1}$, and ``$l.o.t.$'' means lower order terms.
Taking horizontal derivative $\nabla'$ on both sides of \eqref{navier4}, we get
\[
\nabla'\pd_{d}^2 E_{\ep_3}(\bm{v}),\ \nabla'\pd_{d} E_{\ep_3}(\pi)=\hbox{span}\Big\{(\nabla')^2\nabla E_{\ep_3}(\bm{v}),\ (\nabla'^2) E_{\ep_3}(\pi), \nabla'\pd_t E_{\ep_3}(\bm{v})\Big\}+ l.o.t..
\]
Then using \eqref{eq5.41}, we have 
\begin{equation}\label{eq5.43}
\nabla' \pd_{d}^2 E_{\ep_3}(\bm{v}),\ \nabla'\pd_{d} E_{\ep_3}(\pi)\in \cH^{q,r}
\end{equation}
with their norms bounded by $ {\mathcal B}_\rho + \|\nabla\pd_{t}E_{\ep_3}\left(\bm{v}\right)\|_{\cH^{q,r}}$.
Then we take vertical derivative $\pd_{d}$ on both sides of \eqref{navier4}, and use \eqref{eq5.41}, \eqref{eq5.43} to get
\begin{equation}\label{galdi-trick4}
 \pd_{d}^3 E_{\ep_3}(\bm{v}),\ \pd_{d}^2 E_{\ep_3}(\pi)\in \cH^{q,r}
\end{equation}
with the same bound
\begin{equation}\label{eq0929}
\|\nabla^3  E_{\ep_3}(\bm{v})\|_{\cH^{q,r}}+\|\nabla^2 E_{\ep_3}(\pi)\|_{\cH^{q,r}}\lec {\mathcal B}_\rho + \|\nabla\pd_{t}E_{\ep_3}\left(\bm{v}\right)\|_{\cH^{q,r}}.
\end{equation}
Hence we get \eqref{claim4}.

\bigskip

\noindent \textbf{\textit{Step 2: }}  The estimates of third order derivatives. 

\smallskip

We will first derive the bound of  $\|\nabla\pd_{t}E_{\ep_3}\left(\bm{v}\right)\|_{\cW^{q,r}}$. In view of \eqref{div4}, we have
\begin{equation}\label{eq5.45}
\pd_{d}\pd_{t}E_{\ep_1,\ep_2,\ep_3}^{-}\left(v_d\right)=-\sum_{i=1}^{d-1}\pd_{i}\pd_{t}E_{\ep_1,\ep_2,\ep_3}^{+}\left(v_i\right)+\pd_{t}J_0. 
\end{equation}
Its LHS is less desired as we cannot move $\pd_d$ inside mollification. By \eqref{2nd-est3}, \eqref{eq5.45} and
Lemma \ref{lem-C-Z est} we get
\begin{align}\label{use6}
\|\nabla_{x}\pd_{t} E^{\#}_{\ep_1,\ep_2,\ep_3}(\bm{v})\|_{\cW^{q,r}}
&\lec \|\pd_{t}J_0\|_{\cW^{q,r}} +\bigg\|\int_{-1}^{t}e^{(t-s)\Delta}\pd_{s}\nabla_{x'}\bm{J} \ds\bigg\|_{\cW^{q,r}}\notag\\[6pt]
&\quad\ +\sum_{i=1}^{d-1}\bigg\|\int_{-1}^{t}e^{(t-s)\Delta}\pd_{s}\pd_{x_{d}} J_{i}\ds\bigg\|_{\cW^{q,r}}. 
\end{align}
Its first term is already bounded in \eqref{use4}, $\|\pd_{t}J_0\|_{\cW^{q,r}}\lec \mathcal B_\rho$.
For the second term in the RHS of \eqref{use6}, we move $\nabla'$ into mollification in $M_1,M_3,M_5, M_6$ terms, and move $\pd_s $  into mollification in $M_2,M_4$ terms.
As to the third term in the RHS,
notice that for $J_i$ with $i\leq d-1$, $M_1,M_3,M_5, M_6$ only contain  $E_{\ep_1,\ep_2,\ep_3}^{+}\left(\cdots\right)$,   
so we can move $\pd_d$ into mollification in $M_1,M_3,M_5, M_6$, while we move $\pd_s$ into mollification in $M_2,M_4$ terms.
Similar to the estimates in \eqref{M1}-\eqref{M6}, and using \eqref{eq5.15}, we get
\begin{align}\label{usehigh}
\text{RHS of \eqref{use6}} 
&\lec {\mathcal B}_\rho +\|\nabla\gamma\|_{L^{\infty}\left(B_{\rho}'\right)}\cdot  \|\nabla\pd_{t}E_{\ep_3}\left(\bm{v}\right),\nabla'\nabla E_{\ep_3}\left(\pi\right)\|_{\cH^{q,r}}. 
\end{align}
Taking $\ep_1,\ep_2\to 0$ in
\eqref{use6}--\eqref{usehigh}, we achieve
\begin{equation}\label{eq5.50}
 \| \nabla\pd_{t} E_{\ep_3}(\bm{v})\|_{\cH^{q,r}} \lec {\mathcal B}_\rho +\|\nabla\gamma\|_{L^{\infty}\left(B_{\rho}'\right)}\cdot  \|\nabla'\nabla  E_{\ep_3}\left(\pi\right)\|_{\cH^{q,r}}. 
\end{equation}

Combining \eqref{eq0929} and \eqref{eq5.50}, we get
\begin{equation}
\quad\,\|\nabla\pd_{t} E_{\ep_3}(\bm{v})\|_{\cH^{q,r}}+\|\nabla^3  E_{\ep_3}(\bm{v})\|_{\cH^{q,r}}+\|\nabla^2 E_{\ep_3}(\pi)\|_{\cH^{q,r}}
\lec {\mathcal B}_\rho.
\end{equation}
Taking $\ep_3 \to 0$, we get
$\nabla\pd_{t} \bm{v}, \ \nabla^3 \bm{v},  \ \nabla^2 \pi \in \cH^{q,r}$,
and
\begin{equation}
\quad\,\|\nabla\pd_{t} \bm{v}\|_{\cH^{q,r}}+\|\nabla^3  \bm{v}\|_{\cH^{q,r}}+\|\nabla^2 \pi\|_{\cH^{q,r}}\lec {\mathcal B}_\rho. 
\end{equation}
This completes the proof of Lemma \ref{lem52}. 
\end{proof}

\section{Auxiliary lemmas}\label{sec6}

We give several useful lemmas in the following.
\begin{lem}\label{lem-C-Z est}
\sl{ Suppose that $1<q,r<\infty$, $0<R,T\lec1$ and $q_*$ defined in \eqref{q_*}.
\begin{enumerate}
\item[(\rm i)]
For any $f\in L^q(\R^d)$, we have
\begin{equation*}
\|\nabla^2(-\Delta)^{-1}f\|_{L^q(\R^d)}\lec \|f\|_{L^q(\R^d)}.
\end{equation*}

\item[(\rm  ii)]
For any $f\in L^q_c(B_R)$, we have
\begin{equation*}
\|\nabla(-\Delta)^{-1}f\|_{L^q(B_R)}\lec R\|f\|_{L^q(B_R)}.
\end{equation*}

\item[(\rm  iii)]
For any $f\in L^r((0,T);L^q(\R^d))$, we have
\begin{equation*}
\|(\nabla^2,\pd_t)\int_{0}^{t} e^{(t-s)\Delta} f\ds\|_{L^r((0,T);L^q(\R^d))}\lec \|f\|_{L^r((0,T);L^q(\R^d))}.
\end{equation*}

\item[(\rm  iv)]
For any $f\in L^r((0,T);L^q_c(B_R))$, we have
\begin{equation*}
\|\int_{0}^{t}\nabla e^{(t-s)\Delta} f\ds\|_{L^r((0,T);L^q(B_R))}\lec R\|f\|_{L^r((0,T);L^q(B_R))},
\end{equation*}
\begin{equation*}
\|\int_{0}^{t}e^{(t-s)\Delta} f\ds\|_{L^r((0,T);L^q(B_R))}\lec R^2\|f\|_{L^r((0,T);L^q(B_R))}.
\end{equation*}
\begin{equation*}
\|\int_{0}^{t}\nabla e^{(t-s)\Delta} f\ds\|_{L^r((0,T);L^q(B_R))}\lec \|f\|_{L^r((0,T);L^{q_*}(B_R))},
\end{equation*}

\item[(\rm  v)]
For any $f\in L^{r_1}((0,T);L^{q_1}_c(B_R))$, we have
\begin{equation*}
\|\int_{0}^{t}\nabla e^{(t-s)\Delta} f\ds\|_{L^{r_2}((0,T);L^{q_2}(B_R))}\lec  \|f\|_{L^{r_1}((0,T);L^{q_1}(B_R))},
\end{equation*}
and
\begin{equation*}
\|\int_{0}^{t}e^{(t-s)\Delta} f\ds\|_{L^{r_2}((0,T);L^{q_2}(B_R))}\lec R  \|f\|_{L^{r_1}((0,T);L^{q_1}(B_R))},
\end{equation*}
where $1\leq q_1 \leq q_2$ and $1\leq r_1\leq r_2$ satisfy $\frac{1}{q_2}+1=\frac{1}{q_1}+\frac{1}{q_0}$, $\frac{1}{r_2}+1=\frac{1}{r_1}+\frac{1}{r_0}$ for some $q_0,r_0\geq 1$ satisfying $\frac{d}{q_0}+\frac{2}{r_0}>d+1$.

\item[(\rm  vi)] For any $f\in L^r((0,T);L^q_c(B_R))$, we have
\begin{equation*}
\|\int_{0}^{t}\nabla^3 (-\Delta)^{-1}e^{(t-s)\Delta} f\ds\|_{L^r((0,T);L^q(B_R))}\lec R\|f\|_{L^r((0,T);L^q(B_R))}.
\end{equation*}
\end{enumerate}
}
\end{lem}

\begin{proof}[Proof (sketch)]
(i) and (iii) are well-known Calder\'on-Zygmund estimates. The proof of (ii), (iv), (v), and (vi) relies on the estimates of Laplace, heat and Oseen kernel and Young's convolution inequality. 
See e.g.~\cite[Proposition 11.1 on page 107]{Lem02} for Oseen kernel.
When $q>\frac{d}{d-1}$, third estimate of (iv) is by (iv)$_1$, (iii), and Sobolev embedding.
\end{proof}

\begin{lem}[Lemma 3.1 in Chapter V  of \cite{Gia83}, page 161]\label{lemA2}
\sl{
Let $f(t)$ be  a non-negative bounded function defined in $[r_0,r_1]$, $r_0\geq 0$. Suppose that for $r_0\leq t<s\leq r_1$, we have
\begin{equation}
f(t)\leq [A(s-t)^{-\alpha}+B]+\theta f(s),
\end{equation}
where $A,B,\alpha,\theta$ are non-negative constants with $0\leq \theta<1$. Then for all $r_0\leq \rho<R\leq r_1$, we have
\begin{equation}
f(\rho)\leq c\, [A(R-\rho)^{-\alpha}+B],
\end{equation}
where $c$ is a constant depending on $\alpha$ and $\theta$. 
}
\end{lem}
\begin{lem}[Interior estimates for Stokes system]\label{lemA3}
\sl{
Let $0<R<1$, $1<q,r<\infty$, $q_*$ satisfy  \eqref{q_*}, and $\bm{f}\in L^{q_{*},r}(B_{R}(x_0)\times(-1,0))$, $\CF\in L^{q,r}(B_{R}(x_0)\times(-1,0))$ with $B_{R}(x_0)\subset \Omega$. Assume that $\bm{u}\in L^{q,r}\left(B_{R}(x_0)\times(-1,0)\right)$ is a very weak solution of the Stokes system \eqref{Stokes-Eqn}. Then $\nabla\bm{u}\in L^{q,r}\left(B_{R/2}(x_0)\times(-\frac34,0)\right)$ and
\begin{align}\label{interior1}
&\quad\ \|\nabla\bm{u}\|_{L^{q,r}\left(B_{R/2}(x_0)\times(-\frac34,0)\right)}\notag\\[6pt]
&\lec R^{-1}\|\bm{u}\|_{L^{q,r}\left(B_{R}(x_0)\times(-1,0)\right)}+\|\bm{f}\|_{L^{q_*,r}\left(B_{R}(x_0)\times(-1,0)\right)}+\|\CF\|_{L^{q,r}\left(B_{R}(x_0)\times(-1,0)\right)}.
\end{align}
If in addition $\bm{f}+\div \CF\in L^{r}\left((-1,0);W^{m,q}\left(B_{R}(x_0)\right)\right)$ for some integer $m\geq0$, then $\nabla^{m+2}\bm{u}\in L^{q,r}\left(B_{R/2}(x_0)\times(-\frac34,0)\right)$ and
\begin{align}\label{interior2}
&\quad\ \|\nabla^{m+2}\bm{u}\|_{L^{q,r}\left(B_{R/2}(x_0)\times(-\frac34,0)\right)}\notag\\[6pt]
&\lec R^{-m-2}\|\bm{u}\|_{L^{q,r}\left(B_{R}(x_0)\times(-1,0)\right)}+\sum_{i=0}^{m}R^{-m+i}\|\nabla^i(\bm{f}+\div \CF)\|_{L^{q,r}\left(B_{R}(x_0)\times(-1,0)\right)}.
\end{align}
If in addition $(\bm{u},p)$ is a very weak solution pair of the Stokes system \eqref{Stokes-Eqn} with $p\in L^{q,r}\left(B_{R}(x_0)\times(-1,0)\right)$, then $\pd_t\nabla^{m}\bm{u},\nabla^{m+2}\bm{u}, \nabla^{m+1}p \in L^{q,r}\left(B_{R/2}(x_0)\times(-\frac34,0)\right)$ and
\begin{align}\label{interior3}
&\|\pd_t\nabla^{m}\bm{u},\,\nabla^{m+2}\bm{u},\, \nabla^{m+1}p  \|_{L^{q,r}\left(B_{R/2}(x_0)\times(-\frac34,0)\right)}\lec R^{-m-2}\|\bm{u}\|_{L^{q,r}\left(B_{R}(x_0)\times(-1,0)\right)}\notag\\[6pt]
&\,\quad +R^{-m-1}\|p\|_{L^{q,r}\left(B_{R}(x_0)\times(-1,0)\right)}+\sum_{i=0}^{m}R^{-m+i}\|\nabla^i(\bm{f}+\div \CF)\|_{L^{q,r}\left(B_{R}(x_0)\times(-1,0)\right)}.
\end{align}
The constants depend on $q,r,m$ and is independent of $x_0$.
}
\end{lem}

\begin{proof}[Proof (sketch)]
\eqref{interior1} and \eqref{interior2} follow from the same proof of Lemma A.2 in \cite{CSYT08}. \eqref{interior3} is achieved by standard elliptic estimate of the Laplace equation of pressure, and standard parabolic estimate of the heat equation of velocity.
\end{proof}

\section*{Acknowledgments}
Chen was supported in part by National Natural Science Foundation of China under grant [12101556] and Zhejiang Provincial Natural Science Foundation of China under  grant [LY24A010015] and Fundamental Research Funds of ZUST [2025QN062]. 
The research of both Liang and Tsai was partially supported by Natural Sciences and Engineering Research Council of Canada (NSERC) under grant RGPIN-2023-04534.

\bibliography{NavierBC3}
\bibliographystyle{abbrv}

\end{document}